\documentclass[reqno]{amsart}

\usepackage{amssymb,amsthm}
\usepackage{mathrsfs}
\usepackage{longtable}

\usepackage{lunadiagrams}

\newtheorem{theorem}{Theorem}[subsection]
\newtheorem{lemma}{Lemma}[subsection]
\newtheorem{proposition}{Proposition}[subsection]
\newtheorem{corollary}{Corollary}[subsection]
\newtheorem{conjecture}{Conjecture}[subsection]

\theoremstyle{definition}
\newtheorem{definition}{Definition}[subsection]
\newtheorem{example}{Example}[subsection]

\numberwithin{equation}{section}
\numberwithin{table}{subsection}

\newcommand{\weights}{\mathcal X}
\newcommand{\supp}{\mathrm{supp}}

\newcommand{\defect}{\mathrm{d}}

\newcommand{\Hom}{\mathrm{Hom}}
\newcommand{\rank}{\mathrm{rank}}
\newcommand{\C}{\mathbb C}
\newcommand{\Z}{\mathbb Z}
\newcommand{\IN}{\mathbb N}
\newcommand{\Q}{\mathbb Q}
\newcommand{\A}{\mathbf A}
\newcommand{\p}{\mathbb P}
\newcommand{\s}{\mathscr S}
\renewcommand{\t}{\mathscr T}

\newcommand{\card}{\mathrm{card}}
\newcommand{\Lie}{\mathrm{Lie}}

\newcommand{\N}{\mathrm N}
\newcommand{\GL}{\mathrm{GL}}
\newcommand{\SL}{\mathrm{SL}}
\newcommand{\SO}{\mathrm{SO}}
\newcommand{\Sp}{\mathrm{Sp}}

\newcommand{\Ad}{\mathrm{Ad}}
\newcommand{\diag}{\mathrm{diag}}

\newcommand{\loccit}{{\em loc.cit.}}

\title{Wonderful subgroups of reductive groups\\ and spherical systems}
\author{P.~Bravi and G.~Pezzini}

\address{Dipartimento di Matematica G.~Castelnuovo, Sapienza Universit\`a di Roma\\ P.le Aldo Moro 5, 00185 Roma, Italy}
\email{bravi@mat.uniroma1.it}
\address{Department Mathematik, Emmy-Noether-Zentrum, FAU Erlangen-N\"urnberg\\ Cauerstr.~11, 91058 Erlangen, Germany}
\email{pezzini@mi.uni-erlangen.de}

\begin{document}

\begin{abstract}
Let $G$ be a semisimple complex algebraic group, and $H\subseteq G$ a wonderful subgroup. 
We prove several results relating the subgroup $H$ to the properties of a combinatorial invariant $\s$ of $G/H$, 
called its spherical system. 
It is also possible to consider a spherical system $\s$ as a datum defined by purely combinatorial axioms, 
and under certain circumstances our results prove the existence of a wonderful subgroup $H$ associated to $\s$. 
As a byproduct, we reduce for any group $G$ the proof of the classification of wonderful $G$-varieties, 
known as the Luna conjecture, to its verification on a small family of cases, called primitive.
\end{abstract}

\maketitle

\section{Introduction}

Luna's theory of wonderful varieties establishes a connection between {\em wonderful $G$-varieties} 
for a semisimple complex algebraic group $G$, 
and certain combinatorial objects called {\em spherical $G$-systems} (\cite{Lu01}). 
The main conjecture, i.e.\ that spherical systems classify these varieties, 
has previously been verified in many particular cases (see \cite{Lu01,Pe04,BP05,B07,BC10}). 
Here we develop some further general techniques related to this correspondence; 
this allows us to complete the proof in full generality in \cite{BP11p2}. 
It follows from Losev's results in \cite{Lo07} that
a wonderful $G$-variety having a given spherical $G$-system is uniquely determined 
(if it exists) up to $G$-isomorphism, 
therefore it remains to prove that every spherical system is {\em geometrically realizable}, 
i.e.\ that there actually exists a corresponding wonderful variety. 
A proof of the Luna conjecture has also been proposed by Cupit-Foutou in \cite{C09}, 
but her approach is different from ours: it is intrinsic (via invariant Hilbert schemes and equivariant deformation theory) 
but is not constructive, that is, does not explicitly provide the wonderful variety associated with a given spherical system.

Beyond the classification, the theory describes many geometric constructions and their combinatorial counterparts. For example, morphisms between wonderful varieties correspond to {\em quotients} of spherical systems. 

In this article we provide some results relating the stabilizer $H$ (called a {\em wonderful subgroup} of $G$) of a point in the open $G$-orbit of a wonderful $G$-variety $X$, to the spherical system $\s$ of $X$.

We propose two series of such results. 
The first is described in \S\ref{s:fromto}, and extracts information on $H$ only from combinatorial properties of $\s$. 
This approach works in great generality, but is yet unable to provide a complete description of $H$. 
More precisely, our results can be applied even if $\s$ is not assumed to be geometrically realizable, 
but in this case they do not provide a candidate for $H$.
A Levi part of $H$ (as already pointed out in \cite{Lo07}), 
can be determined using the classification of reductive spherical subgroups of $G$ 
and the knowledge of their spherical $G$-systems (\cite{BP11p3}). 
The unipotent radical of $H$ is more elusive, we only outline a possible strategy to determine it. 
%Our study of $H$ is also in some sense non-canonical, 
%since it requires the choice of a parabolic subgroup of $G$ minimal containing $H$: 
%examples show that such a choice is not unique (see \cite{}).

The second series of results is developed in \S\S\ref{s:fiberprod}--\ref{s:tails}. Here we analyze the relationship between $\s$, certain of its quotients, and the relative wonderful subgroups. We prove that under certain circumstances the geometric realizability of $\s$ is a consequence of the geometric realizability of these quotients, and we produce at the same time a precise description for $H$.

We also review some results from \cite{Lu01}, refining in particular those on fiber products of wonderful varieties. We discuss the so-called projective fibrations (see \cite[\S3.6]{Lu01}) without restrictions on the acting group $G$ as a particular case of a more general situation: {\em quotients of higher defect}.

The above mentioned constructions lead to the definition of {\em primitive} spherical systems, a direct generalization of the notion defined in \cite{Lu01} and subsequent works. The wonderful varieties corresponding to primitive spherical systems, called again {\em primitive}, play a significative role in the theory. In an analogy with representation theory they would correspond to indecomposable representations. In addition, the smaller class of primitive varieties admitting no quotient of higher defect would correspond to simple representations.

Primitive spherical systems can be used to provide a proof of the existence part of Luna's conjecture, in the same spirit of \cite{Lu01}. Namely, we prove that the conjecture follows from a case-by-case analysis on the primitive spherical systems without quotients of higher defect. This case-by-case analysis is the subject of \cite{BP11p2}.

\subsubsection*{Acknowledgments}
We are grateful to Prof.~D.~Luna for communicating to us his ideas which strongly influenced this paper and for pointing out some errors in a previous version, and to Prof.~F.~Knop for useful discussions. We thank the anonymous referee for his/her suggestions and corrections. The second-named author was supported by the DFG Schwerpunktprogramm 1388 -- Darstellungstheorie.

\section{Definitions}

Throughout this work, $G$ is a semisimple complex algebraic group, $T$ a maximal torus of $G$, $B$ a Borel subgroup containing $T$, $S$ the corresponding set of simple roots and $B_-$ the Borel subgroup opposite to $B$ with respect to $T$. We identify $B$-weights with $T$-weights. Unless otherwise stated, we number the simple roots of a connected component of the Dynkin diagram of $G$ as in \cite{Bou}.

For all $S'\subset S$, $P_{S'}$ denotes the corresponding parabolic subgroup containing $B$, $L_{S'}$ the Levi subgroup of $P_{S'}$ containing $T$, and $P_{-S'}$ the opposite parabolic subgroup of $P_{S'}$ with respect to $T$. We also denote by $G_{S'}$ the semisimple part of $L_{S'}$. For any parabolic subgroup $P$ containing $T$, $L_P$ is the standard Levi subgroup containing $T$.

In general, if $H$ is any affine algebraic group, 
$\Lie\,H$ denotes its Lie algebra and $H^\circ$ its connected component containing the unit element $e\in H$, 
$H^r$ is its radical and $H^u$ its unipotent radical, $Z(H)$ is its center, 
and if $K\supseteq H$ is another affine algebraic group, then $\N_K H$ is the normalizer of $H$ in $K$. 
We denote by $\mathcal X(H)$ the set of characters of $H$, and if $K\subseteq H$ is a subgroup, 
then $\mathcal X(H)^K$ denotes the set of characters of $H$ that are trivial on $K$. 
For any $H$-module $V$ and any $\chi\in\mathcal X(H)$, we define
\[
V^{(H)}_\chi = \{ v\in V \;|\; hv = \chi(h)v \;\forall h\in H \},
\]
and
\[
V^{(H)} = \bigcup_{\chi\in\mathcal X(H)} V_\chi^{(H)}.
\]

\subsection{Wonderful varieties}
We recall here and in the next sections the basic definitions about wonderful varieties, see \cite{Lu01} and \cite{BL11} for further details, proofs and references.

\begin{definition}
A $G$-variety is called {\em wonderful of rank $r$} if it is smooth, complete, with $r$ smooth prime $G$-divisors $D_1,\ldots,D_r$ with normal crossings, such that the $G$-orbit closures are exactly all the intersections $\cap_{i\in I}D_i$, for any $I\subseteq\{1,\ldots,r\}$.
\end{definition}

We will always assume that $Z(G)$ acts trivially on a wonderful $G$-variety. This is justified by the results of Luna in \cite{Lu01}, where the classification of wonderful varieties is reduced to the case where $G$ is adjoint.

A wonderful $G$-variety is known to be projective and spherical. Recall that a normal $G$-variety is called spherical if it contains an open $B$-orbit. The non-$G$-stable prime $B$-divisors of a spherical $G$-variety $X$ are called colors, their set is denoted by $\Delta_X$. A simple (i.e.\ with a unique closed $G$-orbit) spherical $G$-variety is wonderful if and only if it is smooth complete and toroidal (i.e.\ no color contains a $G$-orbit).

Let $X$ be a wonderful $G$-variety of rank $r$. 
There exists a unique point $z\in X$ stabilized by $B_-$. 
The orbit $G.z$ is the unique closed $G$-orbit, equal to $\cap_{i=1}^rD_i$; 
the parabolic subgroup opposite to $G_z$, usually denoted by $P_X$, is the stabilizer of the open $B$-orbit of $X$. 
Let $S^p_X$ denote the set of simple roots such that $P_X=P_{S^p_X}$.

The $T$-weights of the normal space $\mathrm T_zX/\mathrm T_zG.z$ of $G.z$ at $z$ are called spherical roots
%\footnote{For an alternative definition of spherical roots, valid for any spherical variety, see \cite{Lo07}.}
, their set is denoted by $\Sigma_X$. 
Spherical roots are naturally in bijective correspondence with prime $G$-divisors, in the following way: 
we denote by $D^\sigma$ the prime $G$-divisor such that $\sigma$ is the $T$-weight of $\mathrm T_zX/\mathrm T_zD^\sigma$, 
for all $\sigma\in\Sigma_X$.

In general, we define the rank of a spherical variety $Y$ to be the rank of the lattice $\Lambda(Y)$ of the $B$-weights of $B$-semiinvariants in $\C(Y)$. If $X$ is wonderful the spherical roots form a basis of $\Lambda(X)$  (see e.g.\ \cite[\S1]{Kn96}), therefore the rank of $X$ as wonderful variety and as spherical variety coincide.

The colors of $X$ are representatives of a basis of $\mathrm{Pic}X$. One can thus define a $\Z$-bilinear pairing, called Cartan pairing, $c_X\colon \Z\Delta_X\times\Z\Sigma_X\to\Z$ such that $[D^\sigma]=\sum_{D\in\Delta_X}c_X(D,\sigma)[D]$ in $\mathrm{Pic}X$, for all $\sigma\in\Sigma_X$.

Sometimes it is notationally convenient to use the following map: $\rho_X\colon\Delta_X\to\Hom_\Z(\Z\Sigma_X,\Q)$, where $\rho_X(D)(-)=c_X(D,-)$. Then, the image of $\rho_X$ spans the vector space $\Hom_\Z(\Z\Sigma_X,\Q)$.

For all $\alpha\in S$, define the set $\Delta_X(\alpha)=\{D\in\Delta_X:P_{\{\alpha\}}D\neq D\}$ of colors {\em moved} by $\alpha$. One has $\cup_{\alpha\in S}\Delta_X(\alpha)=\Delta_X$ and, for all $\alpha\in S$, $\card(\Delta_X(\alpha))\leq2$. Clearly, $\alpha\in S^p_X$ if and only if $\Delta_X(\alpha)=\emptyset$. Moreover, if $D\in\Delta_X(\alpha)$ and $\card(\Delta_X(\alpha))=1$ then $\rho_X(D)$ is uniquely determined by $\alpha$. One has $\card(\Delta_X(\alpha))=2$ if and only if $\alpha\in S\cap\Sigma_X$, in this case, say $\Delta_X(\alpha)=\{D^+,D^-\}$, the elements $\rho_X(D^+)$ and $\rho_X(D^-)$ are not always determined by $\alpha$, but their sum is. Let $\A_X$ denote the set of colors $D\in\cup_{\alpha\in S\cap\Sigma_X}\Delta_X(\alpha)$ endowed with the $\Z$-linear functionals $\rho_X(D)$.

\begin{definition}
The datum of $(S^p_X,\Sigma_X,\A_X)$, also denoted by $\s_X$, is called the {\em spherical system} of $X$.
\end{definition} 

%If $X$ has generic stabilizer $H$, then the normalizer $\N_GH$ is wonderful too. Its spherical system will be denoted by $\N_G(\s_X)$, and differs from $\s$ only for some spherical roots, which are replaced by their doubles. See \cite{Lo07} or \cite{BL11} for details.

Any $G$-orbit closure $X'$ of $X$ is a wonderful $G$-variety itself. Its set of spherical roots $\Sigma'=\Sigma_{X'}$ is included in $\Sigma_X$, and $X'$ is called the localization of $X$ with respect to $\Sigma'$ and denoted by $X_{\Sigma'}$. Indeed, $X_{\Sigma'}=\cap_{\sigma\in\Sigma_X\setminus\Sigma'}D^\sigma$. The spherical system $\s_{X'}$ of $X'$ is given by $S^p_X$, $\Sigma'$ and $\A_{X'}$, where the latter can be identified with the set of colors $D\in\cup_{\alpha\in S\cap\Sigma'}\Delta_X(\alpha)$. One has $c_{X'}(D,\sigma)=c_X(D,\sigma)$ for all $\sigma\in\Sigma'$.

In particular, any spherical root of $X$ is the spherical root of a wonderful $G$-variety of rank 1. 
The wonderful $G$-varieties of rank 1 are well known, for all $G$ (\cite{A83}). 
The finite set of spherical roots of wonderful $G$-varieties of rank 1 is denoted by $\Sigma(G)$. 
Any spherical root $\sigma$ is a linear combination of simple roots with non-negative integer coefficients: 
we define $\supp\,\sigma$ the set of simple roots whose coefficient is non-zero. 
For the spherical roots we follow the labeling of \cite[Definition~1.1.1]{B09}. 
If $\Sigma'\subset\Sigma(G)$ then $\supp\,\Sigma'$ is the union of $\supp\,\sigma$ for all $\sigma\in\Sigma'$.

The spherical system of $X$ is determined by the spherical systems of all the localizations $X_{\Sigma'}$ of rank 2 (actually, it is enough to restrict to the localizations of rank 1 and those of rank 2 with $\Sigma'\cap S\neq\emptyset$). Furthermore, the wonderful $G$-varieties of rank 2 are known, for all $G$ (\cite{W96}).

The stabilizer $H$ of a point of the open $G$-orbit of a wonderful $G$-variety $X$ 
is called a {\em generic stabilizer} of $X$, 
in this case $X$ is also called a {\em wonderful compactification} or a {\em wonderful embedding} of $G/H$. 
A subgroup $H$ of $G$ is called {\em wonderful} if $G/H$ admits a wonderful compactification $X$, 
in this case the spherical system $\s$ of $X$ is also called the spherical system of $H$ and, vice versa, 
$H$ is called {\em the} wonderful subgroup of $G$ associated with $\s$ 
(it is uniquely determined only up to $G$-conjugation). 

A subgroup $H$ of $G$ is called {\em spherical} if $G/H$ is a spherical $G$-variety. In particular, every wonderful subgroup of $G$ is spherical.

\subsection{Spherical systems}
\begin{definition} \label{def:sphsys}
Let $(S^p,\Sigma,\mathbf A)$ be a triple such that $S^p \subset S$, $\Sigma \subset \Sigma(G)$ 
and $\A$ is a finite set endowed with a pairing $c\colon \Z\A\times\Z\Sigma\to\Z$. 
For every $\alpha \in \Sigma \cap S$, let $\A (\alpha)$ denote the set $\{D \in \A \colon c(D,\alpha)=1 \}$. 
Such a triple is called a \textit{spherical $G$-system} if: 
\begin{itemize} 
\item[(A1)] for every $D \in \A$ we have $c(D,-)\leq 1$, 
and if $c(D, \sigma)=1$ for some $\sigma\in\Sigma$ then $\sigma \in S\cap\Sigma$; 
\item[(A2)] for every $\alpha \in \Sigma \cap S$, $\A(\alpha)$ contains exactly two elements, 
and denoting by $D_\alpha^+$ and $D_\alpha^-$ these elements, 
it holds $c(D_\alpha^+,-) + c(D_\alpha^-,-) = \langle -,  \alpha^\vee \rangle$;
%\footnote{By $\langle-,-\rangle$ we denote the usual scalar product induced by the Killing form on the root space 
%and by $\alpha^\vee$ the co-root associated with the simple root $\alpha$.} 
\item[(A3)] the set $\A$ is the union of $\A(\alpha)$ for all $\alpha\in\Sigma \cap S$; 
\item[($\Sigma 1$)] if $2\alpha \in \Sigma \cap 2S$ 
then $\frac{1}{2}\langle \sigma, \alpha^\vee\rangle$ is a non-positive integer 
for every $\sigma \in \Sigma \setminus \{ 2\alpha \}$; 
\item[($\Sigma 2$)] if $\alpha, \beta \in S$ are orthogonal and $\alpha + \beta \in \Sigma$ 
then $\langle \sigma ,  \alpha ^\vee\rangle = \langle \sigma ,  \beta ^\vee\rangle$ for every $\sigma \in \Sigma$; 
\item[(S)] for every $\sigma \in \Sigma$, there exists a wonderful $G$-variety $X$ of rank 1 
with $S^p_X=S^p$ and $\Sigma_X=\{\sigma\}$.
\end{itemize}
\end{definition}
With abuse of language two spherical $G$-systems $(S^p_1,\Sigma_1,\mathbf A_1)$ and $(S^p_2,\Sigma_2,\mathbf A_2)$ are defined to be equal if $S^p_1=S^p_2$, $\Sigma_1=\Sigma_2$, and if there is a bijection $\mathbf A_1\to \mathbf A_2$ compatible with the respective Cartan pairings.

The set of rank $1$ wonderful $G$-varieties is finite and known for any $G$, 
and the last axiom admits an equivalent and more explicit combinatorial formulation, see  \cite[\S1.1.6]{BL11}.

The set $\Sigma$ (resp.\ its cardinality) is the set of {\em spherical roots} (resp.\ the {\em rank}) of the spherical $G$-system.

For all wonderful $G$-varieties $X$, the spherical system of $X$ is a spherical $G$-system in the sense of the above definition.

\subsection{The set of all colors}\label{s:allcolors}
Let $\s =(S^p,\Sigma,\A)$ be a spherical $G$-system, the set $\Delta \supseteq \A$ of its {\em colors} (endowed with an extension of $c$ to $\Z\Delta\times\Z\Sigma$) is the disjoint union of $\A$ and two finite disjoint sets $\Delta^{2a}$ and $\Delta^b$, defined as follows:
\begin{itemize}
\item[-] $\Delta^{2a}=\{D_\alpha \colon \alpha\in S\mbox{ and } 2\alpha\in\Sigma\}$ where $c(D_\alpha,-)=\frac12 \langle-, \alpha^\vee\rangle$, and $D_\alpha=D_\beta$ iff $\alpha=\beta$;
\item[-] $\Delta^b=\{D_\alpha \colon \alpha\in\left(S \setminus \left(\Sigma\cup\frac12\Sigma\cup S^p \right) \right)\}$
where $D_\alpha=D_\beta$ iff $\alpha=\beta$, or $\alpha, \beta$ are orthogonal and $\alpha+\beta\in\Sigma$; here $c(D_\alpha,-)= \langle-, \alpha^\vee\rangle$.
\end{itemize}

As for spherical systems of wonderful varieties, we will also denote the (extended) Cartan pairing $c\colon\Z\Delta\times\Z\Sigma\to\Z$ of $\s$ as a map $\rho\colon\Delta\to\Hom_\Z(\Z\Sigma,\Z)$.

The definition of $\Delta(\alpha)\subseteq\Delta$ for any $\alpha\in S$, called again the set of colors {\em moved} by $\alpha$, is the following:
\begin{itemize}
\item[-] if $\alpha\in S\cap\Sigma$, then $\Delta(\alpha) = \{D\in \A\;|\; c(D,\alpha)=1\}$,
\item[-] if $\alpha\in S^p$, then $\Delta(\alpha)=\emptyset$,
\item[-] otherwise, $\Delta(\alpha)=\{D_\alpha\}$, according to the above definitions of $\Delta^{2a}$ and $\Delta^b$.
\end{itemize}

For all wonderful $G$-varieties $X$, the set of colors of $X$ can be identified with the set of colors of $\s_X$ (considered as abstract sets) in such a way that the Cartan pairing of $X$ coincides with the pairing $c$ defined here, and $\Delta(\alpha)=\Delta_X(\alpha)$ for all simple roots $\alpha$.

Finally, we say that a color $D$ is {\em free} if there exists at most one spherical root $\sigma$ such that $c(D,\sigma)>0$. If $D\in\A$, it is free if and only if it is moved by only one simple root.

\subsection{Normalizers, spherical closure and augmentations}

\begin{definition}\cite[\S2.2]{BL11}\label{def:loose}
Let $\sigma$ be a spherical root of a spherical $G$-system $\s$. Then $\sigma$ is {\em loose} if one of the following holds:
\begin{itemize}
\item[($\mathsf A$)] $\sigma\in\Sigma\cap S$ and $\rho(D^+)=\rho(D^-)$ where $\Delta(\alpha)=\{D^+,D^-\}$;
\item[($\mathsf B$)] $\sigma=\alpha_1+\ldots+\alpha_m$ where $\{\alpha_1,\ldots,\alpha_m\}\subseteq S$ has type $\mathsf B_m$ with short root $\alpha_m$, and $\Delta(\alpha_i)=\varnothing$ for $i=2,\ldots,m$;
\item[($\mathsf G$)] $\sigma=2\alpha_1+\alpha_2$ where $\{\alpha_1,\alpha_2\}\subseteq S$ has type $\mathsf G_2$ with short root $\alpha_1$.
\end{itemize}
\end{definition}

\begin{lemma}\label{lemma:normalizer}
Let $\s=(S^p,\Sigma,\A)$ be a spherical $G$-system. Let $\Sigma^{N}$ be defined as $\Sigma$ where we replace each loose spherical root $\sigma$ with $2\sigma$, and let $\A^{N}$ be obtained from $\A$ by removing all pairs $\A(\alpha)$ for $\alpha\in\Sigma\cap S$ loose of type $\mathsf A$. Then $\N_G(\s) = (S^p, \Sigma^{N}, \A^{N})$ is a spherical $G$-system. If it is geometrically realizable, then $\s$ is geometrically realizable too.
\end{lemma}
\begin{proof}
The fact that $\N_G(\s)$ is a spherical system follows directly from the definition of $\Sigma^{N}$. Suppose that $\N_G(\s)$ is geometrically realizable, with generic stabilizer $K$: then \cite[Theorem 2]{Lo07} implies that $K$ is equal to its normalizer in $G$.

Let $c^{N}$ be the Cartan pairing of $\N_G(\s)$, and $\Delta^{N}$ its set of colors. It is easy to check that for each loose spherical root $\sigma$ of $\s$, we have $c^{N}(D,2\sigma)\in2\Z$ for all $D\in\Delta^{N}$. It follows that $c^{N}$ can be extended to a $\Z$-bilinear pairing $c^{N}\colon \Z\Delta^{N}\times\Z\Sigma\to\Z$. Then, we define the following map:
\[
\begin{array}{ccc}
\Z\Sigma & \to &  \Z^{\Delta^{N}} \\
\sigma & \mapsto & \sum_{D\in\Delta} c^{N}(D,\sigma) D.
\end{array}
\]
Thanks to \cite[Lemme 6.3.1, Lemme 6.3.3]{Lu01}, the image of this map corresponds to a spherical subgroup $H$ of $G$ such that $H$ is a normal subgroup of $K$ and $\Z\Sigma=\Lambda(G/H)$.

Consider the sets $V$ and $V^{N}$ of $G$-invariant, $\Q$-valued discrete valuations resp.\ of the fields $\C(G/H)$ and $\C(G/K)$. These sets are naturally identified with two finitely generated convex cones of resp.\ $\Hom_\Z(\Z\Sigma,\Q)$ and $\Hom_\Z(\Z\Sigma^{N},\Q)$.

Since the quotient $\Z\Sigma/\Z\Sigma^{N}$ is finite we have $K^{\circ}=H^{\circ}$, and we deduce from \cite[\S4]{Kn96} that the natural map $G/H\to G/K$ induces an identification of the vector spaces $\Hom_\Z(\Z\Sigma,\Q)$ and $\Hom_\Z(\Z\Sigma^{N},\Q)$ in such a way that $V=V^{N}$.

Moreover $V^{N}$ is the set of elements of $\Hom_\Z(\Z\Sigma^{N},\Q)$ that are non-negative on $\Sigma^{N}$ (see \cite[\S1.2]{Lu01}). Since $|\Sigma|=|\Sigma^{N}|$ and for all $\sigma\in\Sigma^{N}$ either $\sigma$ or $\sigma/2$ is in $\Sigma$, the set $V$ is the subset of elements of $\Hom_\Z(\Z\Sigma,\Q)$ that are non-negative on $\Sigma$.

Since $\Sigma$ is a basis of $\Z\Sigma=\Lambda(G/H)$ we conclude that $H$ is wonderful and that the wonderful compactification $X$ of $G/H$ has set of spherical roots $\Sigma$ (see \cite[\S1]{Kn96}). Moreover, since $H^{\circ}=K^{\circ}$ the groups $H$ and $K$ have the same normalizer, whence $K=\N_GH$. Then it follows from \cite[Theorem 2]{Lo07} that $\N_G(\s_X) = \N_G(\s)$, and finally from $\Sigma_X=\Sigma$ we deduce $\s=\s_X$.
\end{proof}

\begin{definition}
A spherical $G$-system $\s$ is {\em spherically closed} if it has no loose spherical roots of type $\mathsf B$ nor type $\mathsf G$.
\end{definition}

Let $H$ be a spherical subgroup of $G$. The group $\N_GH/H$ is naturally identified with the group of $G$-equivariant automorphisms of $G/H$, therefore it acts on the set of $B$-stable prime divisors of $G/H$. The kernel of this action is called the {\em spherical closure} of $H$, and is denoted by $\overline H$. The group $H$ is called {\em spherically closed} if $H=\overline H$.

\begin{lemma}{\cite[\S2.4.1]{BL11}}
Let $X$ be a wonderful variety with generic stabilizer $H$. Then $\s_X$ is spherically closed if and only if $H$ is spherically closed.
\end{lemma}

We recall that a spherically closed subgroup of $G$ is wonderful thanks to \cite[Corollary 7.6]{Kn96}.

The combinatorial notion of {\em augmentations} is used to describe all spherical subgroups of $G$ having the same spherical closure.

\begin{definition}{\cite[\S2.2]{Lu01}}
Let $\s=(S^{p}, \Sigma, \A)$ be a spherical $G$-system, consider a subgroup $\Xi'\subseteq\weights(T)$ containing $\Sigma$ and a map $\rho'\colon\A\to\Hom_\Z(\Xi',\Z)$. The couple $(\Xi',\rho')$ is an {\em augmentation} of $\s$ if:
\begin{itemize}
\item[(a1)] the map $\rho'$ coincides on $\Z\Sigma$ with the map $\rho$ of $\s$;
\item[(a2)] for all $\alpha\in\Sigma\cap S$ we have $\langle\rho'(\delta_\alpha^{+}),\gamma\rangle + \langle\rho'(\delta_\alpha^{+}),\gamma\rangle = \langle\gamma,\alpha^{\vee}\rangle$ for all $\gamma\in\Xi'$, where $\Delta(\alpha) = \{\delta_\alpha^{+}, \delta_\alpha^{-}\}$;
\item[($\sigma$1)] if $2\alpha\in\Sigma\cap 2S$, then $\alpha\notin\Xi'$ and $\langle\gamma,\alpha^{\vee}\rangle$ is an even integer for all $\gamma\in\Xi'$;
\item[($\sigma$2)] if $\alpha+\beta\in\Sigma$ or $\frac12(\alpha+\beta)\in\Sigma$, where $\alpha$ and $\beta$ are orthogonal simple roots, then $\alpha^{\vee}$ and $\beta^{\vee}$ coincide on $\Xi'$;
\item[(s)] for all $\alpha\in S^{p}$ the coroot $\alpha^{\vee}$ vanishes on $\Xi'$.
\end{itemize}
\end{definition}

Let $H$ be a spherical subgroup of $G$ and denote by $X$ the wonderful compactification of $G/\overline H$. The set of colors $X$ is naturally identified with the set of $B$-stable prime divisors of $G/H$.

Define the map $\rho_H\colon\Delta_X\to \Hom_\Z(\Lambda(G/H),\Z)$ as follows: $\langle\rho_H(D),\gamma\rangle$ is the evaluation of a $B$-eigenvector $f_\gamma\in\C(G/H)^{(B)}_\gamma$ along the $B$-stable prime divisor $D$ of $G/H$.

\begin{proposition}{\cite[Proposition 6.4]{Lu01}}\label{prop:augm}
The map $H\mapsto (\Lambda(G/H),\rho_H)$ is a bijection between the set of spherical subgroups of $G$ having the same spherical closure $\overline H$ and the set of all augmentations of the spherical $G$-system $\s_X$, where $X$ is the wonderful compactification of $G/\overline H$.
\end{proposition}

\subsection{Quotients}
As a general reference for this section see \cite[\S3.3]{Lu01}.

Let $f\colon X\to Y$ be a surjective $G$-equivariant morphism between two spherical $G$-varieties, and define $\Delta_f=\{ D\in\Delta_X \;|\; f(D)=Y \}$. 

Let $X$ be a wonderful $G$-variety. A subset $\Delta'$ of $\Delta_X$ is called {\em distinguished} if there exists a linear combination $D$ of elements of $\Delta'$ with positive coefficients such that $c_X(D,\sigma)\geq 0$ for all $\sigma\in\Sigma_X$.

We say that a distinguished subset $\Delta'$ of $\Delta_X$ is {\em good}, 
if the monoid $\{ \sigma\in\IN\Sigma_X \colon c_X(D,\sigma)=0 \;\forall D\in\Delta'\}$ is free; 
we denote in this case its indecomposable elements by $\Sigma_X/\Delta'$.
\begin{proposition}\label{prop:quot}
For every good distinguished subset $\Delta'$ of $\Delta_X$ there exists a couple $(f,Y)$, 
where $Y$ is a wonderful $G$-variety 
and $f\colon X\to Y$ is a $G$-equivariant surjective map with connected fibers
%\footnote{A $G$-equivariant surjective morphism with connected fibers between wonderful varieties is also called wonderful, 
%see \cite[\S2.3]{BL11}.}
, such that $\Delta'=\Delta_f$.
For two such couples $(f,Y)$ and $(f',Y')$, $\Delta_f=\Delta_{f'}$ if and only if 
there exists a $G$-equivariant isomorphisms $\phi\colon Y\to Y'$ such that $f'=\phi\circ f$.
\end{proposition}
%In principle, a distinguished but not (*)-distinguished set of colors corresponds to a $G$-equivariant morphism $X\to Y$ too. In this case $Y$ satisfies the definition of a wonderful variety, except for the smoothness condition on $Y$ and its prime $G$-divisors. Such a variety is also called the {\em canonical embedding} of its open $G$-orbit.

%A homogeneous space $G/H$ has at most one wonderful embedding and at most one canonical embedding (equal, it they both exist). The above proposition can be also stated in terms of the open orbit $G/H$ of $X$, or equivalently in terms of its generic stabilizer $H$.
We will say that an inclusion $K\supseteq H$ between subgroups of $G$ is {\em co-connected} 
if the quotient $K/H$ is connected. We can rephrase Proposition~\ref{prop:quot} as follows.
\begin{proposition}\label{prop:quotH}
Let $X$ be a wonderful $G$-variety with generic stabilizer $H$. For any subgroup $K$ of $G$ containing $H$, we define $\Delta_K=\Delta_f$ where $f\colon G/H\to G/K$ is the natural projection. Then $K\mapsto \Delta_K$ is a bijection between the co-connected inclusions $K\supseteq H$ such that $K$ is of finite index in $\N_G K$ and the set of distinguished subsets of $\Delta_X$. The homogeneous space $G/K$ admits a wonderful compactification if and only if $\Delta_K$ is good.
\end{proposition}
A good distinguished subset $\Delta'\subseteq \Delta_X$ is 
\begin{itemize}
\item[-] {\em smooth} if $\Sigma_X/\Delta'\subseteq \Sigma_X$;
\item[-] {\em homogeneous} (or {\em parabolic}) if $\Sigma_X/\Delta'=\emptyset$.
\end{itemize}
\begin{proposition}
A $G$-equivariant, surjective map with connected fibers between wonderful varieties $f:X\to Y$ is smooth if and only if $\Delta_f$ is smooth. Moreover, $Y$ is homogeneous if and only if $\Delta_f$ is homogeneous.
\end{proposition}
Obviously, the above definition of distinguished set of colors can also be given for any spherical system (with no assumptions on its geometric realizability).
\begin{definition}
Let $\s=(S^p,\Sigma,\A)$ be a spherical $G$-system with set of colors $\Delta$, and let $\Delta'\subseteq\Delta$ be a distinguished set of colors. Set:
\begin{itemize}
\item[-] $S^p/\Delta'= \{ \alpha\in S \;|\; \Delta(\alpha)\subseteq \Delta' \}$;
\item[-] $\Sigma/\Delta' = $ the indecomposable elements (or equivalently the minimal generators) of the monoid $\{ \sigma\in\IN\Sigma \;|\; c_X(D,\sigma) = 0 \; \forall D\in\Delta' \}$;
\item[-] $\A/\Delta' = \bigcup_{\alpha\in S\cap\Sigma/\Delta'} \A(\alpha)$.
\end{itemize}
We say that the distinguished set $\Delta'$ is {\em good} if and only if the monoid $\{ \sigma\in\IN\Sigma \;|\; c_X(D,\sigma) = 0 \; \forall D\in\Delta' \}$ is free and $\s/\Delta'=(S^p/\Delta',\Sigma/\Delta',\A/\Delta')$ is a spherical $G$-system. In this case, the bilinear pairing $c$ for $\s/\Delta'$ is induced from that of $\s$ in a natural way.\footnote{We point out that all distinguished sets of colors are good, but the only proof we know relies on a case-by-case checking (\cite{B09}). On the other hand we do not make use of that result in this paper.}
\end{definition}

If $f\colon X\to Y$ is a surjective $G$-equivariant morphism between two wonderful $G$-varieties as above, then $\s_Y=\s_X/\Delta_f$.

\subsection{Localizations}\label{s:loc}
Let $X$ be a wonderful $G$-variety and $S'$ a subset of $S$. The variety $X$ contains a well defined wonderful $L_{S'}$-subvariety $X_{S'}$, called the {\em localization} of $X$ in $S'$ (see \cite[\S3.2]{Lu01}). One has $S^p_{X_{S'}}=S^p_X \cap S'$. The set $\Sigma_{X_{S'}}$ is the set of spherical roots of $X$ whose support is contained in $S'$. The set $\A_{X_{S'}}$ is the union of $\Delta_X(\alpha)$ for all $\alpha\in S'\cap\Sigma_X$. The localization $\s_{S'}$ of a spherical $G$-system $\s$ in $S'$ is defined analogously.

Recall that we have also defined the localization $X_{\Sigma'}$ of $X$ on a set of spherical roots%
\footnote{In cases where $\Sigma'=S'\subseteq \Sigma\cap S$ our notation may be ambiguous: the context will then indicate whether we are localizing to a subset of simple roots or of spherical roots.}
$\Sigma'\subseteq \Sigma_X$. The same can be done for a spherical $G$-system $\s=(S^p,\Sigma,\A)$ and $\Sigma'\subseteq\Sigma$: we define $\s_{\Sigma'}=(S^p,\Sigma',\A_{\Sigma'})$ where
\[
\A_{\Sigma'} = \bigcup_{\alpha\in S\cap\Sigma'} \Delta(\alpha).
\]
The bilinear pairing of $\s_{\Sigma'}$ is the one of $\s$ restricted to $\Z\Sigma'$.

Given a subset $\widetilde\Delta$ of the colors $\Delta$ of $\s$, 
we can also define the {\em restriction} of $\widetilde\Delta$ to a localization $\s_{S'}$ (for $S'\subseteq S$) 
and to a localization $\s_{\Sigma'}$ (for $\Sigma'\subseteq \Sigma$). In the first case, we define
\[
\widetilde\Delta|_{S'}=\widetilde\Delta\cap\left(\bigcup_{\alpha\in S'}\Delta(\alpha)\right).
\]
It is immediate to show that if $\s=\s_X$ is geometrically realizable, then
\[
\widetilde\Delta|_{S'} = \left\{ D \cap X_{S'} \;|\; D\in\widetilde\Delta \right\} \subseteq \Delta_{X_{S'}}.
\]

In the second case, a more careful analysis is needed, since the full set of colors of a localization $\s_{\Sigma'}$ is not related to colors of $\s$ as easily as its subset $\A_{\Sigma'}$ is.

For notational convenience, let us denote by $\Delta^b_1\subseteq\Delta^b$ the set of colors moved by only one simple root, and by $\Delta^b_2$ the set $\Delta^b\setminus\Delta_1^b$. We also write $\Delta_{\Sigma'}$ for the whole set of colors of the spherical system $\s_{\Sigma'}$, and correspondingly $\Delta_{\Sigma'}(\alpha)$ denotes the colors of $\s_{\Sigma'}$ moved by $\alpha\in S$ (recall that $\Delta_{\Sigma'}(\alpha)$ is identified with $\Delta(\alpha)$ if $\alpha\in\Sigma'\cap S$).

Then we define the restriction of $\widetilde\Delta$ to $\s_{\Sigma'}$ as
\[
\widetilde\Delta|_{\Sigma'} = \left(\widetilde\Delta_{\Sigma',1}\right)\cup\left(\widetilde\Delta_{\Sigma',2}\right)\cup\left(\widetilde\Delta_{\Sigma',3}\right)\cup\left(\widetilde\Delta_{\Sigma',4}\right),
\]
where:
\[
\widetilde\Delta_{\Sigma',1} =  \left(\bigcup_{\alpha\in\Sigma'\cap S} \Delta(\alpha)\cap\widetilde\Delta \right) \cup \left(\bigcup_{\begin{array}{c}\scriptstyle{\alpha\in(\Sigma\setminus\Sigma')\cap S }\\  \scriptstyle{\textrm{with } \Delta(\alpha)\subseteq \widetilde\Delta}\end{array}} \Delta_{\Sigma'}(\alpha)\right),
\]
\[
\widetilde\Delta_{\Sigma',2} =  \bigcup_{\begin{array}{c}\scriptstyle{\alpha\in\frac12\Sigma\cap S }\\  \scriptstyle{\textrm{with } \Delta(\alpha)\subseteq \widetilde\Delta}\end{array}} \Delta_{\Sigma'}(\alpha),
\]
\[
\widetilde\Delta_{\Sigma',3} =  \bigcup_{\begin{array}{c}\scriptstyle{\{\alpha\}\in\Delta^b_1 }\\  \scriptstyle{\textrm{with } \Delta(\alpha)\subseteq \widetilde\Delta}\end{array}} \Delta_{\Sigma'}(\alpha),
\]
\[
\widetilde\Delta_{\Sigma',4} =  \bigcup_{\begin{array}{c}\scriptstyle{\{\alpha,\beta\}\in\Delta^b_2 }\\  \scriptstyle{\textrm{with } \Delta(\alpha)=\Delta(\beta)\subseteq \widetilde\Delta}\end{array}} \left(\Delta_{\Sigma'}(\alpha)\cup\Delta_{\Sigma'}(\beta)\right).
\]

With this definition we can state the following lemma and corollary.
\begin{lemma}\label{lem:restr}
Let $\widetilde\Delta$ be a good distinguished set of colors of $\s$. 
Then the restriction $\widetilde\Delta|_{\Sigma'}$ is a good distinguished set of colors of $\s_{\Sigma'}$, and we have
\begin{equation}\label{eqn:reseq}
(\s/\widetilde\Delta)_{\Sigma''} = (\s_{\Sigma'}) / (\widetilde\Delta|_{\Sigma'}),
\end{equation}
where $\Sigma''$ consists of the elements of $\Sigma/\widetilde\Delta$ that are linear combination of elements of $\Sigma'$.
\end{lemma}
\begin{proof}
Let us fix 
$E=\sum a_D\,D$, for $D\in\widetilde\Delta$, 
with positive coefficients such that $c(E,\sigma)\geq0$ for all $\sigma\in\Sigma$, and define
$F=\sum b_D\,D$, for $D\in\widetilde\Delta|_{\Sigma'}$,
as follows ($S_D$ is the set of simple roots moving $D\in\A$).
\[b_D=\left\{\begin{array}{cl}
\frac{|S_{D}\cap\Sigma'|}{|S_{D}|}a_D & \mbox{if }D\in\Delta_{\Sigma'}(\alpha)\mbox{ and }\alpha\in\Sigma'\cap S \\
\frac1{|S_{D^+_\alpha}|}a_{D^+_\alpha}+\frac1{|S_{D^-_\alpha}|}a_{D^-_\alpha} & \mbox{if }D\in\Delta_{\Sigma'}(\alpha)\mbox{ and }\alpha\in(\Sigma\setminus\Sigma')\cap S \\
a_D & \mbox{if }D\in\Delta_{\Sigma'}(\alpha)\mbox{ and }\alpha\in\frac12\Sigma'\cap S \\
\frac12a_D & \mbox{if }D\in\Delta_{\Sigma'}(\alpha)\mbox{ and }\alpha\in\frac12(\Sigma\setminus\Sigma')\cap S \\
a_D & \mbox{if }D\in\Delta_{\Sigma'}(\alpha)\mbox{ and }\alpha\in\Delta^b_1 \\
a_D & \mbox{if }D\in\Delta_{\Sigma'}(\alpha)\cap\Delta_{\Sigma'}(\beta)\mbox{ and }\{\alpha,\beta\}\in\Delta^b_2 \\
\frac12a_D & \mbox{if }D\in\Delta_{\Sigma'}(\alpha)\setminus\Delta_{\Sigma'}(\beta)\mbox{ and }\{\alpha,\beta\}\in\Delta^b_2 
\end{array}\right.\]
For all $\sigma\in\Sigma'$ one has $c(F,\sigma)=c(E,\sigma)\geq0$, therefore $\widetilde\Delta|_{\Sigma'}$ is distinguished in $\s_{\Sigma'}$. 

The equality (\ref{eqn:reseq}), without assuming that the triple $(\s_{\Sigma'})/(\widetilde\Delta|_{\Sigma'})$ is a spherical system, is easily checked. Finally, $\widetilde\Delta|_{\Sigma'}$ is good because $\widetilde\Delta$ is good.
\end{proof}

\begin{corollary}
Let $f\colon X\to Y$ be a $G$-equivariant morphism with connected fibers between wonderful $G$-varieties, and let $X'\subseteq X$ be a wonderful $G$-subvariety. 
Then $Y'=f(X')$ is a wonderful $G$-subvariety of $Y$, and
\[
\s_{Y'} = \s_{X'} / (\Delta_f|_{\Sigma_{X'}}).
\]
\end{corollary}

\begin{definition}
A spherical $G$-system $\s = (S^p, \Sigma, \A)$ is {\em cuspidal} if $\supp\,\Sigma=S$.
\end{definition}
\begin{lemma}[{\cite[\S3.4]{Lu01}}]
Let $\s$ be a geometrically realizable spherical $G$-system with generic stabilizer $H$, and suppose that $H$ does not contain any simple factor of $G$. Then $\s$ is not cuspidal if and only if there exists a proper parabolic subgroup $P$ of $G$ such that $P^r\subseteq H \subseteq P$.
\end{lemma}

If the spherical system of a wonderful $G$-variety $X$ is not cuspidal, 
then $X$ is obtained from its localization to $S'=(S^p\cup \supp\,\Sigma) \subset S$ by {\em parabolic induction}: 
$X=G\times_{P_{S'}}X_{S'}$ (see \cite[\S3.4]{Lu01} for details). 
We also say that $X$ is obtained by parabolic induction {\em by means} of the parabolic subgroup $P_{S'}$ of $G$.

We end the section with a remark on the behaviour under localization of the {\em defect} of a spherical system.

\begin{definition}
Let $\s=(S^p,\Sigma,\A)$ be a spherical $G$-system with set of colors $\Delta$. The {\em defect} of $\s$ is the difference $|\Delta|-|\Sigma|$, also denoted by $\defect(\s)$.
\end{definition}

\begin{lemma}\label{lemma:defloc}
Let $\s=(S^p,\Sigma,\A)$ be a spherical $G$-system and $\Sigma'\subset \Sigma$ with $|\Sigma'|=|\Sigma|-1$. Then $0 \leq\defect(\s_{\Sigma'}) - \defect(\s)\leq 2$.
\end{lemma}
\begin{proof}
Denote by $\Delta_{\Sigma'}$ the set of colors of $\s_{\Sigma'}$. It is enough to show that $-1\leq |\Delta_{\Sigma'}|- |\Delta|\leq 1$. Let $\Sigma\setminus\Sigma'=\{\sigma\}$. If $\sigma$ is the sum of two orthogonal simple roots, then in the localization the color it is associated to becomes two distinct colors in $(\Delta_{\Sigma'})^{b}$, which implies $|\Delta_{\Sigma'}|-|\Delta| = 1$. If $\sigma$ is not a simple root nor the sum of two orthogonal simple roots, then $|\Delta_{\Sigma'}|=|\Delta|$. Finally, if $\sigma$ is a simple root, then either both its colors are free, which implies that they both become a single color in the localization, i.e.\ $|\Delta_{\Sigma'}|-|\Delta| = -1$, or at least one of them is not free, which implies that $|\Delta_{\Sigma'}|=|\Delta|$.
\end{proof}

\section{From the spherical system to the wonderful subgroup}\label{s:fromto}

\subsection{The Levi part}\label{s:levi}
Let $\s$ be a spherical $G$-system. 
The problem of determining a Levi part of the wonderful subgroup $H$ of $G$ associated with $\s$ 
(or a candidate, if geometric realizability of $\s$ is not assured) 
is related to the study of a special class of spherical subgroups of $G$: 
those of the form $K = L_K P^u$, where $P$ is a parabolic subgroup of $G$ containing $K$, 
and $L_K$ is a very reductive (i.e.\ contained in no proper parabolic subgroup) spherical subgroup of a Levi of $P$.

\begin{definition}\label{def:parared}
Let $\s$ be a spherical $G$-system with colors $\Delta$. Two distinguished sets $\Delta'$ and $\Delta''$ with $\Delta'\subseteq\Delta''\subseteq \Delta$ give a {\em decomposition of $\Delta$ of type ($\mathscr L$ $\mathscr R$ $\mathscr P$)}\footnote{For a motivation of this terminology see \cite[\S2.3.5]{BL11}.} if:
\begin{enumerate}
\item $\Delta''$ is parabolic, and minimal having this property;
\item\label{def:parared:ind} no simple root in the support of $\Sigma/\Delta'$ moves a color in $\Delta\setminus\Delta''$;
\item\label{def:parared:red} there exists a linear combination of elements in $\Sigma/\Delta'$, with positive coefficients, that takes non-negative values on all colors in $\Delta''\setminus\Delta'$;
\item $\Delta'$ is minimal with the above properties.
\end{enumerate}
\end{definition}

\begin{definition}\label{def:combL}
A distinguished set of colors $\Delta'$ of a spherical $G$-system is {\em of type ($\mathscr L$)} 
if it is contained in sets of colors $\widetilde\Delta'$ and $\Delta''$ 
such that $\widetilde\Delta'$ and $\Delta''$ give a decomposition of type ($\mathscr L$ $\mathscr R$ $\mathscr P$). 
If in addition $\Delta'$ is good, then we say that the quotient $\s/\Delta'$ is {\em of type ($\mathscr L$)}. 
If $\widetilde\Delta'$ can be chosen to be $\Delta'$, then $\Delta'$ is {\em maximal of type ($\mathscr L$)}.
\end{definition}

\begin{proposition}\label{prop:parared}
Suppose that the spherical $G$-system $\s$, with set of colors $\Delta$, is geometrically realizable with generic stabilizer $H$. Then maximal subsets $\Delta'\subseteq \Delta$ of type ($\mathscr L$) are in bijection with co-connected $K \supseteq H$ with the following properties:
\begin{enumerate}
\item\label{prop:parared:ind} there exists a minimal parabolic subgroup $P$ of $G$ containing $H$, such that $P\supseteq K \supseteq P^r$;
\item\label{prop:parared:red} there exist Levi components $M\supseteq L_K\supseteq L$ of resp.\ $P$, $K$ and $H$, such that
\begin{enumerate}
\item\label{prop:parared:red:veryred} $L_K$ and $L$ are very reductive subgroups of $M$, so in particular $K = L_K P^u$,
\item\label{prop:parared:center} $Z(L)^\circ\subseteq Z(L_K)^\circ=L_P^{r}$ and $Z(L_K)^{\circ}L = L_K$.
\end{enumerate}
\end{enumerate}
\end{proposition}
\begin{proof}
Recall the bijection between co-connected inclusions and distinguished sets of colors of Proposition~\ref{prop:quotH}. Since $H$ is only defined up to conjugation, we may assume in the whole proof that $P\supseteq B$ and that $M$ is the standard Levi subgroup $L_P$ of $P$.

If we start with a co-connected inclusion $P\supseteq K\supseteq H$ with the above properties, 
it corresponds to a minimal parabolic set $\Delta''$ and a maximal set $\Delta'\subseteq \Delta''$ of type ($\mathscr L$). 
In particular, condition (\ref{def:parared:ind}) of Definition~\ref{def:parared} follows from the fact 
that $G/K$ is the parabolic induction by means of $P$ of the homogeneous space $L_P/L_K$, 
and condition (\ref{def:parared:red}) follows from the fact that $L_K$ is reductive. 
Indeed, the characterization of the reductive spherical subgroups of a reductive group given in \cite[Theorem~6.7]{Kn} 
can be stated in the following way: 
\begin{equation}\label{eq:red}
\begin{array}{l}\mbox{there exists a linear combination of spherical roots with positive}\\ 
\mbox{coefficients that takes non-negative values on all colors.}\end{array} 
\end{equation}

Let now $\Delta'$ and $\Delta''$ give a decomposition of $\Delta$ of type ($\mathscr L$ $\mathscr R$ $\mathscr P$). Let $K\supseteq H$ and $P\supseteq K \supseteq H$ be the corresponding co-connected inclusions. From Definition~\ref{def:parared} and the above condition (\ref{eq:red}) it follows that $K$ is obtained by parabolic induction by means of $P$, and the part (\ref{prop:parared:ind}) of the proposition is proved.

The minimality of $P$ containing $H$ implies that we can choose Levi components $L_K$ and $L$ such that $L_P\supseteq L_K\supseteq L$, and such that $L$ is very reductive in $L_P$. Hence $L_K$ is very reductive in $L_P$ too: this shows property (\ref{prop:parared:red:veryred}), and also implies $Z(L)^\circ\subseteq Z(L_K)^\circ = L_P^{r}$.

If now $Z(L_K)^{\circ}L \neq L_K$, then the group $\widetilde K = L P^{r}$ is strictly contained in $K$, it has finite index in its normalizer, and the inclusion $H\subseteq\widetilde K$ is co-connected. Then the corresponding subset of colors $\widetilde{\Delta'}\subsetneq \Delta'$ is strictly contained in $\Delta'$. By an argument similar to the first part of the proof, the sets $\widetilde {\Delta'}$ and $\Delta''$ give a decomposition of $\Delta$ of type ($\mathscr L$ $\mathscr R$ $\mathscr P$) thus contradicting the minimality of $\Delta'$.
\end{proof}

\begin{definition}
Let $M, N$ be reductive groups. If $M\subset N$ and $Z(N)^{\circ}M = N$, then we say that $M$ and $N$ {\em differ only by the connected center}. 
\end{definition}

The above proposition gives a procedure to determine a Levi part of $H$, up to its connected center, 
using the spherical $G$-system $\s$. 
Namely, let $\Delta'$ be a maximal good distinguished set of colors of type ($\mathscr L$). 
The quotient $\s/\Delta'$ is the parabolic induction of a spherical $M$-system $\t$ 
satisfying the above condition (\ref{eq:red}).
Once $\t$ is computed, its wonderful subgroup $L_K\subseteq M$ 
is obtained from the classification of the reductive wonderful subgroups of reductive groups 
and their spherical systems (\cite{BP11p3}).

Notice that this procedure is well-defined even if $\s$ is not assumed to be geometrically realizable.

\subsection{The connected center}\label{s:c}
Let $\s=(S^p, \Sigma, \A)$ be a geometrically realizable spherical $G$-system with generic stabilizer $H$ and set of colors $\Delta$.

Conjugating $H$ in $G$ if necessary, let us choose a parabolic subgroup $Q_-$ containing $B_-$ and $H$, such that $Q_-$ is minimal having this property. The inclusion $Q_-\supseteq H$ is associated to a distinguished minimal parabolic set of colors $\Delta_{Q_-}$.

From minimality, it follows that the projection of $H$ in $Q_-/Q_-^u$ is not contained in any parabolic subgroup of $Q_-/Q_-^u$. The projection of $H$ in $Q_-/Q_-^u$ is therefore reductive, or equivalently: $Q_-^u\supseteq H^u$.

Recall that $L_{Q_-}$ denotes the standard Levi subgroup of $Q_-$, and let us choose a Levi part $L$ of $H$ inside $L_{Q_-}$. The considerations above tell us that $L$ is very reductive inside $L_{Q_-}$, and that $Z(L)^\circ$ is a subtorus of $Z(L_{Q_-})^\circ$.

The dimension of $Z(L)^\circ$ can be computed directly from the spherical system. Indeed, it is equal to $\rank\,\weights(L)$, and the latter is given by the well known formula (see \cite[\S5.2]{Lu01}):
\[
\rank\,\weights(L) = \defect(\s)=|\Delta|-|\Sigma|.
\]

We define the following subsets of $\Z\Sigma$:
\[
\rho(\Delta_{Q_-})^\perp = \{ \gamma\in\Z\Sigma \;|\; c(D,\gamma)=0 \;\forall D\in\Delta_{Q_-} \}
\]
\[
\rho(\Delta_{Q_-})^\vee = \{ \gamma\in\Z\Sigma \;|\; c(D,\gamma)\geq0 \;\forall D\in\Delta_{Q_-} \}.
\]

%For each color $D$ in $\Delta\setminus\Delta_{Q_-}$ define $\alpha_D$ as the unique simple root moving $D$. This root is unique because the projection $G/H\to G/Q_-$ induces a bijection between $\Delta\setminus\Delta_{Q_-}$ and the set of colors of $G/Q_-$, and each color of $G/Q_-$ is moved by exactly one simple root. Define also $\omega_D$ to be the fundamental dominant weight associated to $\alpha_D$.

\begin{lemma}\label{lemma:c}
The map
\[
\rho(\Delta_{Q_-})^\perp \to \mathcal X\left(Z(L_{Q_-})^\circ\right)^{L\cap Z(L_{Q_-})^\circ} %\cong \mathcal X\left(\frac{Z(L_{Q_-})^\circ}{L\cap Z(L_{Q_-})^\circ}\right)
\]
given by restricting weights to $Z(L_{Q_-})^\circ$ is an isomorphism.
\end{lemma}
\begin{proof}
Let $\gamma\in\rho(\Delta_{Q_-})^\perp$ and consider a $B$-eigenvector $f_\gamma\in\C(G/H)$ having $B$-eigenvalue $\gamma$. Call $F_\gamma$ the pull-back of $f_\gamma$ on $G$ along the projection $G\to G/H$.

Since $f_\gamma$ has neither zeros nor poles on colors of $\Delta_{Q_-}$, $F_\gamma$ has neither zeros nor poles on the pull-backs on $G$ of these colors. Its zeros or poles must then lie only on the pullbacks on $G$ of colors of $G/Q_-$ along the projection $G\to G/Q_-$.

This means that $F_\gamma$ is a $Q_-$-eigenvector under right translation. Since by construction
\[
F_\gamma \in V(\gamma)^{B^u}\otimes \left(V(\gamma)^*\right)^{B^u_-} \subset \C[G],
\]
we conclude that the $Q_-$-eigenvalue of $F_\gamma$ is $-\gamma$. But $F_\gamma$ is also of course $H$-stable under right translation, and this implies that  $(-\gamma)|_{Z(L_{Q_-})^\circ}$ is constantly $1$ on $L\cap Z(L_{Q_-})^\circ$.

From this we also deduce that the map $\rho(\Delta_{Q_-})^\perp \to \mathcal X\left(Z(L_{Q_-})^\circ\right)^{L\cap Z(L_{Q_-})^\circ}$ induced by restriction to $Z(L_{Q_-})^\circ$ is injective.

To show surjectivity, it is enough to notice that any character $\chi$ of $Z(L_{Q_-})^\circ$ is the $Q_-$-eigenvalue of some $Q_-$-eigenvector $F\in\C[G]$ (under right translation), such that $F$ is also a $B$-eigenvector under left translation, with $B$-eigenvalue $-\chi$. If $\chi$ is trivial on $L\cap Z(L_{Q_-})^\circ$ then $F$ is also $H$-stable under right translation, and descends to a $B$-eigenvector $f\in\C(G/H)$ with $B$-eigenvalue $-\chi$. By construction, then $-\chi\in\rho(\Delta_{Q_-})^\perp$.
\end{proof}
\begin{corollary}\label{cor:c1}
We have the following equality between subgroups of $Z(L_{Q_-})^\circ$:
\[
Z(L)^\circ = \left(\bigcap_{\gamma\in\rho(\Delta_{Q_-})^\perp} \ker (\gamma|_{Z(L_{Q_-})^\circ})\right)^\circ.
\]
\end{corollary}
\begin{proof}
It directly follows from Lemma~\ref{lemma:c}.
\end{proof}
The right hand side of the above equality is well-defined even if we do not assume that $\s$ is geometrically realizable. 
Hence this result, together with Section~\ref{s:levi}, provides a candidate Levi subgroup $L$ of $H$ 
defined for any spherical $G$-system $\s$.
\begin{corollary}\label{cor:c2}
The codimension of $Z(L)^\circ$ inside $Z(L_{Q_-})^\circ$ equals the codimension of the convex cone $\Q_{\geq0}\rho(\Delta_{Q_-})$ inside $\Hom_\Z(\Z\Sigma,\Q)$.
\end{corollary}
\begin{proof}
From Lemma~\ref{lemma:c} and Corollary~\ref{cor:c1} we have that the rank of the lattice $\rho(\Delta_{Q_-})^\perp$ is equal to the codimension of $Z(L)^\circ$ inside $Z(L_{Q_-})^\circ$.

It remains to notice that $\rho(\Delta_{Q_-})^\perp$ is the linear part of the convex cone $\rho(\Delta_{Q_-})^\vee$, and that the latter is the dual of the convex cone $\Q_{\geq0}\rho(\Delta_{Q_-})$.
\end{proof}

\subsection{The unipotent radical}
We maintain the assumptions of the previous section. We know that the $L$-module $\Lie\, Q^u_-/\Lie\,H^u$ is spherical (see \cite[\S5.4, auxiliary result (*)]{Lu01}), and the results of \cite{Lo07} imply that the spherical system $\s$ determines the structure of $\Lie\, Q^u_-/\Lie\,H^u$ as an $L$-module (see \cite[Proof of Theorem 3, step 5]{Lo10} for more details).

However, there is no known direct technique to compute the $L$-module structure of $\Lie\,Q^u_-/\Lie\,H^u$ from $\s$ unless one uses case-by-case considerations, for example based on the list of affine smooth spherical varieties of \cite{KVS06}.

On the other hand, the spherical system provides immediate information on the fiber $Q_-/H$ of the map $G/H\to G/Q_-$. It is easy to see that this fiber is $L_{Q_-}$-equivariantly isomorphic to the homogeneous vector bundle $L_{Q_-}\times_L \Lie\,Q^u_-/\Lie\,H^u$, so it is a smooth affine variety. It is also spherical under the action of $L_{Q_-}$, because $BQ_-/Q_-$ is dense in $G/Q_-$, therefore $Q_-/H$ intersects the dense $B$-orbit of $G/H$. It follows that the Borel subgroup $B\cap Q_-$ of $L_{Q_-}$ has a dense orbit on $Q_-/H$. 
\begin{definition}
We denote by $\Lambda^+(Q_-/H)$ the {\em weight monoid} of $Q_-/H$ with respect to the Borel subgroup $B\cap Q_-$ of $L_{Q_-}$, i.e.\ the set of $\chi\in\mathcal X(B\cap L_{Q_-})$ such that $\C[Q_-/H]^{(B\cap L_{Q_-})}_\chi \neq \{0\}$.
\end{definition}
\begin{lemma}\label{lemma:lambda}
We have $\Lambda^+(Q_-/H) = \rho(\Delta_{Q_-})^\vee$.
\end{lemma}
\begin{proof}
Let $f\in\C(G/H)^{(B)}$, and consider its restriction $f|_F$ on the fiber $F=Q_-/H\subseteq G/H$. The restriction is a $(B\cap L_{Q_-})$-eigenvector, and it is non-zero because $Q_-$ is $B$-spherical, therefore no color of $G/H$ contains $F$.

We obtain in this way a bijection between $\C(G/H)^{(B)}$ and $\C(F)^{(B\cap L_{Q_-})}$. This follows from \cite[Proposition~I.2]{Br87}; let us provide here a direct proof. Injectivity is obvious, since $T$ is a maximal torus of $(B\cap L_{Q_-})$, therefore the $(B\cap L_{Q_-})$-eigenvalue of $f|_F$ is equal to the $B$-eigenvalue of $f$. Let now $f_0\in\C(Q_-/H)^{(B\cap L_{Q_-})}$, and let us extend it to a function $f\in\C(G/H)^{(B)}$ with the following definition:
\[
f(uqH)=f_0(qH)
\]
for $u\in B^u$ and $q\in Q_-$. From the local structure theorem of spherical varieties (see \cite{BL87}) it follows that $f$ is a well defined rational function on the open set $B^uQ_-/H$ of $G/H$, and that it is a $B$-eigenvector.

Now, functions in $\C[F]^{(B\cap L_{Q_-})}$ are exactly the restrictions of functions on $G/H$ having no pole on colors in $\Delta_{Q_-}$, and the lemma follows.
\end{proof}

A further step, going beyond the scope of the present work, 
would be to study all possible spherical quotients of $\Lie\,Q^u_-$ 
under the action of spherical very reductive subgroups of $L_{Q_-}$. 
Taking into account the results of \S\ref{s:higherdefect}, 
it would be enough to assume that the very reductive subgroup contains $Z(L_{Q_-})^\circ$: 
this would provide a strategy to guess $H^u$ inside $Q_-^u$ based only on $\s$.

\section{Wonderful fiber products}\label{s:fiberprod}

\subsection{}
Let $X_1,X_2,X_{1,2}$ be wonderful $G$-varieties and $\varphi_1\colon X_1\to X_{1,2}$, $\varphi_2\colon X_2\to X_{1,2}$ surjective equivariant maps with connected fibers. Then the fiber product $X=X_1\times_{X_{1,2}}X_2$ is a (not necessarily spherical) $G$-variety, with $\psi_1\colon X\to X_1$, $\psi_2\colon X\to X_2$ surjective equivariant maps with connected fibers such that $\varphi_1\circ\psi_1=\varphi_2\circ\psi_2$. 

Let $Z_1,Z_2,Z_{1,2}$ be the corresponding closed $G$-orbits, with the corresponding restricted maps. 

\begin{definition}
The $G$-variety $X=X_1\times_{X_{1,2}}X_2$ is called a {\em wonderful fiber product} if it is wonderful with closed $G$-orbit $Z\cong Z_1\times_{Z_{1,2}}Z_2$.
\end{definition}

\begin{definition}\label{def:dec}
Let $\s$ be a spherical $G$-system. Two good distinguished sets of colors $\Delta_1,\Delta_2$ are said to {\em decompose} $\s$ if:
\begin{itemize}
\item ($S^p/\Delta_1\setminus S^p) \perp (S^p/\Delta_2\setminus S^p$),
\item $\Sigma\subset(\Sigma/\Delta_1\cup\Sigma/\Delta_2)$. 
\end{itemize}
If $\s$ admits two non-empty such subsets of colors then $\s$ is {\em decomposable}.
\end{definition}

If $\Delta_1,\Delta_2$ decompose $\s$ then there is no $\sigma\in\Sigma$, $D_1\in\Delta_1$, $D_2\in\Delta_2$ such that both $c(D_1,\sigma)$ and $c(D_2,\sigma)$ are non-zero. Therefore, $\card(\Sigma)+\card(\Sigma/\Delta_1\cup\Delta_2)=\card(\Sigma/\Delta_1)+\card(\Sigma/\Delta_2)$. Moreover, $P_{S^p}=P_{S^p/\Delta_1}\cap P_{S^p/\Delta_2}$ and $P_{S^p/(\Delta_1\cup\Delta_2)}=P_{S^p/\Delta_1}P_{S^p/\Delta_2}=P_{S^p/\Delta_2}P_{S^p/\Delta_1}$.

\begin{theorem}\label{thm:fiberprod}
Let $X$ be a wonderful $G$-variety with spherical system $\s$. Then $X$ is isomorphic to a wonderful fiber product $X_1\times_{X_{1,2}}X_2$ (with notation as above) if and only if the corresponding good distinguished sets $\Delta_{\psi_1},\Delta_{\psi_2}$ decompose $\s$.
\end{theorem}

\begin{proof}
Recall the local structure of a wonderful $G$-variety $X$ (\cite{BL87}): $X_B=X\setminus\cup_{D\in\Delta}D$ is an open affine subvariety $P_{S^p}$-isomorphic to $P_{S^p}^u\times M$, where $M$ is a smooth affine $T$-variety isomorphic to $\C^{\Sigma}$. A surjective $G$-equivariant map with connected fibers between wonderful $G$-varieties $X\to X'$ restricts to a $B$-equivariant map $X_B\to X'_B$ and, moreover, to a $T$-equivariant map between affine $T$-varieties $\C^\Sigma\to\C^{\Sigma/\Delta'}$ corresponding to the inclusion $\IN\,\Sigma/\Delta'\subset\IN\,\Sigma$. 

Therefore, if $X$ is a wonderful fiber product then one has 
\begin{equation}\label{eq:canch}
\C^\Sigma=\C^{\Sigma/\Delta_1}\times_{\C^{\Sigma/\Delta_1\cup\Delta_2}}\C^{\Sigma/\Delta_2},
\end{equation}
which, together with the condition on the closed orbits 
\begin{equation}\label{eq:clorb}
Z\cong Z_1\times_{Z_{1,2}}Z_2,
\end{equation}
is equivalent to the conditions given in Definition~\ref{def:dec}.

Vice versa, once one has (\ref{eq:canch}) and (\ref{eq:clorb}), then consider the fiber product $X'=X_1\times_{X_{1,2}}X_2$. 
It contains the open affine $B$-subvariety $(X_1)_B\times_{(X_{1,2})_B}(X_2)_B$ which has an open $B$-orbit, 
then $X'$ is spherical. 
Furthermore, it is smooth since $M'=\C^\Sigma$ is, it is complete, 
and it is toroidal since every color of $X'$ maps non-dominantly into $X_1$ or $X_2$ (which are both toroidal). 
Therefore, $X'$ is wonderful and isomorphic to $X$.  
\end{proof}

From the above proof, it is clear that we have also the following.

\begin{theorem}
Let $\s$ be a spherical system, and let $\Delta_1,\Delta_2$ be good distinguished sets of colors that decompose $\s$. If $\s/\Delta_1$ and $\s/\Delta_2$ are geometrically realizable then $\s$ is geometrically realizable. 
\end{theorem}

If $X=X_1\times_{X_{1,2}}X_2$ is a wonderful fiber product, 
choose $x$ in the open orbit and set $x_1=\psi_1(x)$, $x_2=\psi_2(x)$, $x_{1,2}=\varphi_1(x_1)=\varphi_2(x_2)$. 
Let $H,H_1,H_2,H_{1,2}$ be the corresponding stabilizers, then $H=H_1\cap H_2$ and $H_{1,2}=H_1H_2=H_2H_1$, 
namely $G/H\cong G/H_1 \times_{G/H_{1,2}} G/H_2$.

\section{Quotients of type ($\mathscr L$)}\label{s:higherdefect}

\subsection{Combinatorial characterization of quotients of type ($\mathscr L$)}

\begin{definition}\label{def:geomL}
Let $H$ be a wonderful subgroup of $G$. A minimal co-connected inclusion 
$K\supset H$ is {\em of type ($\mathscr L$)} if $H^u$ is strictly contained in $K^u$, 
$\Lie\, K^u/\Lie\,H^u$ is a simple $H$-module and the Levi components of $H$ and $K$ differ only by the connected center.
\end{definition}

In this case $L\subseteq\N_{L_K}(H^u)$, where $L$ and $L_K$ denote the Levi components of $H$ and $K$, respectively. 
We show in \S\ref{s:wondsubL} that the equality holds, under certain combinatorial assumptions (Definition~\ref{def:combqhd}, see the proof of Theorem~\ref{thm:qhd}).

\begin{proposition}\label{prop:L}
Suppose that $\s$ is a geometrically realizable spherical $G$-system with set of colors $\Delta$ and generic stabilizer $H$. 
Then the bijection of Proposition~\ref{prop:quotH} restricts to a bijection between minimal distinguished subsets 
$\Delta'\subseteq \Delta$ of type ($\mathscr L$) and minimal co-connected inclusions $K\supset H$ of type ($\mathscr L$).
\end{proposition}
\begin{proof}
If we have a minimal co-connected inclusion $K\supset H$ of type ($\mathscr L$), the corresponding set of colors $\Delta_K$ is minimal distinguished, so we have only to show that it is of type ($\mathscr L$).

Let $P$ be a minimal parabolic subgroup of $G$ containing $H$. 
Since the minimal inclusion $K\supset H$ is of type ($\mathscr L$), 
up to conjugation we can suppose that $P\supseteq K \supset H$, that $P^u\supseteq K^u \supset H^u$, 
and that we have chosen Levi parts $L_P\supseteq L_K \supseteq L$ 
in such a way that $L_K$ and $L$ differ only by their connected centers.

Consider $K'=L P^r$. We have $P\supseteq K'\supseteq K\supseteq H$, and one checks using Proposition~\ref{prop:parared} that $K'\supseteq H$ corresponds to $\Delta_{K'}$, a maximal set of colors of type ($\mathscr L$) containing $\Delta_K$.

Vice versa, suppose we have a maximal set of type ($\mathscr L$) containing a minimal distinguished set of colors. 
From geometric realizability we deduce that they are associated to $K'\supseteq K\supset H$. 
We can thus call our sets of colors resp.\ $\Delta_{K'}$ and $\Delta_{K}$.

Since $\Delta_{K'}$ is maximal of type ($\mathscr L$), we know that Levi parts of $H$ and $K'$ differ only by their connected centers. 
This must happen for $K$ and $H$ too, since we can suppose to have Levi parts $L_{K'}\supseteq L_K\supseteq L$.
It also follows that $K^u\supset H^u$. The minimality of $\Delta_K$ assures that $H^u$ contains $(K^u, K^u)$: 
indeed, the projection of $H^u$ on $K^u/(K^u,K^u)$ cannot be surjective because otherwise we would have $H^u= K^u$. 
Then $H^u\supseteq (K^u,K^u)$.
Hence $\Lie\, H^u$ is an ideal of $\Lie\,K^u$ and $\Lie\,K^u /\Lie\,H^u$ is abelian: 
again the minimality of $\Delta_K$ implies that it is a simple $L$-module. 
Therefore $K\supset H$ is of type ($\mathscr L$).
\end{proof}

\subsection{The wonderful subgroup associated with a spherical system admitting a quotient of type ($\mathscr L$)}
Let $\s$ be a spherical $G$-system admitting a geometrically realizable quotient $\s/\Delta'$ of type ($\mathscr L$). 
We propose here a conjectural construction of the wonderful subgroup $H$ of $G$ associated with $\s$ 
starting from the wonderful subgroup of $G$ associated with $\s/\Delta'$.

The interest of this procedure is its great generality, 
and the fact that its proof reduces to the verification of the conjecture below, 
which can be easily checked on any given example. 
%At present, however, the only known approach to prove the two conjectures would be case-by-case 
%after reducing the problem to the list of primitive spherical systems.

We introduce an additional hypothesis which actually makes the description of $H$ easier 
and excludes only a few special minimal quotients of type ($\mathscr L$). 

A color of $\s$ is called {\em negative} if is $\leq0$ on the set of spherical roots of $\s$ 
(in this case it is moved by a unique simple root). 
A negative color is called {\em exterior} 
if is moved by a simple root not belonging to the support of a spherical root. 

Let $\Delta'$ be any minimal good distinguished set of colors of $\s$, such that 
\begin{itemize}
\item[-] there exist exterior negative colors of $\s/\Delta'$ that are not exterior or not negative as colors of $\s$,
\item[-] every negative color of $\s/\Delta'$, that is not negative as color of $\s$, is exterior.
\end{itemize} 
Then $\Delta'$ is of type ($\mathscr L$), see \cite[\S2.3.5]{BL11}.

Let $\Delta'$ be a minimal good distinguished set of colors of $\s$ fulfilling the above conditions, 
and assume that $\s/\Delta'$ is geometrically realizable.
Let $S_{\Delta'}$ be the set of simple roots moving exterior negative colors of $\s/\Delta'$ 
that are not exterior or not negative as colors of $\s$. 
Let $Q_-$ be the parabolic subgroup of $G$ corresponding to $S\setminus S_{\Delta'}$ 
and let $K$ be the wonderful subgroup of $G$ associated with $\s/\Delta'$: 
we can assume $Q^r_-\subset K\subset Q_-$. 
Set a decomposition $Q_-=Q^r_-M$, we have $K=Q^r_-K_{\Delta'}$ with $K_{\Delta'}\subset M$ 
(notice that $K_{\Delta'}$ is not necessarily reductive).

Let us consider the spherical $M$-systems $\t$ and $\t'$ obtained by localization of $\s$ and $\s/\Delta'$, resp., 
in the set of simple roots $S\setminus S_{\Delta'}$. 
The spherical system $\t'$ equals the quotient spherical system $\t/\Delta'|_{S\setminus S_{\Delta'}}$,
and is the spherical system of $K_{\Delta'}$ which is a wonderful subgroup of $M$. 

Let us assume that $\t$ is geometrically realizable and 
let $H_{\Delta'}$, included in $K_{\Delta'}$, be the wonderful subgroup of $M$ associated with $\t$.

The following determines $H$ by describing the intersection of its unipotent radical $H^u$ with the unipotent radical of $Q_-$.
For all $\alpha\in S_{\Delta'}$ let $V(-\alpha)$ be the simple $M$-module of highest weight $-\alpha$, 
and let $V$ be equal to $\oplus_{\alpha\in S_{\Delta'}}V(-\alpha)$,
notice that $\Lie\,Q^u_-=V\oplus [\Lie\,Q^u_-,\Lie\,Q^u_-]$.

\begin{conjecture}\
\begin{enumerate}
\item\label{conj:W}
There exists a simple $K_{\Delta'}$-module $W$ 
such that $\p(W)$ contains an open $K_{\Delta'}$-orbit isomorphic to $K_{\Delta'}/H_{\Delta'}$.
\item\label{conj:V}
Let $W^\ast$ be the dual of $W$.
There exists a $K_{\Delta'}$-equivariant inclusion $W^\ast\subseteq V(-\alpha)^\ast$, for all $\alpha\in S_{\Delta'}$.
\end{enumerate}
\end{conjecture}

It would follow that $\Lie\,H^u\cap \Lie\,Q^u_-$ is the co-simple $K_{\Delta'}$-submodule of $\Lie\,Q^u_-$ 
dual to $W^\ast$ diagonally included in $V^\ast$ (compare with Theorem~\ref{thm:qhd}).

Without going into technical details, the above conjecture has a combinatorial counterpart
in terms of spherical systems and spherical orbits in simple projective spaces, as follows. 

The space $\p(W^\ast)$ would still have an open spherical $K_{\Delta'}$-orbit. 
Let $H^\ast_{\Delta'}\subset K_{\Delta'}$ be the generic stabilizer, 
namely the open $K_{\Delta'}$-orbit of $\p(W^\ast)$ is isomorphic to $K_{\Delta'}/H^\ast_{\Delta'}$. 
Provided $H^\ast_{\Delta'}$ is a wonderful subgroup of $M$, 
it would correspond to a spherical system $\t^\ast$ admitting $\t/\Delta'$ as quotient.
Moreover, in this case, 
$\p(V(-\alpha)^\ast)$ would contain a $K_{\Delta'}$-orbit isomorphic to $K_{\Delta'}/H^\ast_{\Delta'}$,
and this is equivalent to a combinatorial statement on the spherical system $\t^\ast$, 
see \cite[\S2.4]{BL11} for details.

%Therefore, in this case, $H=(H^u\cap Q^u_-)\,K_{\Delta'}\,\N_{Z(Q_-)}(H^u\cap Q^u_-)$. 

Let us give some examples of minimal quotients of type ($\mathscr L$), further examples can be found in \cite[\S3.7 and \S3.8]{BL11}.

\subsubsection{}
Here and in the next sections we use certain graphical representations of spherical $G$-systems called {\em Luna diagrams}, we refer to \cite{B09} for their definition.

Let $G$ be $\Sp(6)$. 
Consider the spherical $G$-system $\s$ with $S^p=\emptyset$, $\Sigma=\{\alpha_1+\alpha_2, \alpha_2+\alpha_3\}$ 
and $\A=\emptyset$. 
The set of colors is $\Delta=\{D_{\alpha_1}, D_{\alpha_2}, D_{\alpha_3}\}$ with Cartan pairing as follows.
{\small
\[\begin{array}{l|rr}
c(-,-)&\alpha_1+\alpha_2&\alpha_2+\alpha_3\\
\hline
D_{\alpha_1}& 1\rule{10pt}{0pt}&-1\rule{10pt}{0pt}\\
D_{\alpha_2}& 1\rule{10pt}{0pt}& 0\rule{10pt}{0pt}\\
D_{\alpha_3}&-1\rule{10pt}{0pt}& 1\rule{10pt}{0pt}
\end{array}\]
}
The set $\Delta'=\{D_{\alpha_2}\}$ is distinguished of type ($\mathscr L$). 
The corresponding quotient $\s/\Delta'$ is such that 
$S^p/\Delta'=\{\alpha_2\}$, $\Sigma/\Delta'=\{\alpha_2+\alpha_3\}$ and $\A=\emptyset$.
\[\begin{picture}(12600,1800)(-300,-900)
\put(0,0){\usebox{\atwo}}
\put(1800,0){\usebox{\leftbiedge}}
\put(3600,0){\usebox{\gcircle}}
\put(4500,0){\vector(1,0){3000}}
\put(8700,0){
\put(0,0){\usebox{\dynkincthree}}
\put(0,0){\usebox{\wcircle}}
\put(3600,0){\usebox{\gcircle}}
}
\end{picture}\]
Here $S_{\Delta'}=\{\alpha_2\}$, $K_{\Delta'}\cong\SL(2)\times\SL(2)$, 
$H_{\Delta'}$ is a codimension 1 parabolic subgroup of $K_{\Delta'}$ 
and $W$ is isomorphic to the standard representation of one simple factor of $K_{\Delta'}$.

\subsubsection{} Let $G$ be $\Sp(10)$. 
Consider the spherical $G$-system $\s$ with $S^p=\emptyset$, $\Sigma=S$, 
$\A=\Delta=\{D^+_{\alpha_1}=D^+_{\alpha_3}=D^+_{\alpha_5}, D^+_{\alpha_2}, D^+_{\alpha_4}, D^-_{\alpha_1}=D^-_{\alpha_4}, D^-_{\alpha_2}=D^-_{\alpha_5}, D^-_{\alpha_3}\}$ 
and Cartan pairing as follows.
{\small
\[\begin{array}{l|rrrrr}
c(-,-)&\alpha_1&\alpha_2&\alpha_3&\alpha_4&\alpha_5\\
\hline
D^+_{\alpha_1}& 1&-1& 1&-1& 1\\
D^+_{\alpha_2}& 0& 1& 0& 0&-1\\
D^+_{\alpha_4}&-1& 0& 0& 1&-1\\
D^-_{\alpha_1}& 1& 0&-1& 1&-1\\
D^-_{\alpha_2}&-1& 1&-1& 0& 1\\
D^-_{\alpha_3}&-1& 0& 1& 0&-1
\end{array}\]
}
The set $\Delta'=\Delta\setminus\{D^-_{\alpha_3}\}$ is distinguished of type ($\mathscr L$). 
The corresponding quotient $\s/\Delta'$ is such that 
$S^p/\Delta'=S\setminus\{\alpha_3\}$, $\Sigma/\Delta'=\emptyset$ and $\A=\emptyset$.
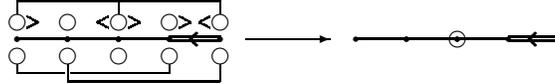
\[\begin{picture}(19800,2850)
\put(0,0){\diagramcesevenfive}
\put(8400,1500){\vector(1,0){3000}}
\put(12300,1500){
\put(0,0){\usebox{\dynkincfive}}\put(3600,0){\usebox{\wcircle}}
}
\end{picture}\]
In this case $S_{\Delta'}=\{\alpha_3\}$, $Q_-=K$ and $K_{\Delta'}=M\cong\SL(3)\times\Sp(4)$. 
One gets $W=V(-\alpha_3)$ (this would directly follow from the conjecture by dimension reasons). 
To find $H_{\Delta'}$ one has to reapply the same construction to a quotient of type ($\mathscr L$) of $\t$. 
Let us reset our notation.

Let $G$ be $\SL(3)\times\Sp(4)$. 
Consider the spherical $G$-system $\s$ with $S^p=\emptyset$, $\Sigma=S$, 
$\A=\Delta=\{D^+_{\alpha'_1}=D^+_{\alpha''_2}, D^+_{\alpha'_2}, D^+_{\alpha''_1}, D^-_{\alpha'_1}=D^-_{\alpha''_1}, D^-_{\alpha'_2}=D^-_{\alpha''_2}\}$ 
and Cartan pairing as follows.
{\small
\[\begin{array}{l|rrrr}
c(-,-)&\alpha'_1&\alpha'_2&\alpha''_1&\alpha''_2\\
\hline
D^+_{\alpha'_1} & 1&-1&-1& 1\\
D^+_{\alpha'_2} & 0& 1& 0&-1\\
D^+_{\alpha''_1}&-1& 0& 1&-1\\
D^-_{\alpha'_1} & 1& 0& 1&-1\\
D^-_{\alpha'_2} &-1& 1& 0& 1
\end{array}\]
}
The set $\Delta'=\{D^-_{\alpha'_1}, D^-_{\alpha'_2}\}$ is distinguished of type ($\mathscr L$). 
The corresponding quotient $\s/\Delta'$ is such that 
$S^p/\Delta'=\emptyset$, $\Sigma/\Delta'=\{\alpha'_1+\alpha''_2\}$ and $\A=\emptyset$.
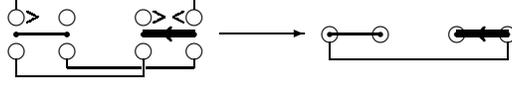
\[\begin{picture}(18000,2850)(-300,-1500)
\put(0,0){\usebox{\edge}}
\put(4500,0){\usebox{\leftbiedge}}
\multiput(0,0)(4500,0){2}{\multiput(0,0)(1800,0){2}{\usebox{\aone}}}
\multiput(0,1350)(6300,0){2}{\line(0,-1){450}}
\put(0,1350){\line(1,0){6300}}
\multiput(0,-1500)(4500,0){2}{\line(0,1){600}}
\put(0,-1500){\line(1,0){4500}}
\multiput(1800,-1200)(4500,0){2}{\line(0,1){300}}
\put(1800,-1200){\line(1,0){2600}}
\put(4600,-1200){\line(1,0){1700}}
\multiput(0,600)(4500,0){2}{\usebox{\toe}}
\put(6300,600){\usebox{\tow}}
\put(7200,0){\vector(1,0){3000}}
\put(11100,0){
\put(0,0){\usebox{\edge}}
\put(4500,0){\usebox{\vertex}}
\put(4500,0){\usebox{\leftbiedge}}
\multiput(0,0)(4500,0){2}{\multiput(0,0)(1800,0){2}{\usebox{\wcircle}}}
\multiput(0,-900)(6300,0){2}{\line(0,1){600}}
\put(0,-900){\line(1,0){6300}}
}
\end{picture}\]
Therefore, $S_{\Delta'}=\{\alpha'_2, \alpha''_1\}$, $M\cong\SL(2)\times\SL(2)$ 
and $K_{\Delta'}^\circ\cong\SL(2)$ included diagonally in $M$. 
Here $H_{\Delta'}$ is a Borel subgroup of $K_{\Delta'}$, $W$ is isomorphic to the standard representation of $K_{\Delta'}$ 
and $W\cong V(-\alpha'_2)\cong V(-\alpha''_1)$, hence $W$ is included diagonally in $V(-\alpha'_1)\oplus V(-\alpha''_2)$.

\subsubsection{} Let $G$ be $\SO(9)$. 
Consider the spherical $G$-system $\s$ with $S^p=\emptyset$, $\Sigma=\{\alpha_1,\alpha_2+\alpha_3,\alpha_4\}$ 
and $\A=\{D_{\alpha_1}^+=D_{\alpha_4}^+,D_{\alpha_1}^-,D_{\alpha_4}^-\}$ with $c(D^+_{\alpha_1},\alpha_2+\alpha_3)=-1$. 
The set of colors is $\Delta=\{D_{\alpha_1}^+, D_{\alpha_1}^-,D_{\alpha_2}, D_{\alpha_3},D_{\alpha_4}^-\}$ with Cartan pairing as follows.
{\small
\[\begin{array}{l|rrr}
c(-,-)&\alpha_1&\alpha_2+\alpha_3&\alpha_4\\
\hline
D_{\alpha_1}^+ & 1&-1\rule{10pt}{0pt}& 1\\
D_{\alpha_1}^- & 1& 0\rule{10pt}{0pt}&-1\\
D_{\alpha_2}   &-1& 1\rule{10pt}{0pt}& 0\\
D_{\alpha_3}   & 0& 1\rule{10pt}{0pt}&-1\\
D_{\alpha_4}^- &-1&-1\rule{10pt}{0pt}& 1
\end{array}\]
}
The set $\Delta'=\{D_{\alpha_1}^+,D_{\alpha_3}\}$ is distinguished of type ($\mathscr L$). 
The corresponding quotient $\s/\Delta'$ is such that 
$S^p/\Delta'=\{\alpha_3\}$, $\Sigma/\Delta'=\{\alpha_2+\alpha_3+\alpha_4\}$ and $\A=\emptyset$.
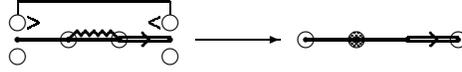
\[\begin{picture}(16200,2250)(-300,-900)
\put(0,0){\usebox{\dynkinbfour}}
\multiput(0,0)(5400,0){2}{\usebox{\aone}}
\multiput(0,900)(5400,0){2}{\line(0,1){450}}
\put(0,1350){\line(1,0){5400}}
\put(0,600){\usebox{\toe}}
\put(5400,600){\usebox{\tow}}
\put(1800,0){\usebox{\atwo}}
\put(6300,0){\vector(1,0){3000}}
\put(10200,0){
\put(0,0){\usebox{\dynkinbfour}}
\multiput(0,0)(5400,0){2}{\usebox{\wcircle}}
\put(1800,0){\usebox{\gcircle}}
}
\end{picture}\]
Here $S_{\Delta'}=\{\alpha_1\}$, $K_{\Delta'}$ is the parabolic subgroup of semisimple type $\mathsf A_2$ of $\SO(6)\subset M\cong\SO(7)$, 
$H_{\Delta'}$ is a codimension 2 parabolic subgroup of $K_{\Delta'}$ 
and $W$ is isomorphic to the standard representation of the semisimple part of $K_{\Delta'}$. 
Notice that $W^\ast$ is $K_{\Delta'}$-stable in $V^\ast\cong V$, the standard representation of $M$.

\subsubsection{} Let $G$ be $\SO(7)$. 
Consider the spherical $G$-system $\s$ with $S^p=\emptyset$, $\Sigma=\{\alpha_1+\alpha_2, \alpha_2+\alpha_3, \alpha_3\}$, 
$\A=\{D^+_{\alpha_3}, D^-_{\alpha_3}\}$ with $c(D^+_{\alpha_3},\sigma)=c(D^-_{\alpha_3},\sigma)$ for all $\sigma\in\Sigma$. 
Its set of colors is $\Delta=\{D_{\alpha_1}, D_{\alpha_2}, D^+_{\alpha_3}, D^-_{\alpha_3}\}$ with Cartan pairing as follows. 
{\small
\[\begin{array}{l|rrr}
c(-,-)&\alpha_1+\alpha_2&\alpha_2+\alpha_3&\alpha_3\\
\hline
D_{\alpha_1}  & 1\rule{10pt}{0pt}&-1\rule{10pt}{0pt}& 0\\
D_{\alpha_2}  & 1\rule{10pt}{0pt}& 1\rule{10pt}{0pt}&-1\\
D^+_{\alpha_3}&-1\rule{10pt}{0pt}& 0\rule{10pt}{0pt}& 1\\
D^-_{\alpha_3}&-1\rule{10pt}{0pt}& 0\rule{10pt}{0pt}& 1
\end{array}\]
}
The set $\Delta'=\{D_{\alpha_2}, D^+_{\alpha_3}\}$ is distinguished of type ($\mathscr L$). 
The corresponding quotient $\s/\Delta'$ is such that 
$S^p/\Delta'=\{\alpha_2\}$, $\Sigma/\Delta'=\{\alpha_1+\alpha_2+\alpha_3\}$ and $\A=\emptyset$.
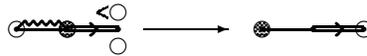
\[\begin{picture}(12600,1800)(-300,-900)
\put(0,0){\usebox{\atwo}}
\put(1800,0){\usebox{\btwo}}
\put(3600,0){\usebox{\aone}}
\put(3600,600){\usebox{\tow}}
\put(4500,0){\vector(1,0){3000}}
\put(8700,0){
\put(0,0){\usebox{\dynkinbthree}}
\put(0,0){\usebox{\gcircle}}
\put(3600,0){\usebox{\wcircle}}
}
\end{picture}\]
In this example, the above hypothesis on the negative colors of $\s/\Delta'$ is not verified. 
One has $H\cong \GL(3)$ included in $K$, which is a parabolic subgroup of semisimple type $\mathsf A_2$ of $\SO(6)\subset G$.

\subsection{The wonderful subgroup associated with a spherical system admitting a quotient of higher defect}\label{s:wondsubL}
Let $\s$ be a spherical $G$-system admitting a {\em minimal quotient of higher defect} $\s/\Delta'$, 
see Definition~\ref{def:combqhd}. 
Assuming the geometric realizability of $\s/\Delta'$ and of certain localizations of $\s$, 
we prove 
the geometric realizability of $\s$ and we describe its associated wonderful subgroup $H$ of $G$. 

\begin{definition}\label{def:combqhd}
Let $\s=(S^p,\Sigma,\A)$ be a spherical $G$-system with set of colors $\Delta$. 
Let $\Delta'$ be a minimal good distinguished set of type ($\mathscr L$) for $\s$ with $k=\defect(\s/\Delta')-\defect(\s)>0$.
The quotient $\s/\Delta'$ is of {\em higher defect} if there exist $k+1$ spherical roots $\sigma_0,\ldots,\sigma_k\in\Sigma$ such that, 
if we set, for all non-empty $I\subset\{0,\ldots,k\}$,  $\Sigma_I = \Sigma\setminus\{\sigma_i:i\not\in I\}$,  
$\s_I=\s_{\Sigma_I}$ and $\Delta'_I=\Delta'|_{\Sigma_I}$, we have: 
\begin{enumerate}
\item $\defect(\s_I)=\defect(\s)+k+1-|I|$,
\item $\Delta'_I$ is minimal of type ($\mathscr L$),
\item $\s_I/\Delta'_I=\s/\Delta'$.
\end{enumerate}
\end{definition}

Notice that if $\Delta'$ is minimal good distinguished and $\defect(\s/\Delta')>\defect(\s)$ 
then the quotient is of type ($\mathscr L$) (see \cite[\S2.3.5]{BL11}).

A direct (but long) combinatorial verification on the list of primitive spherical systems would actually show that all minimal good distinguished sets $\Delta'$ of type ($\mathscr L$) for a spherical $G$-system $\s$ with $\defect(\s/\Delta')-\defect(\s)>0$ always satisfy the above definition. We avoid it.

In the following, $\s/\Delta'$ is a quotient of higher defect, 
with notation as in Definition~\ref{def:combqhd}.
For all $i\in \{0,\ldots,k\}$ we will simply write $\Sigma_i$, $\s_i$ and $\Delta'_i$ 
instead of $\Sigma_{\{i\}}$, $\s_{\{i\}}$ and $\Delta'_{\{i\}}$, respectively. 
Furthermore, we set $\hat j = \{0,\ldots,k\}\setminus\{j\}$.

\begin{lemma}\label{lemma:combqhd}
Let $\Delta_I$ denote the set of colors of $\s_I=(S^p_I, \Sigma_I, \A_I)$ and $c_I$ its Cartan pairing. 
Let $\Delta''$ be a minimal parabolic subset of $\Delta$ containing $\Delta'$. 
Then for all non-empty $I\subset\{0,\ldots,k\}$ the following hold:
\begin{enumerate}
\item\label{lemma:combqhd:quot} the restriction $\Delta''_I$ of $\Delta''$ to $\s_I$ is a minimal parabolic set containing $\Delta'_I$ and $\s_I/\Delta''_I = \s/\Delta''$;
\item\label{lemma:combqhd:values} the natural inclusion $\A_I\to \A$ extends in a unique way to a bijection $\varphi_I\colon \Delta_I \to \Delta$ such that a simple root $\alpha$ moves $D\in\Delta_I\setminus\A_I$ if and only if $\alpha$ moves $\varphi_I(D)$;
\item\label{lemma:combqhd:bij} the bijection $\varphi_I$ satisfies
\begin{enumerate}
\item\label{lemma:combqhd:bij:images} $\varphi_I(\Delta''_I)=\Delta''$ and $\varphi_I(\Delta'_I)=\Delta'$,
\item\label{lemma:combqhd:bij:comp1} $c_I(D,\sigma) = c(\varphi_I(D),\sigma)$ for all $D\in\Delta''_I$ and all $\sigma\in\Sigma_I$.
\item\label{lemma:combqhd:bij:comp2} for all $D\in\Delta_I$ either $c_I(D,\sigma)=  c(\varphi_I(D),\sigma)$ for all $\sigma\in\Sigma_I$, or there exists $D'\in\Delta'$ such that
\[
c_I(D,\sigma)=  c(\varphi_I(D),\sigma) + c(D',\sigma)
\]
for all $\sigma\in\Sigma_I$.
\end{enumerate}
\end{enumerate}
\end{lemma}

\begin{proof}
First of all notice that, by Definition~\ref{def:combqhd} part (3), 
$\Sigma/\Delta'$ lies in the lattice generated by 
$\Sigma\setminus\{\sigma_0,\ldots,\sigma_k\}$.
Therefore, the first statement follows from Lemma~\ref{lem:restr}. 

From Definition~\ref{def:combqhd}, part (1),  
since $\s_I$ has $|\Sigma|-k-1+|I|$ spherical roots, 
we deduce that $|\Delta_I|=|\Delta|$ for all $I$. 
From the proof of Lemma~\ref{lemma:defloc} we also deduce that 
for all $j$ the spherical root $\sigma_j$ is not the sum of two orthogonal simple roots, 
and that if $\sigma_j$ is a simple root then at least one of the colors it moves is not free in $\s_I$, 
if $j\in I$.

Notice also that if $\sigma_j$ is a simple root then it cannot move two non-free colors in $\s$, 
otherwise $\s_{\hat j}$ would have defect equal to $\defect(\s)+2$, which contradicts Definition~\ref{def:combqhd}, part (1).

Therefore the map $\varphi_I$ can be defined as follows: 
a color $D$ of $\s_I$ moved by $\alpha\in S\setminus\Sigma_I$ is mapped to the unique free color of $\s$ moved by $\alpha$. 
The uniqueness of $\varphi_I$ 
and the equalities $\varphi_I(\Delta''_I)=\Delta''$ and $\varphi_I(\Delta'_I)=\Delta'$ are obvious, 
it remains to show the compatibility of $\varphi_I$ with the Cartan pairing.

If $\varphi_I(D)\in\Delta^{b}$ then $D\in(\Delta_I)^{b}$, 
hence both $c(\varphi_I(D),-)$ and $c_I(D,-)$ take the same values of the coroot of a simple root moving $D$.

If $\varphi_I(D)\in\Delta^{2a}$, and $\alpha$ is the (unique) simpe root moving it (so $2\alpha\in\Sigma$), 
both $c(\varphi_I(D),-)$ and $c_I(D,-)$ are equal to $\frac12\alpha^{\vee}$, 
since no spherical root $\sigma_j$ can be the double of a simple root. 
Indeed, if $2\alpha=\sigma_j$, only the color moved by $\alpha$ would take positive values on $\sigma_j$, 
and the latter is not in $\Sigma/\Delta'$, then the color would be in $\Delta'$. 
Hence it would correspond to a color in $\Delta'_{\hat j}$, non-positive on any spherical root of $\s_{\hat j}$. 
This contradicts the minimality of $\Delta'_{\hat j}$. 

If $\varphi_I(D)\in\A$ and $D$ is in $\A_I$, then the compatibility with the Cartan pairing is true thanks to the definition of localization.

It remains the case where $\varphi_I(D)\in\A$ but $D\notin\A_I$. Then $\varphi_I(D)$ is the free color moved by some simple root $\sigma$ in $\{\sigma_i:i\not\in I\}$. Therefore $\sigma$ is not in $S^{p}_I$ and not in the support of $\Sigma_I$: from the minimality of $\Delta''_I$ we deduce that $D\notin\Delta''_I$, and the proof of (\ref{lemma:combqhd:bij:comp1}) is complete.

Let now $D'$ be the non-free color moved by $\sigma$. Since $\sigma$ is not a spherical root of $\s/\Delta'$, we deduce that $D'\in\Delta'$. We also know that $D\in(\Delta_I)^{b}$, hence $c_I(D,-)$ coincides with $\sigma^{\vee}$ on $\Sigma_I$. The same holds for the sum $c(\varphi_I(D),-)+c(D',-)$, since in $\s$ they are both moved by the simple root $\sigma$.
\end{proof}

\begin{lemma}\label{lemma:combqhd2}\
\begin{enumerate}
\item\label{lemma:combqhd2:free} If $\sigma\in\{\sigma_0,\ldots,\sigma_k\}$ is a simple root 
then exactly one of the colors of $\Delta(\sigma)$ is free, and it doesn't belong to $\Delta'$.
\item\label{lemma:combqhd2:colors}Either 
\begin{itemize}
\item[(a)]\label{lemma:combqhd2:colors:pos} all $\sigma_0,\ldots,\sigma_k$ are simple roots, 
we have $\Delta(\sigma_0)\cap\ldots\cap\Delta(\sigma_k) = \{ E \}$ for some color $E\in\Delta$, 
and for all simple root $\alpha$ we have $E\in\Delta(\alpha)$ only if $\alpha\in\{\sigma_0,\ldots,\sigma_k\}$,
\end{itemize}
or
\begin{itemize}
\item[(b)]\label{lemma:combqhd2:colors:npos} if $\sigma\in\{\sigma_0,\ldots,\sigma_k\}$ is a simple root 
then there exists a simple root $\alpha\notin\{\sigma_0,\ldots,\sigma_k\}$
such that $\Delta(\sigma)\cap \Delta(\alpha)\neq \varnothing$.
\end{itemize}
\end{enumerate}
\end{lemma}

\begin{proof}
To prove part (\ref{lemma:combqhd2:free}) let us fix $\sigma=\sigma_i$.
The first statement stems from the proof of Lemma~\ref{lemma:combqhd}. 
For the second statement, let $j\in\{0,\ldots,k\}$ be such that $i\neq j$. 
In $\s_j$ the simple root $\sigma_i$ moves a color $D$ which is non-positive on all spherical roots. 
If the free color moved by $\sigma_i$ in $\s$ is in $\Delta'$ 
then the set $\Delta_j'\setminus\{D\}$ is properly contained in $\Delta_j'$ and distinguished, 
hence contradicting the minimality of $\Delta_j'$.

We prove part (\ref{lemma:combqhd2:colors}). Suppose that there is some $\sigma_i$ which is a simple root 
and such that no color in $\Delta(\sigma_i)$ is moved by any simple root not in $S'=S\setminus\{\sigma_0,\ldots,\sigma_k\}$. 
We must then show that $\sigma_0,\ldots,\sigma_k$ are all simple roots, 
and the existence of a color $E$ with the required properties.

Suppose that $\sigma_j$ is not a simple root for some $j\neq i$, 
and consider $\s_{\{i,j\}}$. 
Then $\s_{\{i,j\}}$ and $\Delta'_{\{i,j\}}$ satisfy Definition~\ref{def:combqhd} 
with difference of defects equal to $1$ (the required spherical roots are namely $\sigma_i$, $\sigma_j$). 
In $\s_{\{i,j\}}$ both colors moved by $\sigma_i$ are free, 
which contradicts part (\ref{lemma:combqhd2:free}) of the present proposition. 
We conclude that all spherical roots $\sigma_0,\ldots,\sigma_k$ are simple roots. 
For all $i$ let us denote by $D^{+}_i$ the non-free color in $\Delta(\sigma_i)$.

We choose any $i\in\{0,\ldots,k\}$ and we show that $E = D^+_i$ belongs to $\Delta(\sigma_k)$ for all $k$. The color $E$ is non-free, moved by no simple root in $S'$, so $E$ is also equal to $D^+_j$, for some $j\neq i$. If $k=1$ then $\{i,j\} = \{0,1\}$ and we are done. Otherwise, consider $\s_{\hat i}$. 
The set $\Delta'_{\hat i}$ gives a quotient satisfying Definition~\ref{def:combqhd}, and with difference of defects equal to $k-1$.

The colors moved by $\sigma_j$ are again not moved by any root in $S'$, so by induction $\sigma_0,\ldots,\sigma_{i-1},\sigma_{i+1},\ldots,\sigma_k$ all move one color in common (in $\s_{\hat i}$). This color must be the restriction of $E$, otherwise $\sigma_j$ would move two non-free colors, and the proof is complete.
\end{proof}

We make now the crucial assumption that $\s/\Delta'$ and $\s_{i}$ (for all $i$) are geometrically realizable. To simplify the proofs of the next results of this section, we also add the mild assumption that $\s$ and $\s/\Delta'$ are spherically closed.

Let us denote by $K$ and $H_i$ the subgroups of $G$ corresponding to $\s/\Delta'$ and $\s_i$.  
For all $i$, we can suppose that $H_i\subset K$, and we choose corresponding Levi subgroups $L_i\subseteq L_K$.

Recall that $\s/\Delta'$ is a quotient of type ($\mathscr L$) of $\s$, therefore there exists set of colors $\Delta''$ of $\s$ as in Definition~\ref{def:combL}. In particular $\Delta''$ is parabolic and contains $\Delta'$.

Let us denote by $Q_-$ the parabolic subgroup of $G$ (containing $B_-$ and $K$) corresponding to $\Delta''$, and choose a Levi subgroup $L_{Q_-}$ containing $L_K$. 
Thanks to Lemma~\ref{lemma:combqhd}, part (\ref{lemma:combqhd:quot}), for all $i\in\{0,\ldots,k\}$ the group $Q_-$ is parabolic minimal containing $H_i$, the groups $L_i$ and $L_K$ differ only by the connected center and are both very reductive in $L_{Q_-}$, and moreover $H^{u}_i \subseteq K^{u}$. Since $\defect(\s_i) = \defect(\s/\Delta')$, we also have $L_i=L_K$.

\begin{lemma}\label{lemma:Hiprime}
For all $i,j\in\{0,\ldots,k\}$ the following statements hold.
\begin{enumerate}
\item \label{lemma:Hiprime:augm} The spherical $G$-system $\s_i$ is spherically closed; defining $\Xi_i' = \Z\Sigma$ and $\rho_i'=\rho$, the pair $
(\Xi'_i, \rho_i')$ is an augmentation of $\s_i$.
\item \label{lemma:Hiprime:unipotent} Denoting by $H'_i\subseteq H_i$ the spherical subgroup of $G$ 
given by the augmentation $(\Xi_i,\rho'_i)$ of $\s_i$ as in Proposition~\ref{prop:augm}, we have
\[
(H_i')^u = H_i^u.
\]
\item \label{lemma:Hiprime:Levi}  Let $L_i'\subseteq L_K$ be a Levi part of $H_i'$ contained in $L_K$. Then $L_i'$ and $L_K$ differ only by the connected center. Moreover, we have
\[
L_i' = L_j'.
\]
\item \label{lemma:Hiprime:isom} The spherical $L_{Q_-}$-varieties $Q_-/H_i'$ and $Q_-/H_j'$ are equivariantly isomorphic.
\end{enumerate}
\end{lemma}
\begin{proof}
Since $\s$ is spherically closed its localization $\s_i$ has this property too. The combinatorial properties required by $(\Z\Sigma,\rho)$ to be an augmentation are just some of the conditions defining a spherical $G$-system applied to $\s$, whence part (\ref{lemma:Hiprime:augm}).

Part (\ref{lemma:Hiprime:unipotent}) holds because $H_i/H_i'$ is a diagonalizable group. Moreover, the quotient
\[
\frac{L_K}{L'_i Z(L_K)^{\circ}}
\]
has finite character group, since the rank of $\weights(Z(L_K)^{\circ})$ is equal to the rank of $\weights(L_K)$. On the other hand the group $\Z\Sigma/\Z\Sigma_i$ has no torsion, which implies that $H_i/H_i'$ (and thus $L_K/L_i'$) is also connected: we deduce that $L_K = L'_i Z(L_K)^{\circ}$, and the first statement of part (\ref{lemma:Hiprime:Levi}) follows.

We prove the second statement of part (\ref{lemma:Hiprime:Levi}). The pull-back along the projection $G/H_i'\to G/H_i$ induces a bijection between the sets of colors of $G/H_i$ and $G/H'_i$, so we denote both by $\Delta_i$. We also consider the sets $\Delta'_i$ and $\Delta''_i$ as sets of colors of both $G/H_i'$ and $G/H_i$.

Then we apply once again the theory of augmentations, in particular its uniqueness part, applied to the spherical homogeneous space $G/K$ and the subgroups $L_i'K^{u}$ of $K$. The lattice of weights associated to $G/L_i'K^{u}$ is the lattice $\rho_i'(\Delta'_i)^{\perp}$ inside $\Xi_i'$, and its set of colors is the set of colors of $\s_i/\Delta'_i$, with the Cartan pairing induced by $\rho'_i$ on the lattice $\rho_i'(\Delta'_i)^{\perp}$.

We claim that these invariants are equal for all $i$. 
To show the claim, we observe that the compatibility relations of Lemma~\ref{lemma:combqhd} part (\ref{lemma:combqhd:bij}) also hold for the ``extended'' Cartan pairings $\rho_i'$, with the same proof.

Then $\rho_i'(\Delta'_i)^{\perp}$ is equal to $\rho_j'(\Delta'_j)^{\perp}$ as subsets of the lattice $\Xi_i'=\Xi'_j=\Z\Sigma$ thanks to parts (\ref{lemma:combqhd:bij:images}) and (\ref{lemma:combqhd:bij:comp1}) of Lemma~\ref{lemma:combqhd}. The set of colors is also the same for all $i$, since it is identified with the set of colors $\s_i/\Delta'_i$, that is also identified with the set of colors of $G/K$ thanks to part (\ref{lemma:combqhd:quot}) of Lemma~\ref{lemma:combqhd}, and part (\ref{lemma:combqhd:bij:comp2}) of the same lemma assures that the Cartan pairings restricted to $\rho_i'(\Delta'_i)^{\perp}$ are all the same. It follows that $L_i'K^{u}=L_j'K^{u}$ for all $i,j$.

To show part (\ref{lemma:Hiprime:isom}), we use Lemma~\ref{lemma:lambda} and \cite{Lo06}. Thanks to Lemma~\ref{lemma:combqhd} parts (\ref{lemma:combqhd:values}) and (\ref{lemma:combqhd:bij}) applied to the ``extended'' Cartan pairings $\rho_i'$ we have
\begin{equation}\label{eqn:deltasec}
\Q_{\geq0}\rho_i'(\Delta''_i) = \Q_{\geq0}\rho_j'(\Delta''_j).
\end{equation}

Then we notice that Lemma~\ref{lemma:lambda} holds even under the weaker hypothesis that $H$ (of \loccit) be a spherical, not necessarily wonderful subgroup of $G$, with the same proof. Then it may be applied to the group $H_i'$ for all $i$: together with the equalities (\ref{eqn:deltasec}) and $\Xi'_i=\Xi'_j$ for all $i,j$, it implies that $\C[Q_-/H_i']\cong \C[Q_-/H_j']$ as $L_{Q_-}$-modules.

Thanks to \cite{Lo06}, the smooth affine spherical $L_{Q_-}$-varieties $Q_-/H_i'$ for $i$ ranging from $0$ to $k$ are all $L_{Q_-}$-equivariantly isomorphic, which shows part (\ref{lemma:Hiprime:isom}).
\end{proof}
\begin{corollary}\label{cor:isomquot}
Denote by $L'$ the group $L'_0=\ldots=L'_k$. The quotient $\Lie\,K^u/\Lie\,H^u_i$ is an irreducible $L'$-module for all $i\in\{0,\ldots,k\}$. Up to replacing if necessary the groups $H_0,\ldots,H_k$ with conjugates also inside $K$, we have
\[
\Lie\,K^u/\Lie\,H^u_i \cong \Lie\,K^u/\Lie\,H^u_j
\]
as $L'$-modules for all $i,j\in\{0,\ldots,k\}$.
\end{corollary}
\begin{proof}
The quotients are irreducible $L_K$-modules, and also irreducible under the action of $L'$ since $Z(L_K)^{\circ}L' = L_K$.

From part (\ref{lemma:Hiprime:isom}) of Lemma~\ref{lemma:Hiprime} and the standard isomorphism
\[
Q_-/H_i' \cong L_{Q_-}\times_{L'} \Lie\,Q^u_-/\Lie\,H^u_i
\]
as $L_{Q_-}$-varieties, we deduce that the quotients  $\Lie\,Q^u_-/\Lie\,H^u_i$ and $\Lie\,Q^u_-/\Lie\,H^u_j$ are $L'$-equivariantly isomorphic, up to twisting the action of $L'$ on $\Lie\,Q^u_-/\Lie\,H^u_j$ with an element in $\N_{L_{Q_-}}L'$ (see \cite[Lemma~3.6.7]{Lo07}).

Now it is possible to apply \cite[Proof of Theorem 1, step 5]{Lo07}, 
with $H_1$, $H_2$, $\widetilde H$ of \loccit\ equal resp.\ to $H'_i$, $H'_j$ and $K$. 
The proof can be carried out in the same way as in \loccit\ even though $H'_i$ and $H'_j$ are not equal to their normalizers. This assures the existence of an element $g\in\N_GK$ such that the groups $H_i$ and $gH_jg^{-1}$ satisfy the thesis.
\end{proof}

\begin{proposition}\label{prop:higher}
Up to replacing the groups $H_0,\ldots,H_k$ with conjugates inside $K$, there exist an $L_K$-module decomposition $\Lie\,K^u=V\oplus W_0\oplus\ldots\oplus W_k$ and $Z(L_K)^\circ$-characters $\gamma_0,\ldots,\gamma_k$ such that for all $i\in \{0,\ldots,k\}$
\begin{enumerate}
\item\label{prop:higher:iso} $W_i$ is a simple module under the action of $L'$;
\item\label{prop:higher:weight} $Z(L_K)^\circ$ acts on $W_i$ via the weight $\gamma_i$;
\item\label{prop:higher:sum} $\Lie\,H^u_i = V\oplus W_0\oplus\ldots W_{i-1}\oplus W_{i+1}\oplus\ldots\oplus W_k$.
\end{enumerate}
\end{proposition}
\begin{proof}
The proposition stems from Lemma~\ref{lemma:Hiprime} and Corollary~\ref{cor:isomquot}.
\end{proof}

We are ready to state the main result of this section.

\begin{theorem}\label{thm:qhd}
Let $\s$ be a spherically closed spherical $G$-system with a minimal quotient of higher defect $\s/\Delta'$ with
$k=\defect(\s/\Delta')-\defect(\s)$. 
Assume that the spherical $G$-systems $\s/\Delta'$, $\s_0,\ldots,\s_k$ are geometrically realizable, and denote the corresponding wonderful subgroups of $G$ resp.\ by $K$, $H_0,\ldots,H_k$. 
If $K$ is spherically closed then $\s$ is geometrically realizable. Choosing the subgroups $K$, $H_0,\ldots,H_k$ in such a way that Proposition~\ref{prop:higher} holds, then the weights $\gamma_0-\gamma_1,\gamma_0-\gamma_2\ldots,\gamma_0-\gamma_k$ are linearly independent, and there exists a generic stabilizer $H$ of $\s$ such that:
\begin{enumerate}
\item\label{thm:qhd:l} its Levi part $L$ differs from $L_K$ only by the connected center;
\item\label{thm:qhd:c} the connected center $Z(L)^\circ$ is the connected subgroup of $Z(L_K)^\circ$ defined by the equations $\gamma_0=\ldots=\gamma_k$;
\item\label{thm:qhd:u} $\Lie\,H^u$ is a co-simple $L$-submodule of $\Lie\, K^u$ containing $V$ but not any direct summand $W_0,\ldots,W_k$ of $\Lie\, K^u$.
\end{enumerate}
\end{theorem}
\begin{proof}
The unipotent radical $H_i^u$ contains $(K^u,K^u)$ for all $i$, as in the proof of Proposition~\ref{prop:L}. 
Hence there exists a subgroup $H^u$ of $K^u$ such that $\Lie\,H^u$ is a co-simple $L'$-submodule of $\Lie\,K^u$ containing $V$ but not any direct summand $W_0,\ldots,W_k$ of $\Lie\, K^u$.

We define $L=\N_{L_K} H^u$, and claim that $L_K/L$ is connected. Indeed, the quotient $L_K/L'$ is connected by Lemma~\ref{lemma:Hiprime} and we have $L\supseteq L'$, whence the claim.

With this definition the groups $L$, $Z(L)^\circ$ and $H^{u}$ satisfy properties (\ref{thm:qhd:l}), (\ref{thm:qhd:c}) and (\ref{thm:qhd:u}). In addition $L$ is very reductive in $L_{Q_-}$, therefore $L$ and $H^{u}$ are resp.\ a Levi subgroup and the unipotent radical of the product $LH^u$. We set then $H=LH^{u}$, and it remains to show that $H$ is wonderful, with spherical system $\s$.

For the first claim we may apply \cite[\S2.4.2, Corollary 3]{BL11}, if we show that $\N_KH=H$ and that $H$ is spherical. The former equality is true by construction, since $H$ is also equal to $\N_KH^u$. Let us show the latter property. The inclusion $L\supseteq L'$ induces a surjective $L_{Q_-}$-equivariant map $Q_-/(L' H^{u}) \to Q_-/H$. Both these varieties are affine, and $Q_-/(L'H^{u})$ is also $L_{Q_-}$-equivariantly isomorphic to $Q_-/H'_i$.  Since the groups $H'_i$ are spherical, we deduce that the quotient $Q_-/H$ is a spherical variety with respect to the action of $L_{Q_-}$. This implies that $H$ is spherical.

Moreover, from the proof of Lemma~\ref{lemma:lambda} we deduce that the rank (as a spherical $G$-variety) of $G/(L'H^{u})$ is equal to the rank (as a spherical $L_{Q_-}$-variety) of $Q_-/(L'H^{u})$. Since the same holds for $G/H_i'$ and $Q_-/H_i'$, we deduce that the rank of $G/(L'H^{u})$ is equal to the rank of $G/H_i'$, which is $|\Sigma|$ by definition of $H_i'$. Therefore the rank of $G/H$ is not greater than $|\Sigma|$.

Let now $X$ be the wonderful embedding of $G/H$, with spherical system $\s_X$: we claim that $\s_X=\s$.  We write $\s_X = (S^p_X, \Sigma_X, \A_X)$ and $\s = (S^p, \Sigma, \A)$.

We assume for a while that $\gamma_0-\gamma_1,\ldots,\gamma_0-\gamma_k$ are linearly independent. Then for all $i\in\{0,\ldots,k\}$ there exists a 1-PSG $\lambda_i$ of $Z(L_K)^\circ$ such that
\[
\lim_{t\to 0} \lambda_i(t)\Lie\,H^{u} = \Lie\,H_i^{u}
\]
in the Gra\ss mannian of $\dim\Lie\,H^{u}$-dimensional subspaces of $\Lie\,K^{u}$. If $x\in X$ is the point $eH \in G/H$, then consider the limit
\[
x(i) = \lim_{t\to 0} \lambda_i(t)\,x
\]
in $X$ and its stabilizer $H(i)\subseteq G$. The projection $G/H\to G/K$ extends to a $G$-equivariant map from $X$ to the wonderful embedding of $G/K$. It sends $x$ to the point $eK \in G/K$, which is $\lambda_i(\C^*)$-stable, therefore $H(i)\subseteq K$.

On the other hand $\lambda_i(t)x$ is $L$-stable for any $t$, hence $L\subseteq H(i)$, and we also have by construction $H_i^{u}\subseteq H(i)$. Since $L$ is very reductive in $L_{Q_-}$ we have that $H(i)^{u}\subseteq Q_-^{u}$, hence $H(i)^{u}\subseteq K^{u}$.

Then we claim that $H(i)^{u}=H_i^{u}$. The equality stems from the fact that the $G$-equivariant natural projection $G/H\to G/\N_GH$ can be extended to a finite $G$-equivariant map $\delta\colon X\to Y$, where $Y$ is the {\em Demazure embedding} of $G/\N_GH$ inside the Gra\ss mannian of ($\dim\Lie\, H$)-dimensional subspaces of $\Lie\,G$ (see \cite{Br90}).

The Lie algebra of $(LH_i^{u})^{\circ}$ is the limit, in the Gra\ss mannian, of the Lie algebra of the stabilizer of $\lambda_i(t)x$. Hence the stabilizer in $G$ of $\delta(x(i))$ has connected component equal to $(LH_i^{u})^{\circ}$. Now, thanks to the definition of the Demazure embedding itself, the groups $H(i)^{\circ}$ and $(LH_i^{u})^{\circ}$ have the same unipotent radical, which implies $H(i)^{u}=H_i^{u}$.

This also implies that $Z(L_K)^{\circ}$ normalizes $H(i)$: since $H(i)$ is wonderful, thus with finite index in its normalizer, and since $Z(L_K)^{\circ}L'=L_K$, we conclude that $H(i)$ contains $L_K$. Therefore $H(i)=L_KH_i^{u}=H_i$.

To sum up, $\s_X$ has $k+1$ localizations with spherical systems $\s_0,\ldots,\s_k$, and a quotient $\s_X/\Delta'_X$ equal to $\s/\Delta'$, due to the inclusion $H\subseteq K$. It follows that the set of spherical roots of $\s_X$ contains $\Sigma$. Since the rank of $X$ is $\leq|\Sigma|$, we have $\Sigma_X=\Sigma$. 
Moreover $\s_X$ and $\Delta'_X$ satisfy parts (2) and (3) of Definition~\ref{def:combqhd}, thanks to the properties of $H$.

We now prove that the weights $\gamma_0,\ldots,\gamma_k$ are linearly independent.
Let $I$ be any nonempty subset of $\{0,\ldots,k\}$ and recall the notation
\[
\Sigma_I = \Sigma \setminus\{ \sigma_j \;|\; j\in\{0,\ldots,k\}\setminus I\} .
\]

Define $H_I^{u}$ such that $\Lie\,H^{u}_I$ is a co-simple $L'$-submodule of $K^{u}$ containing $V$ and $W_j$ for all $j\notin I$ but not containing $W_j$ for any $j \in I$, and define the groups $L_I$ and $H_I$ with the same above procedure for $L$ and $H$. With this notation $H_{\{i\}}=H_i$ for all $i$.

Choose any $i_0\in I$. If the weights $\gamma_{i_0}-\gamma_i$ for $i\in I$ are linearly independent, the same arguments used above for $H$ show that $H_I$ is a wonderful subgroup of $G$, that the corresponding wonderful variety has a localization equal to $\s_i$ for each $i\in I$, and that its set of spherical roots is $\Sigma_I$.

Suppose now that the weights $\gamma_0-\gamma_1,\ldots,\gamma_0-\gamma_k$ are not all linearly independent, but not all equal to $0$. Then there exist two different subsets $I$, $J$ of $\{0,\ldots,k\}$ such that
\begin{enumerate}
\item for some $i_0\in I$ and some $j_0\in J$ the sets $\{\gamma_{i_0}-\gamma_i\;|\; i\in I, i\neq i_0\}$ and $\{\gamma_{j_0}-\gamma_j\;|\; j\in J, j\neq j_0\}$ are both linearly independent sets of weights;
\item the two above sets span the same subspace of $\weights(Z(L_K)^{\circ})\otimes_\Z\Q$. 
\end{enumerate}
This implies that $L_I=L_J$, and that the $L_I$-module $\Lie\,K^{u}/\Lie\,H_I$ is isomorphic to $\Lie\,K^{u}/\Lie\,H_J$. Then we may apply \cite[Proof of Theorem 1, step 6]{Lo07} with $H_1$ and $H_2$ of \loccit\ equal resp.\ to $H_I$ and $H_J$. This can be done, in particular \cite[Proposition 3.3.1]{Lo07} applies to the Lie algebras of $H_I$ and $H_J$, since these groups have finite index in their respective normalizers.

We deduce that $H_I$ and $H_J$ are conjugated in $G$, which is absurd because the respective wonderful varieties have different sets of spherical roots.

If the weights $\gamma_0,\ldots,\gamma_k$ are all equal, the same proof applies to the groups $H_i$ and $H_j$ instead of $H_I$ and $H_J$, for any $i,j\in\{0,\ldots,k\}$ with $i\neq j$. We conclude that the weights $\gamma_0-\gamma_1,\ldots,\gamma_0-\gamma_k$ are linearly independent.

We prove now that $\s_X$ and $\Delta'_X$ satisfy part (1) of Definition~\ref{def:combqhd}. 
Consider again a nonempty subset $I$ of $\{0,\ldots,k\}$, and let $\lambda_I$ be a $1$-PSG of $Z(L_K)^\circ$ such that
\[
\lim_{t\to 0} \lambda_I(t)\Lie\,H^{u} = \Lie\,H_I^{u}
\]
The stabilizer $H(I)$ of the point
\[
x(I) = \lim_{t\to 0} \lambda_i(t)\,x
\]
satisfies $H(I)^{u}=H_I^{u}$, thanks to a proof similar to that of the equality $H(i)^{u}=H_i^{u}$. Since $H(I)$ has finite index in its normalizer, we also deduce that $H(I)$ has finite index in $H_I=\N_KH_I^{u}$. Now the spherical roots of $G/H_I$ are a subset of those of $X$, which implies that $H(I)=H_I$. It follows that $H(I)$ is a generic stabilizer of the localization of $X$ on the set of spherical roots $\Sigma_I$.

By construction, and since the characters $\gamma_0-\gamma_1,\ldots,\gamma_0-\gamma_k$ are linearly independent, the rank of $\weights(H_I)$ is $\rank\,\weights(K) - |I| + 1$, hence the defect of $G/H(I)$ is $\defect(\s)+|\Sigma|-|\Sigma_I|$. 
This shows that $\s_X$ and $\Delta'_X$ satisfy part (1) of Definition~\ref{def:combqhd}.

We are now ready to show the equality $\s_X=\s$. The sets $\Sigma_X$ and $\Sigma$ are equal, and since $\s_X$ and $\s$ have equal localizations on some sets of spherical roots (e.g.\ $\Sigma_i$ for some $i$) we also have $S^{p}_X = S^{p}$. It remains to compare $\A_X$ and $\A$.

Define
\[
\Sigma^{s} = \Sigma \setminus (\{\sigma_0,\ldots,\sigma_k\} \setminus S).
\]
Let $I$ be a non-empty subset of $\{0,\ldots,k\}$ such that $\Sigma_I\subseteq\Sigma^{s}$. We first show by induction on $|I|$ that the localizations $(\s_X)_{\Sigma_I}$ and $\s_{\Sigma_I}$ are equal. 
Notice that both satisfy Definition~\ref{def:combqhd} hence Lemma~\ref{lemma:combqhd2} with respect to the sets of colors resp.\ $(\Delta'_X)|_{\Sigma_I}$ and $\Delta'|_{\Sigma_I}$. The case $|I|=1$ is true, we show the induction step.

We denote by $(\Delta_X)_{\Sigma_I}$ and by $\Delta_{\Sigma_I}$ the sets of colors of resp.\ $(\s_X)_{\Sigma_I}$ and $\s_{\Sigma_I}$. Consider an element $i\in I$ and let $J = I\setminus\{i\}$. By the induction hypothesis the localizations $(\s_X)_{\Sigma_J}$ and $\s_{\Sigma_J}$ are equal. 
Moreover, Lemma~\ref{lemma:combqhd} provides certain bijections
\[
(\Delta_X)_{\Sigma_I} \leftarrow \Delta_J\to \Delta_{\Sigma_I}
\]
Inverting the first one and composing we obtain a bijection $e_I\colon(\Delta_X)_{\Sigma_I} \to \Delta_{\Sigma_I}$ which is compatible with the Cartan pairings computed on spherical roots in $\Sigma_J$ (i.e.\ all except $\sigma_i$), thanks to Lemma~\ref{lemma:combqhd}, part (\ref{lemma:combqhd:bij:comp1}), and the proof of part (\ref{lemma:combqhd:bij:comp2}).

To complete the induction we must check that the map $e_I$ is also compatible with the Cartan pairing computed on the spherical root $\sigma_i$. Let then $D\in(\Delta_X)_{\Sigma_I}$. If $D$ is moved by some simple root $\alpha\notin \Sigma_I$ then its values on spherical roots are equal to those of $\alpha^{\vee}$, and the same holds for $e_I(D)$: the compatibility of $e_I$ with the Cartan pairing follows. Suppose now that $D$ is moved by some simple root in $\Sigma_I$, and it is moved by no simple root outside $\Sigma_I$. Thanks to Lemma~\ref{lemma:combqhd2}, either $D$ is moved by some simple root in $\Sigma\setminus\{\sigma_0,\ldots,\sigma_k\}$, or $D$ is moved by $\sigma_t$ for some $t\in\{0,\ldots,k\}$ and case (\ref{lemma:combqhd2:colors:pos}) of Lemma~\ref{lemma:combqhd2} occurs.

In the first case the Cartan pairings of $\sigma_i$ with $D$ in $(\s_X)_{\Sigma_I}$ and  with $e_I(D)$ in $\s_{\Sigma_I}$ are equal thanks to the fact that both $\s_X$ and $\s$ localize to $\s_i$. In the second case, consider the color $E$ of Lemma~\ref{lemma:combqhd2} applied to $(\s_X)_{\Sigma_I}$. Then $e_I(E)$ is the color of Lemma~\ref{lemma:combqhd2} applied to $\s_{\Sigma_I}$, and both $E$ and $e_I(E)$ take value $1$ on $\sigma_i$. Therefore the Cartan pairings of both $D$ and $e_I(D)$ with $\sigma_i$ are equal to $\langle\alpha_i,\alpha_t^{\vee}\rangle - 1$.

This completes the proof of $(\s_X)_{\Sigma_I}=\s_{\Sigma_I}$ for all $I$, which also shows in particular $(\s_X)_{\Sigma^{s}}=\s_{\Sigma^{s}}$. Finally, we observe that $(\A_X)_{\Sigma^{s}} =  \A_X$ and that $(\A)_{\Sigma^{s}} =  \A$, which provides a bijection $e\colon \A_X\to \A$. Since $(\s_X)_{\Sigma^{s}}=\s_{\Sigma^{s}}$, we only have to show that $D\in \A_X$ takes the same value as $e(D)$ on any spherical root $\sigma_i\in\Sigma\setminus\Sigma^{s}$.

If case (\ref{lemma:combqhd2:colors:npos}) of Lemma~\ref{lemma:combqhd2} occurs, then this is true because both $\s_X$ and $\s$ localize to $\s_i$. If case (\ref{lemma:combqhd2:colors:pos}) of Lemma~\ref{lemma:combqhd2} occurs then $\Sigma=\Sigma^{s}$, and the proof is complete.
\end{proof}

\subsection{Quotients of higher defect given by positive colors}
A particular case of quotients of higher defect has already been studied in \cite{Lu01} for $G$ of type $\mathsf A$. It is related to a special type of distinguished sets: those containing only one element.

\begin{definition}
Let $\s=(S^p,\Sigma,\A)$ be a spherical $G$-system. 
A color $D$ is {\em positive} if $c(D,\sigma)\geq 0$ for all $\sigma\in\Sigma$. 
In this case we define $\Sigma_D$ to be the set of spherical roots $\sigma\in\Sigma$ such that $c(D,\sigma)>0$.
\end{definition}

It is obvious that $\Delta'=\{D\}$ is good distinguished and smooth whenever $D$ is a positive color. 
Indeed, $\Sigma/\{D\}=\Sigma\setminus\Sigma_D$. Since the colors of the quotient are $\Delta\setminus\{D\}$, we have
\[
\defect(\s/\{D\})-\defect(\s) = |\Sigma_D|-1.
\]
Therefore $\defect(\s/\{D\})=\defect(\s)$ if and only if $D$ is free.

If $D$ is a non-free positive color of a spherical system $\s$, each localization $\s_i$, 
uniquely determined by $\{D\}$ as above, contains a free positive color equal to the restriction of $D$.

As a consequence, the geometric realizability of spherical systems with non-free positive colors follows 
from the geometric realizability of all spherical systems with free positive colors.

Notice that here the assumption that the wonderful subgroup $K$ of $G$ associated with $\s/\{D\}$ is spherically closed 
can be dropped. Indeed, a loose spherical root of type $\mathsf B$ or $\mathsf G$ belongs to $\Sigma/\{D\}$ 
if and only if it belongs to $\Sigma$.

\section{Tails}\label{s:tails}

\subsection{Spherical systems with tails}\label{s:deftails}
%Another general procedure applies to indecomposable spherical $G$-systems $\s$ having a subset of spherical roots $\tau$ of special type. In these cases, we say that $\s$ is obtained from its localization in a certain subset of simple roots $S(\tau)\supseteq S\setminus\supp\,\tau$ by {\em adding the tail} $\tau$.

%If $S'\cap\supp\,\gamma=\emptyset$ we denote by $\{\alpha_1,\ldots,\alpha_p\}$ (numbered from left to right) the simple roots of the connected component of $S'$ where $\gamma$ is glued. There exists at most one spherical root $\gamma'$ non-orthogonal to $\supp\,\gamma$, and $\gamma$ is one of the following:

\begin{definition}\label{def:tailroots}
Let $\s=(S^p,\Sigma,\A)$ be a spherical $G$-system and $\tau=\{\gamma\}\subseteq\Sigma$. 
Suppose that $\supp\,\tau$ lies on a single connected component $\overline{\supp\,\tau}=\{\alpha_1,\ldots,\alpha_n\}$ 
of the Dynkin diagram of $G$. Then $\tau$ is called a {\em tail root} if one of the following cases occurs.
\begin{description}
\item[$\mathbf {b(m)}$] ($2\leq m < n$)
\[
\begin{picture}(7500,600)(-300,-300)
\put(0,0){\usebox{\shortbm}}
\end{picture}
\]
$\gamma = \alpha_{n-m+1}+\ldots+\alpha_{n}$ and $\overline{\supp\,\tau}$ of type $\mathsf B_n$.
\item[$\mathbf {b(1)}$]
\[
\begin{picture}(600,1800)(-300,-900)
\put(0,0){\usebox{\aone}}
\end{picture}
\]
$\gamma = \alpha_{n}$ and $\overline{\supp\,\tau}$ of type $\mathsf B_n$. The following must also hold: $\rho(D_n^+) = \rho(D_n^-)$.
\item[$\mathbf {2b(m)}$] ($1\leq m < n$)
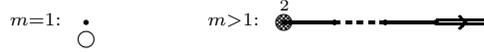
\[
\begin{picture}(15000,2000)(-300,-1000)
\put(2000,0){\usebox{\aprime}}
\put(-700,-100){$\scriptstyle{m=1:}$}
\put(9000,0){\usebox{\shortbprimem}}
\put(6300,-100){$\scriptstyle{m>1:}$}
\end{picture}
\]
$\gamma = 2\alpha_{n-m+1}+\ldots+2\alpha_{n}$ and $\overline{\supp\,\tau}$ of type $\mathsf B_n$.
\item[$\mathbf {c(m)}$] ($2\leq m \leq n$)
\[
\begin{picture}(9300,600)(-300,-300)
\put(0,0){\usebox{\shortcsecondm}}
\end{picture}
\]
$\gamma = \alpha_{n-m+1}+2\alpha_{n-m+2}+\ldots+2\alpha_{n-1}+\alpha_{n}$ with $\alpha_{n-m+1}\notin S^p$, and $\overline{\supp\,\tau}$ of type $\mathsf C_n$. In this case we will say that $\tau$ is {\em overlapping}.
\item[$\mathbf {d(m)}$] ($2\leq m < n$)
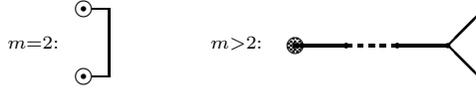
\[
\begin{picture}(14000,3000)(-1500,-300)
\multiput(0,0)(0,2400){2}{
\put(0,0){\circle*{150}}
\put(0,0){\usebox{\wcircle}}
\put(300,0){\line(1,0){600}}}
\put(900,0){\line(0,1){2400}}
\put(-2700,1000){$\scriptstyle{m=2:}$}
\put(7500,1100){\usebox{\shortdm}}
\put(4500,1000){$\scriptstyle{m>2:}$}
\end{picture}
\]
$\gamma =2\alpha_{n-m+1}+\ldots+2\alpha_{n-2}+\alpha_{n-1}+\alpha_{n}$ and $\overline{\supp\,\tau}$ of type $\mathsf D_n$.
\end{description}
\end{definition}

\begin{definition}\label{def:tailpairs}
Let $S=(S^p,\Sigma,\A)$ be a spherical $G$-system and $\tau \subseteq\Sigma$ with $|\tau|=2$. 
Suppose that $\supp\,\tau$ lies on a single connected component 
$\overline{\supp\,\tau}=\{\alpha_1,\ldots,\alpha_n\}$ of the Dynkin diagram of $G$. 
Then $\tau$ is called a {\em tail pair} if one of the following cases occurs.
\begin{description}
\item[$\mathbf {(aa,aa)}$] $\tau = \{ \alpha_1+\alpha_6, \alpha_3+\alpha_5 \}$ and $\overline{\supp\,\tau}$ of type $\mathsf E_6$.
\item[$\mathbf {(d3,d3)}$] $\tau = \{ \alpha_2+2\alpha_4+\alpha_5, \alpha_5+2\alpha_6+\alpha_7 \}$ and $\overline{\supp\,\tau}$ of type $\mathsf E_7$.
\item[$\mathbf {(d5,d5)}$] $\tau = \{ 2\alpha_1+\alpha_2+2\alpha_3+2\alpha_4+\alpha_5, \alpha_2+\alpha_3+2\alpha_4+2\alpha_5+2\alpha_6 \}$ and $\overline{\supp\,\tau}$ of type $\mathsf E_8$.
\item[$\mathbf {(2a,2a)}$] $\tau = \{ 2\alpha_3, 2\alpha_4\}$ and $\overline{\supp\,\tau}$ of type $\mathsf F_4$.
\end{description}
\end{definition}

\begin{example}
Let us give some examples of spherical $G$-systems with a tail pair, each with localization on $S\setminus \supp\,\tau$ having Luna diagram as follows:
\[
\begin{picture}(3600,600)(-300,-300)
\put(0,0){\usebox{\atwo}}
\end{picture}
\]
The examples are:
\begin{description}
\item[$\mathbf {(aa,aa)}$]
\[
\begin{picture}(7800,3600)(-300,-2000)
\multiput(0,0)(1800,0){4}{\usebox{\edge}}
\put(3600,0){\usebox{\atwos}}
\multiput(0,0)(5400,0){2}{\multiput(0,0)(1800,0){2}{\circle{600}}}
\multiput(0,300)(7200,0){2}{\line(0,1){600}}
\put(0,900){\line(1,0){7200}}
\multiput(1800,300)(3600,0){2}{\line(0,1){300}}
\put(1800,600){\line(1,0){3600}}
\end{picture}
\]
\item[$\mathbf {(d3,d3)}$]
\[
\begin{picture}(9600,3600)(-300,-2000)
\put(0,0){\usebox{\dynkineseven}}
\put(0,0){\usebox{\atwo}}
\multiput(3600,0)(3600,0){2}{\usebox{\gcircle}}
\end{picture}
\]
\item[$\mathbf {(d5,d5)}$]
\[
\begin{picture}(11400,3600)(-300,-2000)
\put(0,0){\usebox{\dynkineeight}}
\put(9000,0){\usebox{\atwo}}
\multiput(0,0)(7200,0){2}{\usebox{\gcircle}}
\end{picture}
\]
\item[$\mathbf {(2a,2a)}$]
\[
\begin{picture}(6000,2000)(-300,-1000)
\put(0,0){\usebox{\dynkinffour}}
\put(0,0){\usebox{\atwo}}
\multiput(3600,0)(1800,0){2}{\usebox{\aprime}}
\end{picture}
\]
\end{description}
\end{example}

\begin{definition}\label{def:tails}
Let $\s=(S^p,\Sigma,\A)$ be a spherical $G$-system and $\tau\subseteq\Sigma$ a tail root or a tail pair. 
Set $S(\tau) = (S\setminus \supp\,\tau) \cup \{\alpha_{n-m+1}\}$ if $\tau$ is overlapping, 
and $S(\tau) = S\setminus \supp\,\tau$ otherwise. Then $\tau$ is a {\em tail} of $\s$ if the following conditions hold:
\begin{enumerate}
\item $\supp\,(\Sigma\setminus\tau) = S(\tau)$;
\item $\s$ has a smooth quotient $\s/\Delta(\tau)$ such that $\Sigma/\Delta(\tau)=\tau$.
\end{enumerate}
\end{definition}

In what follows we will always suppose that $\Delta(\tau)$ is minimal with respect to its properties. If $|\tau|=1$ (resp.\ $2$) then $\tau$ will also be called a {\em classical} (resp.\ {\em exceptional}) {\em tail}.

We underline that the assumptions on $\Delta(\tau)$ are needed for our discussion: there exist geometrically realizable examples of $\s$ where $\tau$ exists with the above properties except for the existence of $\Delta(\tau)$. In these cases, generic stabilizers do behave differently from cases where $\tau$ is a tail.

A set of colors $\Delta(\tau)$ exists always if $\supp\,\tau=\overline{\supp\,\tau}$. In this case indecomposability implies that $\tau$ is overlapping (we prove this fact at the beginning of \S\ref{s:tailno}), and $\Delta(\tau)$ can be taken to be all the colors of $\s_{S(\tau)}$. Otherwise, if some $\Delta(\tau)$ exists, it is easy to see that a (non-minimal) choice can be given by all the colors of $\s_{S(\tau)}$ except one.

Being of rank $1$ or $2$, the quotient $\s/\Delta(\tau)$ is geometrically realizable: we call its generic stabilizer $K^\tau$.

\begin{theorem}\label{thm:tails}
Let $\s$ be an indecomposable spherical system without minimal quotients of higher defect and with a tail $\tau$, 
and suppose that the localization $\s_{S(\tau)}$ is geometrically realizable with generic stabilizer $K_\tau$. 
Then $\s$ is geometrically realizable, and a generic stabilizer $H$ can be defined using $K_\tau$ and $K^\tau$ 
as described in \S\ref{s:tailno} (for non-overlapping tails) or \S\ref{s:tailC} (for overlapping tails).
\end{theorem}

We prove the theorem in \S\ref{s:tailsproof}.

%%%%% ......... If $S'\cap\supp\,\gamma\neq\emptyset$ we denote by $\{\alpha_1,\ldots,\alpha_{p+1}\}$ (numbered from left to right) the simple roots of the connected component of $S'$ where $\gamma$ is glued. Then $\gamma$ is of the form:

%Let us call $\beta_1,\beta_2$ (and $\beta_3$) resp.\ the first simple root of $\supp\,\gamma$ and the simple root(s) moving the colors whose functional is positive on $\gamma$:
%\begin{itemize}
%\item[-] $b(m)$, $2b(m)$, and $d(m)$ with $m\geq 3$: here $\beta_1=\beta_2=\beta_3$ is the first simple root of $\supp\,\gamma$;
%\item[-] $c(m)$: here $\beta_2=\beta_3$ is the second simple root of $\supp\,\gamma$;
%\item[-] $d(2)$: here $\beta_2$ and $\beta_3$ are the two simple roots of $\supp\,\gamma$, say also $\beta_1=\beta_2$.
%\end{itemize}

%Being of rank $1$ or $2$, the wonderful variety associated to $\s/\Delta''$ exists: we call its generic stabilizer $M\subset G$. If our strategy is based on an induction on the number of tails, we can suppose the same for the wonderful variety associated to $\s'$: we call its generic stabilizer $K\subset G_{S'}$. We will now see how to build a candidate subgroup $H$ for the spherical system $\s$ starting from $M$ and $K$.

\subsection{The wonderful subgroup associated with a spherical system with tails}\label{s:generictails}
In this section $\s$ is an indecomposable spherical $G$-system with a tail in $\tau\subseteq\Sigma$, 
and we suppose that the localization $\s_{S(\tau)}$ is geometrically realizable with generic stabilizer $K_\tau$. 
We use the groups $K_\tau$, $K^\tau$ and the system $\s$ to define a subgroup $H$ of $G$, 
a candidate for the wonderful subgroup of $G$ associated with $\s$.

\begin{example} Consider the following example of a tail of type $b(m)$, for $G= \SO(2n+1)$ (see \cite[Table 2, case 6]{W96}):
\[
\begin{picture}(14700,600)(-300,-300)
\put(0,0){\usebox{\mediumam}}
\put(5400,0){\usebox{\plusbm}}
\end{picture}
\]
\[
\s=(\{ \alpha_2,\ldots,\alpha_{l-1},\alpha_{l+2},\ldots,\alpha_{l+m} \}, \{ \alpha_1+\ldots+\alpha_l, \alpha_{l+1}+\ldots+\alpha_{l+m} \}, \emptyset)
\]
where $l=n-m$. 
The system $\s$ has rank $2$, its wonderful subgroup $H$ of $G$ has Levi factor $\GL(l)\times \SO(2m)$ 
and Lie algebra of the unipotent radical equal to $\C^l\otimes\C^{2m}\oplus\wedge^2\C^l$. 
The quotient $\s/\Delta(\tau)$ has rank $1$
\[
\begin{picture}(14700,600)(-300,-300)
\put(0,0){\usebox{\edge}}
\put(1800,0){\usebox{\shortsusp}}
\put(3600,0){\usebox{\edge}}
\put(5400,0){\usebox{\wcircle}}
\put(5400,0){\usebox{\plusbm}}
\end{picture}
\]
and its wonderful subgroup of $G$ is $K^\tau=(\GL(l)\times \SO(2m))Q^u$. 
It is a parabolic induction by means of the parabolic subgroup $Q = P_{- S\setminus\{ \alpha_l \}}$ of $G$. 
The set $S(\tau)$ is $\{ \alpha_1,\ldots,\alpha_l\}$, and the localization $\s_{S(\tau)}$ is
\[
\begin{picture}(14700,600)(-300,-300)
\put(0,0){\usebox{\mediumam}}
%\put(5400,0){\usebox{\plusbm}}
\end{picture}
\]
\[
\s_{S(\tau)}=(\{ \alpha_1+\ldots+\alpha_l \}, \{ \alpha_2,\ldots,\alpha_{l-1} \}, \emptyset)
\]
whose wonderful subgroup of $\SL(l+1)$ is $K_\tau=\GL(l)$. 
The morphism $\phi_{\Delta(\tau)}$ restricted to $\s_{S(\tau)}$ corresponds to the inclusion $\GL(l)\subset R$ 
where $R$ is the parabolic subgroup of $\SL(l+1)$ associated to $S_\tau\setminus\{\alpha_l\}$.
\end{example}

\subsubsection{Non-overlapping tails}\label{s:tailno}
Maintaining the notation of Definitions~\ref{def:tailroots} and \ref{def:tailpairs}, let $l$ be such that $\alpha_l$ is the unique simple root in $S(\tau)$ non-orthonogal to $\tau$. It exists: we claim that the indecomposability of $\s$ implies that $\supp\,\tau$ is not a connected component of $S$. In the case of an exceptional tail this is immediate from Definition~\ref{def:tailpairs}. In the case of a non-overlapping classical tail one can easily see that no other spherical root can have support on $\supp\,\tau$, so $\supp\,\tau=\overline{\supp\,\tau}$ would imply that $\s$ is the product of its two localizations $\s_{\tau}$ and $\s_{\Sigma\setminus\tau}$.

The homogeneous space $G/K^\tau$ is a parabolic induction by means of the parabolic subgroup $Q$ of $G$ containing $B_-$ 
and associated to the set of simple roots $S_Q=(S^p/\Delta(\tau))\cup\supp\,\tau$. 
Moreover, $Q$ is minimal among the parabolic subgroups of $G$ containing $K^\tau$. 
Notice that $\alpha_l$ is not in $S_Q$, and that $G_{\supp\,\tau}$ acts transitively on $G/K^\tau$: 
we call its stabilizer $K^\tau_{\supp\,\tau}\subset G_{\supp\,\tau}$. 
In other words, $G_{S_Q}=G_{S(\tau)\cap S_Q}\times G_{\supp\,\tau}$ 
and $K^\tau=(G_{S(\tau)\cap S_Q}\times K^\tau_{\supp\,\tau})Q^r$.

On the other hand let us consider the restriction $\Delta(\tau)|_{S(\tau)}$ of the set of colors $\Delta(\tau)$ to $S(\tau)$. From Definition~\ref{def:tails} and the assumption that $\Delta(\tau)$ is minimal, it follows that $\Delta(\tau)|_{S(\tau)}$ is a minimal parabolic set of colors of $\s_{S(\tau)}$. The corresponding parabolic subgroup $R$ of $G_{S(\tau)}$ can be taken to contain $B_-\cap G_{S(\tau)}$, and has then $S(\tau)\cap S_Q$ as set of associated simple roots. Now choose $K_\tau$ contained in $R$, which implies $K_\tau^u\subseteq R^u$, and choose a Levi part $L_{K_\tau}$ of $K_\tau$ contained in $L_{R}$.

By construction there exists a central isogeny $\pi\colon G_{S(\tau)\cap S_Q} Z(L_{S_Q})\to L_{R}$, and we can define $L = \pi^{-1}(L_{K_\tau}) \times K^\tau_{\supp\,\tau}$ as a subgroup of $L_{S_Q}$. Moreover, the $L_{K_\tau}$-submodule $K_\tau^u$ of $R^u$ corresponds to a $L$-stable subalgebra $\Lie\,U$ of $\Lie\,Q^u$, in such a way that the two quotients $\Lie\,R^u/\Lie\,K_\tau^u$, $\Lie\,Q^u/\Lie\,U$ are isomorphic as $\pi^{-1}(L_{K_\tau})$-modules.

The candidate $H$ is then $H=U\,L$. Table~\ref{table:tails} reports the groups $K^\tau_{\supp\,\tau}$ and $G_{\supp\,\tau}$ for all non-overlapping tails.

\begin{table}
\caption{Factors of Levi subgroups for non-overlapping tail systems}\label{table:tails}
\begin{center}
\begin{tabular}{ccc}
tail & $K^\tau_{\supp\,\tau}$ & $G_{\supp\,\tau}$ \\
\hline
$b(m)$ & $\SO(2m)$  & $\SO(2m+1)$ \\
$2b(m)$ & $\N(\SO(2m))$  & $\SO(2m+1)$ \\
$d(m)$ & $\SO(2m-1)$ & $\SO(2m)$ \\
$(aa,aa)$ & $\N(\diag(\SL(3)))$ & $\SL(3)\times\SL(3)$ \\
$(d3,d3)$ & $\N(\Sp(6))$ & $\SL(6)$ \\
$(d5,d5)$ & $\N(\mathsf F_4)$ & $\mathsf E_6$ \\
$(2a,2a)$ & $\N(\SO(3))$ & $\SL(3)$
\end{tabular}
\end{center}
\end{table}

\begin{example} Let us describe more explicitly the choice of $\Lie\,U$ inside $\Lie\,Q^u$ for a tail of type $b(m)$, the other non-overlapping cases being similar. Let us denote by $\tilde R$ the parabolic subgroup of $G_{S(\tau)}$ containing $B_-\cap G_{S(\tau)}$ and associated to the set of simple roots $S(\tau)\setminus\{\alpha_l\}$, and by $\tilde Q$ the parabolic subgroup of $G$ containing $B_-$ and associated to $S\setminus\{\alpha_l\}$. It is harmless to assume here that $G$ is simply-connected, and we denote by $\omega_\alpha\in\mathcal X(T)$ the fundamental dominant weight associated to $\alpha$, for any $\alpha\in S$.

The $L_{\tilde R}$-module $\Lie\,\tilde{R}^u$ is simple of highest weight $-\alpha_l$, that is, the fundamental dominant weight $\omega_{\alpha_{l-1}}$ restricted to the maximal torus of $G_{S(\tau)\setminus\{\alpha_l\}}$. On the other hand the quotient $L_{\tilde Q}$-module $\Lie\,\tilde{Q}^u/[\Lie\,\tilde{Q}^u,\Lie\,\tilde{Q}^u]$ is simple of highest weight $-\alpha_l$, that is, $\omega_{\alpha_{l-1}}+\omega_{\alpha_{l+1}}$ restricted to the maximal torus of $G_{S\setminus\{\alpha_l\}}$. Therefore, $\Lie\,\tilde{Q}^u/[\Lie\,\tilde {Q}^u,\Lie\,\tilde{Q}^u]\cong \Lie\,\tilde{R}^u\otimes V(\omega_{\alpha_{l+1}})$, where $V(\lambda)$ is the simple $\SO(2m+1)$-module of highest weight $\lambda$. The module $V(\omega_{\alpha_{l+1}})$ decomposes into the sum $W(\omega_{\alpha_{l+1}})\oplus W(0)$ as $\SO(2m)$-module, where $W(\mu)$ is the simple $\SO(2m)$-module of highest weight $\mu$. We obtain the desired isomorphism between the quotients $\Lie\,R^u/\Lie\,K_\tau^u$ and $\Lie\,Q^u/\Lie\,U$ since $[\Lie\,\tilde{Q}^u,\Lie\,\tilde{Q}^u]$ and $\Lie\,\tilde{R}^u\otimes W(\omega_{\alpha_{l+1}})$ are included in $\Lie\,U$, and $\Lie\,\tilde{R}^u$ is isomorphic to $\Lie\,\tilde{R}^u\otimes W(0)\subset \Lie\,\tilde{R}^u\otimes V(\omega_{\alpha_{l+1}})$.
\end{example}

\subsubsection{Overlapping tails}\label{s:tailC}
Here we use the same approach as before, but we have to deal with additional complications. 
These are due to the possibility that different spherical systems $\s$ have the same tail $\tau$ 
and the same localization $\s_{S(\tau)}$. 
Essentially, this is reflected in different non-conjugated choices of $R$ and $K_\tau$ that we use to define $H$.

\begin{example} Suppose that $\s_{S(\tau)}$ is the following system.
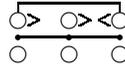
\[
\begin{picture}(4200,2000)(-300,-1000)
\multiput(0,0)(1800,0){2}{\usebox{\edge}}
\multiput(0,0)(1800,0){3}{\usebox{\aone}}
\multiput(0,900)(3600,0){2}{\line(0,1){300}}
\put(0,1200){\line(1,0){3600}}
\put(0,600){\usebox{\toe}}
\put(3600,600){\usebox{\tow}}
\put(1800,600){\usebox{\toe}}
\end{picture}
\]
Its associated wonderful subgroup of $\SL(4)$ can be chosen 
to have Levi subgroup equal to the standard Levi of the parabolic $P_{-\{\alpha_1 \}}$ of $\SL(4)$, 
or of the parabolic $P_{-\{\alpha_3 \}}$ (see \cite[Table 3, Case 9 with $p=1$]{BP05}). 
These two parabolic subgroups of $\SL(4)$ play the role of $R$ in two different tail cases, which are the following.
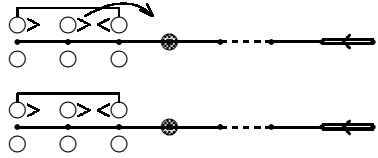
\[ 
\begin{picture}(12900,2400)(-300,-1200)
\multiput(0,0)(1800,0){2}{\usebox{\edge}}
\put(3600,0){\usebox{\shortcm}}
\multiput(0,0)(1800,0){3}{\usebox{\aone}}
\multiput(0,900)(3600,0){2}{\line(0,1){300}}
\put(0,1200){\line(1,0){3600}}
\put(0,600){\usebox{\toe}}
\put(3600,600){\usebox{\tow}}
\put(1800,600){\usebox{\toe}}
\put(1800,600){\usebox{\tobe}}
\end{picture}
\]
\[ 
\begin{picture}(12900,2700)(-300,-1350)
\multiput(0,0)(1800,0){2}{\usebox{\edge}}
\put(3600,0){\usebox{\shortcm}}
\multiput(0,0)(1800,0){3}{\usebox{\aone}}
\multiput(0,900)(3600,0){2}{\line(0,1){300}}
\put(0,1200){\line(1,0){3600}}
\put(0,600){\usebox{\toe}}
\put(3600,600){\usebox{\tow}}
\put(1800,600){\usebox{\toe}}
%\put(1800,-600){\qbezier(600,-300)(2100,-1200)(3000,-300)}
%\multiput(4800,-900)(-5,-10){45}{\line(0,-1){30}}
%\multiput(4800,-900)(-10,-5){45}{\line(-1,0){30}}
\end{picture}
\]
%The two rank $1$ quotients by $\Delta''$ are resp.:
%\[
%\s/\Delta' = 
%\begin{picture}(12900,1200)(-300,-300)
%\multiput(0,0)(1800,0){2}{\usebox{\edge}}
%\put(3600,0){\usebox{\shortcm}}
%\multiput(1800,0)(1800,0){2}{\usebox{\wcircle}}
%\end{picture}
%\]
%\[
%\s/\Delta' = 
%\begin{picture}(12900,1200)(-300,-300)
%\multiput(0,0)(1800,0){2}{\usebox{\edge}}
%\put(3600,0){\usebox{\shortcm}}
%\end{picture}
%\]
\end{example}

The same phenomenon can also involve $K_\tau$. 
For this reason, the next lemma will be used in the definition of $H$ for non-overlapping tails.

\begin{lemma}\label{lemma:Levitails}
In the hypotheses of Theorem~\ref{thm:tails}, there exists a maximal set of colors $\Delta'(\tau)\subseteq\Delta(\tau)$ of $\s$ of type ($\mathscr L$) such that $\Delta'(\tau)|_{S(\tau)}$ and $\Delta(\tau)|_{S(\tau)}$ give a decomposition of $\Delta|_{S(\tau)}$of type ($\mathscr L$ $\mathscr R$ $\mathscr P$).
\end{lemma}
\begin{proof}
Since we assumed that $\Delta(\tau)$ is minimal with respect to its properties, the localization $\Delta(\tau)|_{S(\tau)}$ is minimal parabolic. It contains then a set of colors $\Delta_0$ which has the required properties of $\Delta'(\tau)|_{S(\tau)}$.

Now the colors of $\s_{S(\tau)}$ are identified with a subset of the set of colors of $\s$, and we can take $\Delta'(\tau)$ to be $\Delta_0$ considered as a subset of the set of colors of $\s$.
\end{proof}

\begin{example} If $\s_{S(\tau)}$ is
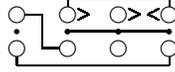
\[
\begin{picture}(6000,2400)(-300,-1200)
\multiput(1800,0)(1800,0){2}{\usebox{\edge}}
\multiput(0,0)(1800,0){4}{\usebox{\aone}}
\multiput(1800,900)(3600,0){2}{\line(0,1){300}}
\put(1800,1200){\line(1,0){3600}}
\multiput(0,-900)(5400,0){2}{\line(0,-1){300}}
\put(0,-1200){\line(1,0){5400}}
\put(300,600){\line(1,0){600}}
\put(1500,-600){\line(-1,0){600}}
\put(900,-600){\line(0,1){1200}}
\put(3600,600){\usebox{\toe}}
\put(1800,600){\usebox{\toe}}
\put(5400,600){\usebox{\tow}}
\end{picture}
\]
then two maximal sets of type $(\mathscr L)$ exist, with the following quotients.
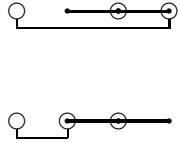
\[
\begin{picture}(6000,2400)(-300,-1200)
\multiput(1800,0)(1800,0){2}{\usebox{\edge}}
\put(0,0){\usebox{\wcircle}}
\multiput(3600,0)(1800,0){2}{\usebox{\wcircle}}
\multiput(0,-300)(5400,0){2}{\line(0,-1){300}}
\put(0,-600){\line(1,0){5400}}
\end{picture}
\]
\[
\begin{picture}(6000,2400)(-300,-1200)
\multiput(1800,0)(1800,0){2}{\usebox{\edge}}
\multiput(0,0)(1800,0){3}{\usebox{\wcircle}}
\multiput(0,-300)(1800,0){2}{\line(0,-1){300}}
\put(0,-600){\line(1,0){1800}}
\end{picture}
\]
These quotients are equal to $\s_{S(\tau)}/(\Delta'(\tau)|_{S(\tau)})$ resp.\ for the following two systems $\s$.
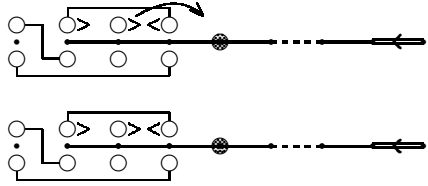
\[
\begin{picture}(14700,2400)(-300,-1200)
\multiput(1800,0)(1800,0){2}{\usebox{\edge}}
\put(5400,0){\usebox{\shortcm}}
\multiput(0,0)(1800,0){4}{\usebox{\aone}}
\multiput(1800,900)(3600,0){2}{\line(0,1){300}}
\put(1800,1200){\line(1,0){3600}}
\multiput(0,-900)(5400,0){2}{\line(0,-1){300}}
\put(0,-1200){\line(1,0){5400}}
\put(300,600){\line(1,0){600}}
\put(1500,-600){\line(-1,0){600}}
\put(900,-600){\line(0,1){1200}}
%\put(1800,1200){\line(1,0){3600}}
\put(3600,600){\usebox{\tobe}}
\put(3600,600){\usebox{\toe}}
\put(1800,600){\usebox{\toe}}
\put(5400,600){\usebox{\tow}}
\end{picture}
\]
\[ 
\begin{picture}(14700,2400)(-300,-1200)
\multiput(1800,0)(1800,0){2}{\usebox{\edge}}
\put(5400,0){\usebox{\shortcm}}
\multiput(0,0)(1800,0){4}{\usebox{\aone}}
\multiput(1800,900)(3600,0){2}{\line(0,1){300}}
\put(1800,1200){\line(1,0){3600}}
\multiput(0,-900)(5400,0){2}{\line(0,-1){300}}
\put(0,-1200){\line(1,0){5400}}
\put(300,600){\line(1,0){600}}
\put(1500,-600){\line(-1,0){600}}
\put(900,-600){\line(0,1){1200}}
%\put(1800,1200){\line(1,0){3600}}
%\put(3600,600){\usebox{\tobe}}
\put(3600,600){\usebox{\toe}}
\put(1800,600){\usebox{\toe}}
\put(5400,600){\usebox{\tow}}
\end{picture}
\]
\end{example}

Let us now define the candidate $H$ in general. 
Denote by $Q'$ the parabolic subgroup of $G$ containing $B_-$ 
and corresponding to the quotient of $\s$ by a minimal parabolic set containing $\Delta(\tau)$. 
Denote by $R'$ the parabolic subgroup of $G_{S(\tau)}$ containing $B_-\cap G_{S(\tau)}$ 
and corresponding to the quotient of $\s_{S(\tau)}$ by $\Delta(\tau)|_{S(\tau)}$. 
With the same notation as in Definition~\ref{def:tailroots}, let $l=n-m$. 
Notice that $\alpha_{l+1}$ belongs to the sets of simple roots associated to $Q'$ and $R'$ 
if and only if $\alpha_{l+1}\in S^p/\Delta(\tau)$. 

Let us also define parabolic subgroups $Q$ and $R$ obtained resp.\ from $Q'$ and $R'$ 
by adding $\alpha_{l+1}$ to their sets of associated simple roots. 
Of course, $\alpha_{l+1}$ may already be there, thus giving $Q'=Q$ and $R'=R$. 
In any case $Q=P_{-S_{Q}}$ where as before $S_Q=S^p/\Delta(\tau)\cup\supp\,\tau$.

We have $G_{S_{Q}}= G_{S(\tau)\cap S^p/\Delta(\tau)}\times\Sp(2m)$, 
and $G_{S(\tau)\cap S_{Q}}=G_{S(\tau)\cap S^p/\Delta(\tau)}\times\SL(2)$. 
We choose a subgroup $\SL(2)\times\Sp(2m-2)$ in $\Sp(2m)\subseteq G_{S_{Q}}$, 
and identify its $\SL(2)$ factor with the one of $G_{S(\tau)\cap S_{Q}}$.

Now we choose $K_\tau$ inside $R'$, in such a way that $K_\tau\subseteq R'$ factors through the co-connected inclusion associated to $\Delta'(\tau)|_{S(\tau)}$. Again we have $K_\tau^u\subseteq (R')^u$, and we choose a Levi component $L_{K_\tau}$ of $K_\tau$ contained in $L_{R'}$.

As before, we set $L=\pi^{-1}(L_{K_\tau})\times \Sp(2m-2)$, and define $U$ in the same way as in \S\ref{s:tailno}. Notice that the simple $\Sp(2m)$-module $V(\omega_{\alpha_{l+1}})$ of highest weight $\omega_{\alpha_{l+1}}$ decomposes into $W_1\oplus W_2$ as $\SL(2)\times\Sp(2m-2)$-module, where $W_1$ and $W_2$ are the defining modules resp.\ of the first and second factor. Finally, put $H=U\,L$.

\subsection{Geometric realizability of a system with tails}\label{s:tailsproof}
We use here all the elements introduced in \S6.2, in particular $\s$ is an indecomposable spherical $G$-system without quotients of higher defect and with tail in $\tau$ and let $H$ be the subgroup of $G$ defined in \S\ref{s:tailno} or \S\ref{s:tailC}.  Denote also by $M$ the factor $\Sp(2m-2)$ of $L$ acting trivially on $Q^u/U$ if the tail is of type $c(m)$, or the factor $K^\tau_{\supp\,\tau}$ otherwise. 

\begin{lemma} \label{lemma:tails1}
Suppose that $H$ is a generic stabilizer of a wonderful $G$-variety $X$ with spherical system $\s_X$. Then $\tau \subseteq \Sigma_X$. In the case of a tail $\tau=\{\gamma\}$ of type $b(1)$, both colors moved by $\gamma$ are free, and the Cartan pairing is equal on them.
\end{lemma}
\begin{proof}
Let us choose a 1-PSG $\lambda\colon\C^*\to T$ adapted to $\supp\,\tau$ (see \cite{Lu97}), 
and consider the definition of the candidate $H$. 
It is not $B$-spherical, in other words it is the stabilizer of a point $x\in X$ not lying on the open $B$-orbit. 
Nevertheless, $uHu^{-1}$ is $B$-spherical for any $u$ varying in a non-empty open subset of $B_-^u$. 
This means that $ux$ lies on the open $B$-orbit of $X$, and thus the limit
\[
x_0=\lim_{t\to0}\lambda(t)ux
\]
lies on the open $G_{\supp\,\tau}$-orbit of the localization $X_{\supp\,\tau}$ of $X$ in $\supp\,\tau$ (see \cite[\S1.1]{Lu97}). It is possible to write any $u\in B_-^u$ as a product $u=vw$ such that $w\in H$ and $v$ commutes with $M$. The consequence is that $M\subseteq uHu^{-1}$ for any $u\in B_-^u$.  The elements $\lambda(t)$ commute with $M$ for any $t$: we conclude that $M$ stabilizes $x_0$, and hence it is contained in a generic stabilizer of $X_{\supp\,\tau}$. Now it is easy to list all spherical subgroups of $G_{\supp\,\tau}$ containing $M$: they have rank $1$ or $2$ and in all cases except for $b(m)$ we have $\tau\subseteq\Sigma_{X_{\supp\,\tau}}\subseteq\Sigma_{X}$. If the tail $\tau=\{\gamma\}$ is of type $b(m)$, then the analysis of the spherical subgroups of $G_{\supp\,\tau}$ containing $M$ shows that either $\gamma$ or $2\gamma$ is a spherical root of $G/H$. To exclude $2\gamma$, it is enough to notice that the wonderful subgroup $K^\tau$ of $\s/\Delta(\tau)$ contains $H$ and the inclusion $K^\tau\supseteq H$ is co-connected: since $G/K^\tau$ has spherical root $\gamma$, and no spherical root can be doubled when passing from $\s_X$ to a quotient, we conclude that $\gamma$ is a spherical root of $\s_X$.

It remains to prove the last assertion in the case of a tail type $b(1)$. 
Consider the spherical system $\t$ equal to $\s$ except for the fact that $\gamma$ is replaced by $2\gamma$. 
Then $\t$ has a tail of type $2b(1)$ in $\tau'=\{2\gamma\}$, 
and that the corresponding candidate subgroup $N$ of $G$ contains $H$, 
is contained in $\N_GH$, and $|N/H| = 2$. 
Moreover, if $H$ is wonderful then $N$ is also wonderful, 
and we can apply the first part of the proof to both spherical systems $\s$ and $\t$. 
We conclude that $\gamma$ (resp.\ $2\gamma$) is a spherical root of the wonderful embedding of $G/H$ (resp.\ $G/N$). 
In this case the simple root $\gamma$ moves two colors of $G/H$, but only one of $G/N$, 
and this implies that the two colors of $G/H$ moved by $\gamma$ are exchanged 
by the $G$-equivariant automorphism of $G/H$ induced by the non-trivial element of $N/H$.

It follows that the Cartan pairing on these two colors is equal, and also that they are both free.
\end{proof}

\begin{corollary}\label{cor:tails}
Suppose that $H$ is a generic stabilizer of a wonderful $G$-variety $X$ with spherical system $\s_X$. Then:
\begin{enumerate}
%\item The quotient $\s_X/\Delta''_X$ of $\s_X$ induced by the inclusion $H\subseteq M$ is smooth.
\item any spherical root $\sigma\in\Sigma_X\setminus\tau$ satisfies $\supp\,\sigma\subseteq S(\tau)$,
\item the spherical system $\s_X$ has a tail in $\tau$.
\end{enumerate}
\end{corollary}
\begin{proof}
The first assertion is verified by checking all possible spherical roots (of any wonderful $G$-variety) that could be supported on simple roots that are non-orthogonal to $\supp\,\tau$. This is easily accomplished using the classification of wonderful varieties of rank $2$.

For the second assertion we must verify Definition~\ref{def:tails}: here it is enough to define $\Delta_X(\tau)$ to be the set of colors associated to the inclusion $H\subseteq K^\tau$.
\end{proof}

\begin{lemma}\label{lemma:tails2}
Under the assumptions of Corollary~\ref{cor:tails}, the localization of $\s_X$ in $S(\tau)$ is equal to $\s_{S(\tau)}$.
\end{lemma}
\begin{proof}
We suppose at first that the tail is not of type $b(1)$, and let us call $H_\tau$ a generic stabilizer of the localization $X_{S(\tau)}$. We claim that $H_\tau$ can be chosen to have $L_{K_\tau}$ as Levi factor.

Corollary~\ref{lemma:tails1} implies that restriction of sets of colors induces a bijection between the set of distinguished sets of colors of $X$ contained in $\Delta_X(\tau)$ and the set of distinguished sets of colors of $X_{S(\tau)}$ contained in $\Delta_X(\tau)|_{S(\tau)}$. Notice that in this case restricting a set of colors of $X$ is just considering is as a set of colors of $X_{S(\tau)}$.

This bijection preserves the property of being of type ($\mathscr L$). On the other hand $K^\tau_{\supp\,\tau}$ is a very reductive subgroup of $L_{\supp\,\tau}$: therefore there exists a maximal set of colors $\Delta'_X(\tau)$ of type ($\mathscr L$) inside $\Delta_X(\tau)$. It follows that $H_\tau$ has Levi factor which differs from $L_{K_\tau}$ only by its connected center. 

Let $\Delta''_X(\tau)\supseteq\Delta_X(\tau)$ be a minimal parabolic set of colors of $\s_X$. We deduce from Lemma~\ref{lemma:tails1} that $\Delta''_X(\tau)$ contains all colors moved by simple roots in $S\setminus S(\tau)$, and that $\Delta''_X(\tau)|_{S(\tau)}$ is minimal parabolic of $\s_X|_{S(\tau)}$. Moreover, on $\tau$ all colors in $\Delta_X(\tau)$ are $0$ and all colors in $\Delta_X\setminus\Delta''_X(\tau)$ are $\leq 0$.

From these facts it is possible to deduce that the formula of Corollary~\ref{cor:c1} gives the same result whether applied to $X$ and $\Delta''_X(\tau)$ or to $X_{S(\tau)}$ and $\Delta''_X(\tau)|_{S(\tau)}$: for each type of tail it is an elementary computation involving the values on $\tau$ of the color(s) moved by the simple root(s) in $S\setminus S(\tau)$. Comparing the connected center of $L$ and of $L_{K_\tau}$ we conclude that $L_{K_\tau}$ and a Levi subgroup of $H_\tau$ also have same connected center, and our claim follows.

Now we prove that $H_\tau$ and $K$ are conjugated. 
It is possible to choose a 1-PSG $\lambda\colon \C^*\to T$ 
such that $\lambda(t)$ commutes with $L_{S(\tau)}$ for all $t\in\C^*$, 
and $\langle \lambda, \alpha\rangle =1$ for all $\alpha\in S\setminus S(\tau)$. 
Corollary~\ref{cor:tails} implies that $\langle\lambda,\sigma\rangle > 0$ for all $\sigma\in\tau$, 
and $\langle\lambda,\sigma\rangle = 0$ for all $\sigma\in\Sigma_X\setminus\tau$. 
The consequence is that the limit
\[
x_0=\lim_{t\to0} \lambda(t)x
\]
lies on the open $G$-orbit of $X_{\Sigma_X\setminus\tau}$, and $\pi^{-1}(L_{K_\tau})$ is contained in the stabilizer $\lambda(t) H \lambda(t)^{-1}$ of $\lambda(t)x$ for all $t\in\C^*$. Hence it is also contained in the stabilizer $H^{\Sigma\setminus\tau}$ of $x_0$.

Thanks to the definition of $H$, for a general choice of $K^\tau_{\supp\,\tau}$ inside $L_{\supp\,\tau}$ we have
\[
\lim_{t \to 0} \Ad(\lambda(t))\left(\Lie\,L_{K_\tau} \oplus \Lie\,U\right) = \Lie\,L_{K_\tau} \oplus \Lie\,K_\tau^u.
\]
where we identify $\Lie\,G_{S(\tau)\cap S_Q} \oplus\Lie\,Z(L_{S_Q})$ with $\Lie\, L_{R}$ via the isogeny $\pi$.

Hence $\Lie\,L_{K_\tau} \oplus \Lie\,K_\tau^u$ is also contained in $\Lie\,H^{\Sigma\setminus\tau}$. Finally, Lemma~\ref{lemma:tails1} and Corollary~\ref{cor:tails} assure that $X_{\Sigma\setminus\tau}$ is the parabolic induction of $X_{S(\tau)}$ by means of a parabolic $P$ containing $P_{-S(\tau)}$. In other words $\Lie\,H^{\Sigma\setminus\tau}$ is the sum of $\Lie\,H_\tau$ and the radical of $\Lie\,P$, up to conjugating $H_\tau$ if necessary. We deduce that $H_\tau\supseteq L_{K_\tau}\,K_\tau^u=K_\tau$, and for dimension reasons they are equal.

If the tail $\tau=\{\gamma\}$ has type $b(1)$, 
we consider the spherical system $\t$ and the corresponding wonderful subgroup $N$ of $G$ 
as defined in the proof of Lemma~\ref{lemma:tails1}. 
Let $Y$ be the wonderful completion of $G/N$, and $\s_Y$ its spherical system. 
The first part of the proof assures that the localization of $\s_Y$ on $S(\tau)$ is equal to $\s_{S(\tau)}$. 
On the other hand, the fact that $|N/H|=2$ together with $\gamma\in\Sigma_X$ and $2\gamma\in\Sigma_Y$ 
implies that the localizations of $\s_X$ and of $\s_Y$ in $S(\tau)$ are equal: the proof is complete.
\end{proof}

\begin{lemma}\label{lemma:tails3}
Suppose that $H$ is a generic stabilizer of a wonderful $G$-variety $X$ with spherical system $\s_X$. 
If $\tau$ is of type $c(m)$, suppose in addition that the candidate subgroup of $G$ 
defined starting from any other spherical $G$-system $\s'$ with the same tail $\tau$ 
and the same localization $\s'_{S(\tau)}=\s_{S(\tau)}$ is a wonderful subgroup of $G$. 
Then the spherical system $\s_X$ is equal to $\s$.
\end{lemma}
\begin{proof}
Applying Lemma~\ref{lemma:tails1} and Lemma~\ref{lemma:tails2} we know that $\tau\subset\Sigma_X$ and that the localization of $\s_X$ in $S(\tau)$ is $\s_{S(\tau)}$. From the co-connected inclusion $H\subseteq K^\tau$, we also deduce that $\tau$ is a tail for $\s_X$. This finishes the proof if the tail is not of type $c(m)$.

If $\tau=\{\gamma\}$ has type $c(m)$, the systems $\s_X$ and $\s$ may both share these properties, and yet be different. More precisely, if we fix such $\tau$ and the localization $\s_{S(\tau)}$, there can be up to two different spherical systems $\s=\s_1$ and $\s_2$, both with localization in $S(\tau)$ equal to $\s_{S(\tau)}$ and set of spherical roots $\Sigma_X$, and differing only for the choice of the color $D$ of $\s_{S(\tau)}$ such that $c(D,\gamma) < 0$. If two such systems exist, then both choices for the color $D$ are moved by $\alpha_l$ (in the notation of \S\ref{s:tailC}) and are both free.

In this situation, either $\s_2$ does not admit $\tau$ as a tail, or it does. 
The first case implies $\s_X = \s$. We must check therefore that, in the second case, 
for the two systems $\s_1$, $\s_2$ our procedure provides two candidates $H_1$ and $H_2$ respectively 
(both wonderful thanks to our assumptions), such that each corresponds to the right spherical system. 
It is convenient here to denote by $D^+_i$ and $D^-_i$ the two colors moved by a spherical root $\alpha_i\in S\cap\Sigma$.

\noindent
{\em Step 1.} Let us suppose that no spherical root $\sigma\neq\gamma$ is supported on $\alpha_{l+1}$. 
This implies that $\s_{S(\tau)}$ must be non-cuspidal, with $\alpha_{l+1}$ outside the support of its spherical roots. 
The two systems $\s_1$ and $\s_2$ in the vicinity of $\gamma$ are as follows.
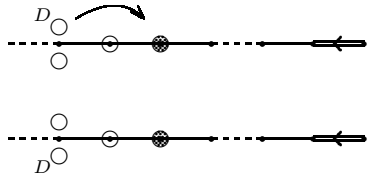
\[
\begin{picture}(12900,2400)(-300,-1200)
\put(0,0){\thicklines\multiput(0,0)(400,0){5}{\line(1,0){200}}}
\put(1800,0){\usebox{\edge}}
\put(3600,0){\usebox{\wcircle}}
\put(1800,0){\usebox{\aone}}
\put(3600,0){\usebox{\shortcm}}
\put(1800,600){\usebox{\tobe}}
\put(900,800){$\scriptstyle{D}$}
\end{picture}
\]
\[
\begin{picture}(12900,2400)(-300,-1200)
\put(0,0){\thicklines\multiput(0,0)(400,0){5}{\line(1,0){200}}}
\put(1800,0){\usebox{\edge}}
\put(3600,0){\usebox{\wcircle}}
\put(1800,0){\usebox{\aone}}
\put(3600,0){\usebox{\shortcm}}
\put(900,-1200){$\scriptstyle{D}$}
\end{picture}
\]
If $D_{l+2}$ is the unique color moved by $\alpha_{l+2}$, then $\Delta(\tau)$ and $\{ D_{l+2} \}$ decompose the system in both cases, which is absurd.

As a consequence, we may assume that some spherical root $\sigma\neq\gamma$ is supported on $\alpha_{l+1}$.

\noindent
{\em Step 2.} Let's consider first the case where $\sigma\notin S$. 
The only possibility is $\gamma = \alpha_{l+1} + \alpha_{l-1}$, 
and the two systems $\s_1$, $\s_2$ in the vicinity of $\gamma$ look like the following.
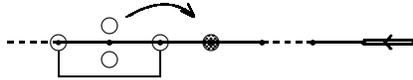
\[
\begin{picture}(14700,2400)(-300,-1200)
\put(0,0){\thicklines\multiput(0,0)(400,0){5}{\line(1,0){200}}}
\put(1800,0){\usebox{\wcircle}}
\multiput(1800,0)(1800,0){2}{\usebox{\edge}}
\put(3600,0){\usebox{\aone}}
\put(5400,0){\usebox{\shortcm}}
\put(5400,0){\usebox{\wcircle}}
\multiput(1800,-300)(3600,0){2}{\line(0,-1){900}}
\put(1800,-1200){\line(1,0){3600}}
\put(3600,600){\usebox{\tobe}}
\end{picture}
\] 
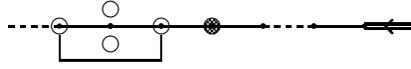
\[
\begin{picture}(14700,2400)(-300,-1200)
\put(0,0){\thicklines\multiput(0,0)(400,0){5}{\line(1,0){200}}}
\put(1800,0){\usebox{\wcircle}}
\multiput(1800,0)(1800,0){2}{\usebox{\edge}}
\put(3600,0){\usebox{\aone}}
\put(5400,0){\usebox{\shortcm}}
\put(5400,0){\usebox{\wcircle}}
\multiput(1800,-300)(3600,0){2}{\line(0,-1){900}}
\put(1800,-1200){\line(1,0){3600}}
%\put(3600,600){\usebox{\tobe}}
\end{picture}
\] 
The required set of colors $\Delta(\tau)$ exists only in the second diagram, so there is no ambiguity.

\noindent
{\em Step 3.} It remains the case where $\alpha_{l+1}\in\Sigma$. 
If $D^+_{l+1}$ and $D^-_{l+1}$ are free, then the two systems $\s_1$, $\s_2$ in the vicinity of $\gamma$ are as follows.
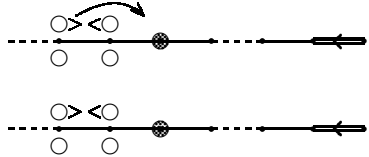
\[
\begin{picture}(12900,2400)(-300,-1200)
\put(0,0){\thicklines\multiput(0,0)(400,0){5}{\line(1,0){200}}}
\put(1800,0){\usebox{\edge}}
\multiput(1800,0)(1800,0){2}{\usebox{\aone}}
\put(3600,0){\usebox{\shortcm}}
\put(3600,600){\usebox{\tow}}
\put(1800,600){\usebox{\toe}}
\put(1800,600){\usebox{\tobe}}
\end{picture}
\]
\[
\begin{picture}(12900,2400)(-300,-1200)
\put(0,0){\thicklines\multiput(0,0)(400,0){5}{\line(1,0){200}}}
\put(1800,0){\usebox{\edge}}
\multiput(1800,0)(1800,0){2}{\usebox{\aone}}
\put(3600,0){\usebox{\shortcm}}
\put(3600,600){\usebox{\tow}}
\put(1800,600){\usebox{\toe}}
%\put(1800,600){\usebox{\tobe}}
\end{picture}
\]
The two quotients $\s/\Delta'(\tau)$ are resp.
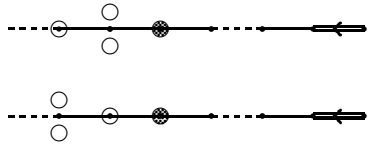
\[
\begin{picture}(12900,2400)(-300,-1200)
\put(0,0){\thicklines\multiput(0,0)(400,0){5}{\line(1,0){200}}}
\put(1800,0){\usebox{\edge}}
\put(1800,0){\usebox{\wcircle}}
\put(3600,0){\usebox{\aone}}
\put(3600,0){\usebox{\shortcm}}
\end{picture}
\] 
\[
\begin{picture}(12900,2400)(-300,-1200)
\put(0,0){\thicklines\multiput(0,0)(400,0){5}{\line(1,0){200}}}
\put(1800,0){\usebox{\edge}}
\put(3600,0){\usebox{\wcircle}}
\put(1800,0){\usebox{\aone}}
\put(3600,0){\usebox{\shortcm}}
\end{picture}
\] 
and they give two different candidates $H_1$, $H_2$. 
The existence of the co-connected inclusions associated to the two quotients above assures that the two subgroups of $G$ 
correspond to the right spherical systems.

We may now assume that $D^+_{l+1}$ and $D^-_{l+1}$ are not both free.

\noindent
{\em Step 4.} Let us suppose that both are non-free: one of them must take the value $-1$ on $\alpha_l$, say $D^+_{l+1}$. 
Consider the colors $D^+_{l}$ and $D^-_{l}$ moved by $\alpha_l$. 
If one of them is zero both on $\alpha_{l-1}$ and on $\alpha_{l+1}$, then we have the following.
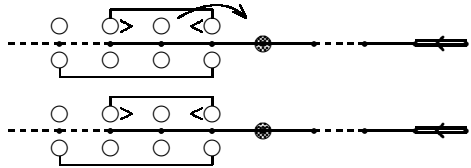
\[
\begin{picture}(14700,2400)(-300,-1200)
\put(0,0){\thicklines\multiput(0,0)(400,0){5}{\line(1,0){200}}}
\put(1800,0){\usebox{\shortsusp}}
\multiput(3600,0)(1800,0){2}{\usebox{\edge}}
\put(7200,0){\usebox{\shortcm}}
\multiput(1800,0)(1800,0){4}{\usebox{\aone}}
\multiput(3600,900)(3600,0){2}{\line(0,1){300}}
\put(3600,1200){\line(1,0){3600}}
\multiput(1800,-900)(5400,0){2}{\line(0,-1){300}}
\put(1800,-1200){\line(1,0){5400}}
\put(5400,600){\usebox{\tobe}}
\put(3600,600){\usebox{\toe}}
\put(7200,600){\usebox{\tow}}
\end{picture}
\]
\[
\begin{picture}(14700,2400)(-300,-1200)
\put(0,0){\thicklines\multiput(0,0)(400,0){5}{\line(1,0){200}}}
\put(1800,0){\usebox{\shortsusp}}
\multiput(3600,0)(1800,0){2}{\usebox{\edge}}
\put(7200,0){\usebox{\shortcm}}
\multiput(1800,0)(1800,0){4}{\usebox{\aone}}
\multiput(3600,900)(3600,0){2}{\line(0,1){300}}
\put(3600,1200){\line(1,0){3600}}
\multiput(1800,-900)(5400,0){2}{\line(0,-1){300}}
\put(1800,-1200){\line(1,0){5400}}
\put(3600,600){\usebox{\toe}}
\put(7200,600){\usebox{\tow}}
\end{picture}
\]
The first system does not admit a set of colors $\Delta(\tau)$, so there is no ambiguity.

We assume now that no color moved by $\alpha_{l}$ is zero on both $\alpha_{l-1}$ and $\alpha_{l+1}$. 
If $D^-_{l-1}$ is free, then we have the following.
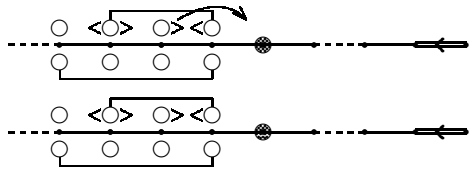
\[
\begin{picture}(14700,2400)(-300,-1200)
\put(0,0){\thicklines\multiput(0,0)(400,0){5}{\line(1,0){200}}}
\multiput(1800,0)(1800,0){3}{\usebox{\edge}}
\put(7200,0){\usebox{\shortcm}}
\multiput(1800,0)(1800,0){4}{\usebox{\aone}}
\multiput(3600,900)(3600,0){2}{\line(0,1){300}}
\put(3600,1200){\line(1,0){3600}}
\multiput(1800,-900)(5400,0){2}{\line(0,-1){300}}
\put(1800,-1200){\line(1,0){5400}}
\put(5400,600){\usebox{\tobe}}
\put(5400,600){\usebox{\toe}}
\put(3600,600){\usebox{\toe}}
\put(3600,600){\usebox{\tow}}
\put(7200,600){\usebox{\tow}}
\end{picture}
\]
\[
\begin{picture}(14700,2400)(-300,-1200)
\put(0,0){\thicklines\multiput(0,0)(400,0){5}{\line(1,0){200}}}
\multiput(1800,0)(1800,0){3}{\usebox{\edge}}
\put(7200,0){\usebox{\shortcm}}
\multiput(1800,0)(1800,0){4}{\usebox{\aone}}
\multiput(3600,900)(3600,0){2}{\line(0,1){300}}
\put(3600,1200){\line(1,0){3600}}
\multiput(1800,-900)(5400,0){2}{\line(0,-1){300}}
\put(1800,-1200){\line(1,0){5400}}
%\put(5400,600){\usebox{\tobe}}
\put(5400,600){\usebox{\toe}}
\put(3600,600){\usebox{\toe}}
\put(3600,600){\usebox{\tow}}
\put(7200,600){\usebox{\tow}}
\end{picture}
\]
Here the two quotients $\s/\Delta'(\tau)$ must be resp.
\[
\begin{picture}(16500,2400)(-300,-1200)
\put(0,0){\thicklines\multiput(0,0)(400,0){5}{\line(1,0){200}}}
\multiput(1800,0)(1800,0){3}{\usebox{\edge}}
\multiput(1800,0)(1800,0){4}{\usebox{\wcircle}}
\put(7200,0){\usebox{\shortcm}}
%\multiput(3600,900)(3600,0){2}{\line(0,1){300}}
%\put(3600,1200){\line(1,0){3600}}
\multiput(1800,-300)(5400,0){2}{\line(0,-1){300}}
\put(1800,-600){\line(1,0){5400}}
\end{picture}
\]
and
\[
\begin{picture}(16500,2400)(-300,-1200)
\put(0,0){\thicklines\multiput(0,0)(400,0){5}{\line(1,0){200}}}
\put(1800,0){\usebox{\atwo}}
\multiput(3600,0)(1800,0){2}{\usebox{\edge}}
\put(5400,0){\usebox{\wcircle}}
\put(7200,0){\usebox{\shortcm}}
\end{picture}
\]
and they give two different candidates $H_1$, $H_2$. 
Again, the existence of the co-connected inclusions associated to the two quotients above 
assures that the two subgroups of $G$ correspond to the right spherical systems.

If $D^-_{l-1}$ is not free, then we have
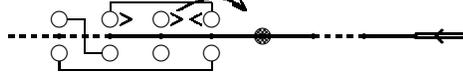
\[
\begin{picture}(16500,2400)(-300,-1200)
\put(0,0){\thicklines\multiput(0,0)(400,0){5}{\line(1,0){200}}}
\put(1800,0){\usebox{\shortsusp}}
\multiput(3600,0)(1800,0){2}{\usebox{\edge}}
\put(7200,0){\usebox{\shortcm}}
\multiput(1800,0)(1800,0){4}{\usebox{\aone}}
\multiput(3600,900)(3600,0){2}{\line(0,1){300}}
\put(3600,1200){\line(1,0){3600}}
\multiput(1800,-900)(5400,0){2}{\line(0,-1){300}}
\put(1800,-1200){\line(1,0){5400}}
\put(2100,600){\line(1,0){600}}
\put(3300,-600){\line(-1,0){600}}
\put(2700,-600){\line(0,1){1200}}
%\put(1800,1200){\line(1,0){3600}}
\put(5400,600){\usebox{\tobe}}
\put(5400,600){\usebox{\toe}}
\put(3600,600){\usebox{\toe}}
\put(7200,600){\usebox{\tow}}
\end{picture}
\]
The two quotients $\s/\Delta'(\tau)$ are resp.
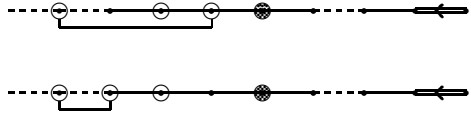
\[
\begin{picture}(16500,2400)(-300,-1200)
\put(0,0){\thicklines\multiput(0,0)(400,0){5}{\line(1,0){200}}}
\put(1800,0){\usebox{\shortsusp}}
\multiput(3600,0)(1800,0){2}{\usebox{\edge}}
\put(7200,0){\usebox{\shortcm}}
\put(1800,0){\usebox{\wcircle}}
\put(5400,0){\usebox{\wcircle}}
\put(7200,0){\usebox{\wcircle}}
\multiput(1800,-300)(5400,0){2}{\line(0,-1){300}}
\put(1800,-600){\line(1,0){5400}}
\end{picture}
\]
\[
\begin{picture}(16500,2400)(-300,-1200)
\put(0,0){\thicklines\multiput(0,0)(400,0){5}{\line(1,0){200}}}
\multiput(3600,0)(1800,0){2}{\usebox{\edge}}
\put(1800,0){\usebox{\shortsusp}}
\multiput(1800,0)(1800,0){3}{\usebox{\wcircle}}
\put(7200,0){\usebox{\shortcm}}
\multiput(1800,-300)(1800,0){2}{\line(0,-1){300}}
\put(1800,-600){\line(1,0){1800}}
\end{picture}
\]
and they give two different candidates $H_1$, $H_2$. 
Again, the existence of the co-connected inclusions associated to the two quotients above 
assures that the two subgroups of $G$ correspond to the right spherical systems.

\noindent
{\em Step 5.} Let us suppose that only one of $D^+_{l+1}$, $D^-_{l+1}$ is not free, say $D^+_{l+1}$. 
We leave for a moment undetermined the values of the colors $D^+_{l}$ and $D^-_{l}$, 
and focus on the value of $D^+_{l+1}$ on $\alpha_{l}$. We have two systems:
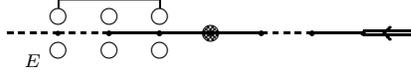
\[
\begin{picture}(12900,2400)(-300,-1200)
\multiput(0,0)(1800,0){2}{\thicklines\multiput(0,0)(400,0){5}{\line(1,0){200}}}
\put(3600,0){\usebox{\edge}}
\put(5400,0){\usebox{\shortcm}}
\multiput(1800,0)(1800,0){3}{\usebox{\aone}}
\multiput(1800,900)(3600,0){2}{\line(0,1){300}}
\put(1800,1200){\line(1,0){3600}}
%\put(3600,600){\usebox{\toe}}
\put(600,-1200){$\scriptstyle{E}$}
\end{picture}
\]
where $c(D^+_{l+1},\alpha_l)=0$, and
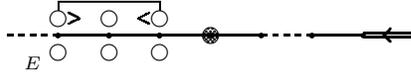
\[
\begin{picture}(12900,2400)(-300,-1200)
\put(0,0){\thicklines\multiput(0,0)(400,0){5}{\line(1,0){200}}}
\multiput(1800,0)(1800,0){2}{\usebox{\edge}}
\put(5400,0){\usebox{\shortcm}}
\multiput(1800,0)(1800,0){3}{\usebox{\aone}}
\multiput(1800,900)(3600,0){2}{\line(0,1){300}}
\put(1800,1200){\line(1,0){3600}}
%\put(3600,600){\usebox{\toe}}
\put(1800,600){\usebox{\toe}}
\put(5400,600){\usebox{\tow}}
\put(600,-1200){$\scriptstyle{E}$}
\end{picture}
\]
where $c(D^+_{l+1},\alpha_l)=-1$. Since $D^-_{l+1}$ is free, the color $E$ must also be free in both cases.

Hence in the first case $D^+_{l+1}$ is a non-free positive color, and the system admits a quotient of higher defect. On the other hand, it is easy to see that if $\s$ has tail of type $c(m)$ and no quotient of higher defect, then $\s_{S(\tau)}$ has no such quotient either. Then, from the definition of $H$ it is clear that $\s_X$ has no quotient of higher defect.

Thus, we may assume that $c(D^+_{l+1},\alpha_l)=-1$.

\noindent
{\em Step 6.} Now let us consider the colors moved by $\alpha_l$. 
If one of them is zero both on $\alpha_{l-1}$ and on $\alpha_{l+1}$, then we have the following.
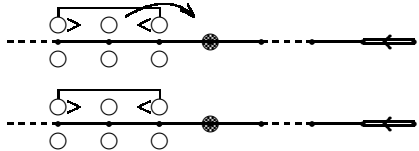
\[
\begin{picture}(12900,2400)(-300,-1200)
\put(0,0){\thicklines\multiput(0,0)(400,0){5}{\line(1,0){200}}}
\multiput(1800,0)(1800,0){2}{\usebox{\edge}}
\put(5400,0){\usebox{\shortcm}}
\multiput(1800,0)(1800,0){3}{\usebox{\aone}}
\multiput(1800,900)(3600,0){2}{\line(0,1){300}}
\put(1800,1200){\line(1,0){3600}}
%\put(3600,600){\usebox{\toe}}
\put(3600,600){\usebox{\tobe}}
\put(1800,600){\usebox{\toe}}
\put(5400,600){\usebox{\tow}}
\end{picture}
\]
\[
\begin{picture}(12900,2400)(-300,-1200)
\put(0,0){\thicklines\multiput(0,0)(400,0){5}{\line(1,0){200}}}
\multiput(1800,0)(1800,0){2}{\usebox{\edge}}
\put(5400,0){\usebox{\shortcm}}
\multiput(1800,0)(1800,0){3}{\usebox{\aone}}
\multiput(1800,900)(3600,0){2}{\line(0,1){300}}
\put(1800,1200){\line(1,0){3600}}
%\put(3600,600){\usebox{\toe}}
\put(1800,600){\usebox{\toe}}
\put(5400,600){\usebox{\tow}}
\end{picture}
\]
The first of the two systems does not admit a set of colors $\Delta(\tau)$, so there is no ambiguity.

\noindent
{\em Step 7.} We may now assume that no color moved by $\alpha_l$ is zero both on $\alpha_{l-1}$ and $\alpha_{l+1}$. 
Our systems are as follows.
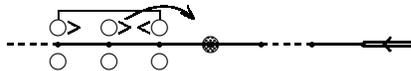
\[
\begin{picture}(12900,2400)(-300,-1200)
\put(0,0){\thicklines\multiput(0,0)(400,0){5}{\line(1,0){200}}}
\multiput(1800,0)(1800,0){2}{\usebox{\edge}}
\put(5400,0){\usebox{\shortcm}}
\multiput(1800,0)(1800,0){3}{\usebox{\aone}}
\multiput(1800,900)(3600,0){2}{\line(0,1){300}}
\put(1800,1200){\line(1,0){3600}}
\put(3600,600){\usebox{\toe}}
\put(3600,600){\usebox{\tobe}}
\put(1800,600){\usebox{\toe}}
\put(5400,600){\usebox{\tow}}
\end{picture}
\]
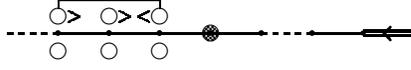
\[
\begin{picture}(12900,2400)(-300,-1200)
\put(0,0){\thicklines\multiput(0,0)(400,0){5}{\line(1,0){200}}}
\multiput(1800,0)(1800,0){2}{\usebox{\edge}}
\put(5400,0){\usebox{\shortcm}}
\multiput(1800,0)(1800,0){3}{\usebox{\aone}}
\multiput(1800,900)(3600,0){2}{\line(0,1){300}}
\put(1800,1200){\line(1,0){3600}}
\put(3600,600){\usebox{\toe}}
\put(1800,600){\usebox{\toe}}
\put(5400,600){\usebox{\tow}}
\end{picture}
\]

Both have a tail in $\tau$. Since $D^-_{l+1}$, $D^+_l$, $D^-_l$ are free, we conclude that $D^-_{l-1}$ is free, and that $D^+_{l-1}$, $D^-_{l-1}$, $D^+_{l-1}$ are all zero on spherical roots different from $\alpha_{l-1}$, $\alpha_l$, $\alpha_{l+1}$, $\gamma$.

It is easy to deduce that the two quotients by $\Delta(\tau)$ are resp.\ as follows.
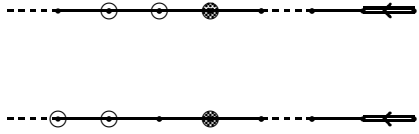
\[
\begin{picture}(12900,2400)(-300,-1200)
\put(0,0){\thicklines\multiput(0,0)(400,0){5}{\line(1,0){200}}}
\multiput(1800,0)(1800,0){2}{\usebox{\edge}}
\put(5400,0){\usebox{\shortcm}}
%\multiput(1800,0)(1800,0){3}{\usebox{\aone}}
%\multiput(1800,900)(3600,0){2}{\line(0,1){300}}
%\put(1800,1200){\line(1,0){3600}}
%\put(3600,600){\usebox{\toe}}
%\put(3600,600){\usebox{\tobe}}
%\put(1800,600){\usebox{\toe}}
%\put(5400,600){\usebox{\tow}}
\put(3600,0){\usebox{\wcircle}}
\put(5400,0){\usebox{\wcircle}}
\end{picture}
\]
\[
\begin{picture}(12900,2400)(-300,-1200)
\put(0,0){\thicklines\multiput(0,0)(400,0){5}{\line(1,0){200}}}
\multiput(1800,0)(1800,0){2}{\usebox{\edge}}
\put(5400,0){\usebox{\shortcm}}
%\multiput(1800,0)(1800,0){3}{\usebox{\aone}}
%\multiput(1800,900)(3600,0){2}{\line(0,1){300}}
%\put(1800,1200){\line(1,0){3600}}
%\put(3600,600){\usebox{\toe}}
%\put(1800,600){\usebox{\toe}}
%\put(5400,600){\usebox{\tow}}
\put(3600,0){\usebox{\wcircle}}
\put(1800,0){\usebox{\wcircle}}
\end{picture}
\]
It is evident that these correspond to two different $R$'s, 
which leads to two different wonderful subgroups $H_1$, $H_2$ of $G$. 
Finally, the existence of the co-connected inclusions associated to the two quotients above 
assures the good choice of our wonderful subgroups of $G$.
\end{proof}

\begin{proof}[Proof of Theorem~\ref{thm:tails}]
We begin by supposing that the tail is not of type $b(m)$, 
and that $\s_{S(\tau)}$ corresponds to a wonderful variety without non-trivial $G$-equivariant automorphisms. 
Then $\N_{L_{S(\tau)}}K_\tau=K_\tau$, and from the definition of $H$ it follows that $\N_GH=H$. 
The subgroup $H$ of $G$ is hence wonderful, call $X$ its wonderful variety. 
Applying Lemma~\ref{lemma:tails1}, Lemma~\ref{lemma:tails2} and Lemma~\ref{lemma:tails3} we know that $\s_X=\s$. 
Notice that if $\tau$ is overlapping, 
we have to check the additional hypothesis required by Lemma~\ref{lemma:tails3} 
regarding (if it exists) a spherical system $\s'\neq\s$ with same tail $\tau$ 
and same localization $\s'_{S(\tau)}=\s_{S(\tau)}$. 
This is easily done: the same argument we used above for $\s$ 
shows that the candidate subgroup of $G$ for $\s'$ is also equal to its normalizer, thus it is wonderful.

Now we drop our assumptions both on $\s_{S(\tau)}$ and on the tail. 
Let $\N_G(\s)$ be defined as in Lemma~\ref{lemma:normalizer}: 
it is a tail spherical system with localization $\N_G(\s_{S(\tau)})$ and same tail as $\s$, 
except for the case where the latter has tail $b(m)$ and $\N_G(\s)$ has tail $2b(m)$. 
By the first part of the proof $\N_G(\s)$ is geometrically realizable (with wonderful subgroup $N$ of $G$), 
and therefore $\s$ too, thanks to Lemma~\ref{lemma:normalizer}. 
Let $\widetilde H$ be the wonderful subgroup of $G$ associated with $\s$: 
it remains to show that $\widetilde H=H$.

As shown in the proof of Lemma~\ref{lemma:normalizer}, the quotient $N/\widetilde H$ is a diagonalizable group, 
hence $\widetilde H$ is defined inside $N$ as the intersection of the kernels of certain characters. 
The same considerations apply to $K_\tau$ and $N_\tau$ associated with resp.\ $\s_{S(\tau)}$ and $\N_G(\s_{S(\tau)})$. 
Moreover, $H$ has (by its definition) finite index in $N$, so it is the intersection of the kernels of some characters too.

At this point, the equality $\widetilde H=H$ follows if we check that $\widetilde H$ and $H$ have Levi subgroups equal up to conjugation. Recall that by definition a Levi subgroup of $H$ is $\pi^{-1}(L_{K_\tau})\times K^\tau_{\supp\,\tau}$, with the notation of \S\ref{s:generictails}. Let $(\Delta'(\tau),\Delta(\tau))$ be as in Lemma~\ref{lemma:Levitails} applied to the spherical system $\s$, and let $\Delta''(\tau)\supseteq \Delta(\tau)$ be a minimal parabolic set of colors of $\s$. Using Proposition~\ref{prop:parared} and Corollary~\ref{cor:c1}, we conclude that a Levi subgroup of $\widetilde H$ has indeed the form $\pi^{-1}(L_{K_\tau})\times K^\tau_{\supp\,\tau}$: the theorem follows.
\end{proof}

\subsection{The wonderful subgroup associated with a spherical system with arcs of type $a(m)$}
Let us now mention a variant of the above situation where the same arguments apply.

Let $\s$ be a spherical $G$-system such that:
\begin{itemize}
\item[-] there exists $\gamma\in\Sigma$ of type $a(m)$, i.e. of the form $\gamma=\alpha_{l+1}+\ldots+\alpha_{l+m}$ where $m\geq 1$ and $\{\alpha_{l+1}+\ldots+\alpha_{l+m}\}$ is a subset of $S$ of type $\mathsf A_m$,
\item[-] there exist $\gamma',\gamma''\in\Sigma$ non-orthogonal to $\supp\,\gamma$ with $\supp\,\gamma'\cap\supp\,\gamma''=\emptyset$ and such that, if $\{D\in\Delta : c(D,\gamma)=1\}=\{D',D''\}$, then $c(D',\gamma')=c(D'',\gamma'')=-1$ and $c(D',\gamma'')=c(D'',\gamma')=0$,
\item[-] there exists a good distinguished set of colors $\Delta(\gamma)$ such that $\Sigma/\Delta(\gamma)=\{\gamma\}$.
\end{itemize}

In this case we say that $\s$ has an {\em arc} in $\{\gamma\}$ of type $a(m)$.

Example:
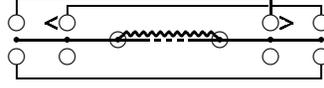
\[\begin{picture}(11400,2850)(-300,-1350)
\multiput(0,0)(1800,0){2}{\usebox{\edge}}
\multiput(0,0)(1800,0){2}{\usebox{\aone}}
\put(3600,0){\usebox{\shortam}}
\put(7200,0){\usebox{\edge}}
\put(9000,0){\usebox{\edge}}
\multiput(9000,0)(1800,0){2}{\usebox{\aone}}
\multiput(0,1500)(9000,0){2}{\line(0,-1){600}}
\put(0,1500){\line(1,0){9000}}
\multiput(1800,1200)(9000,0){2}{\line(0,-1){300}}
\put(1800,1200){\line(1,0){7100}}
\put(9100,1200){\line(1,0){1700}}
\multiput(0,-1350)(10800,0){2}{\line(0,1){450}}
\put(0,-1350){\line(1,0){10800}}
\put(1800,600){\usebox{\tow}}
\put(9000,600){\usebox{\toe}}
\end{picture}\]

Notice the analogies with the tail cases: 
$\s$ is obtained adding $\gamma$ to a spherical system $\s_{S(\gamma)}$ 
(localization of $\s$ in $S(\gamma)=S\setminus\supp\,\gamma$). 
Let us enumerate as above the simple roots of the connected component of $S$ containing $\supp\,\gamma$ in such a way 
that $\gamma=\alpha_{l+1}+\ldots+\alpha_{l+m}$. 
Actually, if $\s$ is indecomposable, without positive $n$-combs with $n>1$, 
the only two simple roots of $S(\gamma)$ non-orthogonal to $\supp\,\gamma$ are $\alpha_l$ and $\alpha_{l+m+1}$. 
As above, $\Delta(\gamma)$ is supposed to be minimal with respect to its properties, 
$K^\gamma\subset G$ and $K_\gamma\subset G_{S(\gamma)}$ are associated with resp.\ $\s/\Delta(\gamma)$ and $\s_{S(\gamma)}$.

The description of the subgroup $H$ of $G$ and the proof that it is actually the wonderful subgroup of $G$ 
associated with $\s$ can be conducted as above. 
Here, take $Q=P_{-S_Q}$ with $S_Q=S^p/\Delta(\gamma)\cup\supp\,\gamma$ 
and notice that $S_Q\not\ni\alpha_l,\alpha_{l+m+1}$. 
Here we have $K^\gamma_{\supp\,\gamma}=\GL(m)$ included in $G_{\supp\,\gamma}=\SL(m+1)$. 
Take $R$ to be the parabolic subgroup of $G_{S(\gamma)}$ containing $B_-\cap G_{S(\gamma)}$ 
and associated to $S(\gamma)\cap S_Q$. 
Choose $K_\gamma\subset R$ (thus $K_\gamma^u\subset R^u$) and $L_{K_\gamma}\subset L_{R}$. 
As above, set $L=\pi^{-1}(L_{K_\gamma})\times K^\gamma_{\supp\,\gamma}\subset L_{S_Q}$ 
and $U\subset Q^u$ such that $\Lie\,R^u/\Lie\,K_\gamma^u$, $\Lie\,Q^u/\Lie\,U$ are isomorphic 
as $\pi^{-1}(L_{K_\gamma})$-modules.  

\bibliographystyle{elsarticle-num}

\end{document}